\documentclass[12pt]{amsart}
\usepackage{latexsym, amssymb, amscd, rotating}
\usepackage{mathrsfs}

\makeatletter
\def\@seccntformat #1{%
    \csname pre#1name\endcsname
    \csname the#1\endcsname
    \csname post#1name\endcsname
    \quad
}
\makeatother

\numberwithin{equation}{section}

\newtheorem{thm}{Theorem}[section]

\newtheorem{defi}[thm]{Definition}
\newtheorem{lem}[thm]{Lemma}
\newtheorem{prop}[thm]{Proposition}
\newtheorem{cor}[thm]{Corollary}
\newtheorem{rem}[thm]{Remark}

\newcommand{\bs}[1]{\boldsymbol{#1}}

\newcommand{\End}{\operatorname{End}}

\newcommand{\Z}{\mathbb{Z}}

\renewcommand{\labelenumi}{(\arabic{enumi})}
\numberwithin{equation}{section}

\newcommand{\Sym}{\mathfrak{S}}

\begin{document}
\setlength{\baselineskip}{4.9mm}
\setlength{\abovedisplayskip}{4.5mm}
\setlength{\belowdisplayskip}{4.5mm}

%%%%%%%%%%%%%%%%%%%%%%%%%%%%%%%%%%%%%%%%%%%%%%%%%%%%%%%%%%%%%%%%%%%%%%%%

\renewcommand{\theenumi}{\roman{enumi}}
\renewcommand{\labelenumi}{(\theenumi)}
\renewcommand{\thefootnote}{\fnsymbol{footnote}}

%%%%%%%%%%%%%%%%%%%%%%%%%%%%%%%%%%%%%%%%%%%%%%%%%%%%%%%%%%%%%%%%%%%%%%%%

\parindent=20pt
\medskip

\begin{center}
\textbf{\Large On decomposition numbers of the cyclotomic $q$-Schur
  algebras}\\  
\vspace{1cm}
Nobuharu Sawada\footnote{The author expresses gratitude to Professor
  Toshiaki Shoji.}\\
\vspace{0.5cm}
Department of Mathematics  \\
Tokyo University of Science \\ 
Noda, Chiba 278-8510, Japan
\end{center}
\title{}
\maketitle

%%%%%%%%%%%%%%%%%%%%%%%%%%%%%%%%%%%%%%%%%%%%%%%%%%%%%%%%%%%%%%%%%%%%%%%%
%%%%%%%%%%%%%%%%%%%%%%%%%%%%%%%%%%%%%%%%%%%%%%%%%%%%%%%%%%%%%%%%%%%%%%%% 

\begin{abstract}
Let $\mathcal{S}(\Lambda)$ be the cyclotomic $q$-Schur algebra associated to
the Ariki-Koike algebra $\mathscr{H}$.   
We construct a certain subalgebra
$\mathcal{S}^{0}(\Lambda)$ of $\mathcal{S}(\Lambda)$, and 
show that it is a standardly based algebra in the sense of Du and
Rui. $\mathcal{S}^{0}(\Lambda)$ has a natural quotient
$\overline{\mathcal{S}^{0}}(\Lambda)$, which turns
out to be a cellular algebra.
In the case where the modified Ariki-Koike algebra $\mathscr{H}^{\flat}$
is defined, $\overline{\mathcal{S}^{0}}(\Lambda)$ coincides with 
the cyclotomic $q$-Schur algebra associated to $\mathscr{H}^{\flat}$.
In this paper, we discuss a relationship among the decomposition 
numbers of $\mathcal{S}(\Lambda)$, $\mathcal{S}^{0}(\Lambda)$ and 
$\overline{\mathcal{S}^{0}}(\Lambda)$.
In particular, we show that some important part of the decomposition 
matrix of $\mathcal{S}(\Lambda)$ coincides with a part of the 
decomposition matrix of $\overline{\mathcal{S}^{0}}(\Lambda)$.
\end{abstract}
\pagestyle{myheadings}
\markboth{NOBUHARU SAWADA}{ON DECOMPOSITION NUMBERS OF THE CYCLOTOMIC
  $q$-SCHUR ALGEBRAS}

%%%%%%%%%%%%%%%%%%%%%%%%%%%%%%%%%%%%%%%%%%%%%%%%%%%%%%%%%%%%%%%%%%%%%%%%
%%%%%%%%%%%%%%%%%%%%%%%%%%%%%%%%%%%%%%%%%%%%%%%%%%%%%%%%%%%%%%%%%%%%%%%% 

\bigskip
\medskip

\addtocounter{section}{-1}

%%%%%%%%%%%%%%%%%%%%%%%%%%%%%%%%%%%%%%%%%%%%%%%%%%%%%%%%%%%%%%%%%%%%%%%%
%%%%%%%%%%%%%%%%%%%%%%%%%%%%%%%%%%%%%%%%%%%%%%%%%%%%%%%%%%%%%%%%%%%%%%%% 

\section{Introduction}

Let $\mathscr{H}$ be the Ariki-Koike algebra over an commutative
integral domain $R$ with parameters $q, q^{-1}, Q_{1}, \ldots , Q_{r} 
\in R$, associated to the complex reflection group $W_{n, r} = G(r,1,n)$. 
In \cite{DJMa} Dipper, James and Mathas introduced the cyclotomic
$q$-Schur algebras $\mathcal{S}(\Lambda)$ with weight poset $\Lambda$
as a tool for studying the representations of the Ariki-Koike algebra 
$\mathscr{H}$.  
It is an important 
problem to determine the decomposition matrix of $\mathcal{S}(\Lambda)$.
In the case where $r =1$, the cyclotomic $q$-Schur algebra
coincides with the $q$-Schur algebra of \cite{DJ}.   
In that case, under the condition that 
$R = \mathbb{C}$ and $q$ is a root of unity, 
Varagnolo and Vasserot  proved in \cite{VV}
the decomposition conjecture due to Leclerc and Thibon \cite{LT}, which
provides us an algorithm of computing the decomposition matrix 
in connection with the canonical basis of the level 1 Fock space of type $A$. 
\par
It is an open problem to determine the decomposition matrix for 
$\mathcal{S}(\Lambda)$ in the case where $r \ge 2$.
It is known by \cite{A} (see also \cite{DM}) that the determination of the
decomposition matrix of $\mathcal{S}(\Lambda)$ is reduced to the case
where $r = 1$ if the parameters satisfy the separation condition
\par\medskip\noindent
(S) \quad $q^{2k}Q_i - Q_j$ are invertible in $R$ for $|k| < n$, 
     $i \ne j$.
\par\medskip\noindent
In \cite{U}, Ugolov constructed the canonical basis of the 
level $r$ Fock space of type $A$, and gave an algorithm of computing them.
Assume that $R = \mathbb C$, and that $Q_i = q^{s_i}$ with 
$s_{i+1} - s_i \ge n$ for a root of unity $q \in \mathbb C$.
In that case, Yvonne \cite{Y} formulated a conjecture 
which makes it possible to determine the decomposition matrix of 
$\mathcal{S}(\Lambda)$ by means of the algorithm of
Ugolov as in the case of $r = 1$.  He showed in \cite{Y} that the
Jantzen sum  
formula for $\mathcal{S}(\Lambda)$ also holds for the Fock space side
in an appropriate sense, which supports the conjecture.
\par
In this paper, we introduce a certain subalgebra 
$\mathcal{S}^0(\Lambda)$ of $\mathcal{S}(\Lambda)$, motivated 
in \cite{SawS}, and show that
$\mathcal{S}^0(\Lambda)$ of $\mathcal{S}(\Lambda)$ is a standardly
based algebra in the sense of Du and Rui \cite{DR1}.
$\mathcal{S}^0(\Lambda)$ has a natural quotient 
$\overline{\mathcal{S}^{0}}(\Lambda)$, which turns out to be 
a cellular algebra.  We discuss a relationship among the 
representations of $\mathcal{S}(\Lambda)$, $\mathcal{S}^0(\Lambda)$ and
$\overline{\mathcal{S}^{0}}(\Lambda)$.  In particular we show, in the
case 
where $R$ is a field, that 
a certain important part of the decomposition matrix of 
$\mathcal{S}(\Lambda)$ coincides with a part of the 
decomposition matrix of $\overline{\mathcal{S}^{0}}(\Lambda)$. 
This reduces a computation of some decomposition numbers of 
$\mathcal{S}(\Lambda)$ to that of 
$\overline{\mathcal{S}^{0}}(\Lambda)$ which seems to have a simpler
structure than $\mathcal{S}(\Lambda)$.
\par
The modified Ariki-Koike algebra
$\mathscr{H}^{\flat}$ over $R$ was discussed in \cite{SawS}, 
which has the same parameter set as  $\mathscr{H}$,  but under the condition that 
\par\medskip\noindent
$(\ast)$ \quad $Q_{i} - Q_{j}$ are invertible in $R$ for any $i\neq j$.  
\par\medskip\noindent
Note that the condition $(\ast)$ is weaker than the separation condition (S).
It is known in [S] that $\mathscr{H}^{\flat}$ is isomorphic to 
$\mathscr{H}$ when $R$ is a field, and the parameters 
satisfies the separation condition. 
Let $\bs{m} = ( m_1, \ldots , m_r )$ be a tuple of non-negative
integers, and 
$\widetilde{\mathcal{P}}_{n,r} = 
\widetilde{\mathcal{P}}_{n,r}(\bs m)$ the set of all
$r$-compositions 
$\lambda = ( \lambda^{( 1 )}, \ldots , \lambda^{( r )} )$
such that $\lambda^{(i)} \in \Bbb Z^{m_i}_{\ge 0}$.
A representation of $\mathscr{H}^{\flat}$
on a tensor space $V^{\otimes n}$ was constructed in [S], where
$V = \bigoplus_{i = 1}^{r} V_{i}$ and $V_{i}$ is a free
$R$-module of rank $m_{i}$ satisfying the condition that 
\par\medskip\noindent
$(\ast\ast)$  \quad $m_{i} \geq n$  for all  $i$.
\par\medskip
In \cite{SawS} we have defined the
cyclotomic $q$-Schur
algebra $\mathcal{S}^{\flat}(\bs{m}, n)$ with the poset 
$\Lambda = \widetilde{\mathcal{P}}_{n,r}$ as the 
endomorphism algebra $\End_{\mathscr{H}^{\flat}}V^{\otimes n}$. 
Moreover, we have constructed a certain subalgebra
$\mathcal{S}^{0}(\bs{m}, n)^{\flat}$ of $\mathcal{S}(\Lambda)$ which has a
surjective algebra homomorphism 
$\widehat{f} : \mathcal{S}^{0}(\bs{m},n)^{\flat} 
        \rightarrow \mathcal{S}^{\flat}(\bs{m}, n)$. 
Note that $\mathcal{S}^{0}(\bs{m},n)^{\flat}$ 
(resp. $\mathcal{S}^{\flat}(\bs{m}, n)$)
coincides with the previous $\mathcal{S}^{0}(\Lambda)$ 
(resp. $\overline{\mathcal{S}^{0}}(\Lambda)$)
whenever the conditions $(\ast)$ and $(\ast\ast)$ are satisfied. 
\par
By a result of \cite{SawS}, the decomposition matrix of
$\mathcal{S}^{\flat}(\bs{m}, n)$ can be described in terms of 
the decomposition 
matrices for $q$-Schur algebras, which are computable by  
the theorem of Varagnolo-Vasserot.     
Thus our results determine a part of the decomposition matrix 
for $\mathcal{S}(\Lambda)$ in the case where $R = \mathbb C$ and 
the parameters $Q_i = q^{s_i}$ are all distinct for a root of unity 
$q \in \mathbb C$. 
\par 
The details of the paper are as follows. In Section 1
we recall the notations and basic results on the representation
theory 
of the Ariki-Koike algebras and the cyclotomic $q$-Schur
algebras. In Section 2 we investigates the subalgebra
$\mathcal{S}^{0}(\Lambda)$ of $\mathcal{S}(\Lambda)$. We show
that $\mathcal{S}^{0}(\Lambda)$ is a
standardly based algebra in the sense of Du and Rui \cite{DR1}, which
is a natural extension of the
cellular algebra. Section 3 discusses
the relationship between $\mathcal{S}^{0}(\Lambda)$ and
$\mathcal{S}(\Lambda)$, and their Weyl modules $Z^{( \lambda, 0)}$ and 
$W^{\lambda}$ for a multipartition (or $r$-partition) $\lambda$ in
$\Lambda$,
respectively. More precisely, we show that $Z^{( \lambda, 0)}
\otimes_{\mathcal{S}^{0}(\Lambda)} \mathcal{S}(\Lambda) \simeq
W^{\lambda}$. Even though $\mathcal{S}^{0}(\Lambda)$ is not cellular,
the Weyl module $Z^{( \lambda, 0)}$ has the radical and
we can discuss almost similar to the cellular theory, i.e, we can
define irreducible modules $L_{0}^{\lambda}$ of
$\mathcal{S}^{0}(\Lambda)$ in a similar way as defining of irreducible
modules 
$L^{\lambda}$ of $\mathcal{S}(\Lambda)$. In turn, in Section 4 we
study the relationship between $\mathcal{S}^{\flat}(\bs{m}, n)$ (or
$\overline{\mathcal{S}^{0}}(\Lambda)$) and
$\mathcal{S}^{0}(\Lambda)$ and their Weyl modules. For a
multipartition $\lambda \in \Lambda$, let
$\overline{Z}^{\lambda}$ be the Weyl module of
$\overline{\mathcal{S}^{0}}(\Lambda)$. Since
$\overline{\mathcal{S}^{0}}(\Lambda)$ is a cellular algebra, we can
also define irreducible modules $\overline{L}^{\lambda}$. Because
$\overline{\mathcal{S}^{0}}(\Lambda)$ is a quotient of
$\mathcal{S}^{0}(\Lambda)$ by $\mathcal{S}^{00}(\Lambda)$, we can
regard $\overline{Z}^{\lambda}$ and $\overline{L}^{\lambda}$ as
$\mathcal{S}^{0}(\Lambda)$-modules. Under this setting, we can verify
that the decomposition number $[ \overline{Z}^{\lambda} :
\overline{L}^{\mu} ]$ is equal to the decomposition number $[
Z^{( \lambda, 0)} : L_{0}^{\mu} ]$ for any multipartitions $\lambda,
\mu$ in $\Lambda$ when $R$ is a field. In section 5, we estimate the
decomposition number $[ W^{\lambda} : Z^{\mu} ]$ by working use of the 
induction of $Z^{( \lambda, 0)}$ and the restriction of
$W^{\lambda}$. Suppose that $R$ is a
field. Our main result asserts that
$[ 
\overline{Z}^{\lambda} : \overline{L}^{\mu} ] =
[ Z^{(\lambda, 0)} : L_{0}^{\mu} ] = [ W^{\lambda} : L^{\mu} ]$ for
all multipartitions $\lambda, \mu$ in $\Lambda$ such that $\lambda =
(\lambda^{( 1 )}, \ldots , \lambda^{( r )} )$, $\mu = (\mu^{( 1 )},
\ldots , 
\mu^{( r )} )$ with $| \lambda^{( i )} | = | \mu^{( i )} |$. Finally,
we denote by $[W^{\lambda^{( i )}} :
L^{\mu^{( i )}} ] ~ ( 1 \le i \le r )$ the decomposition number of
$L^{\mu^{( i )}}$ in $W^{\lambda^{( i )}}$. Again suppose that $R$
is a 
field. Moreover, assume that $\Lambda =
\widetilde{\mathcal{P}}_{n,r}(\bs{m})$ and the previous conditions
$(\ast)$, $(\ast \ast)$. Then, as a corollary, we have that $[
W^{\lambda} : 
L^{\mu} ] = \prod_{i = 1}^{r} [ W^{\lambda^{( i )}} :
L^{\mu^{( i )}} ]$ for all multipartitions $\lambda, \mu$ in
$\widetilde{\mathcal{P}}_{n,r}(\bs{m})$ as above. 
     
%%%%%%%%%%%%%%%%%%%%%%%%%%%%%%%%%%%%%%%%%%%%%%%%%%%%%%%%%%%%%%%%%%%%%%%%
%%%%%%%%%%%%%%%%%%%%%%%%%%%%%%%%%%%%%%%%%%%%%%%%%%%%%%%%%%%%%%%%%%%%%%%%

\par\bigskip\noindent

%%%%%%%%%%%%%%%%%%%%%%%%%%%%%%%%%%%%%%%%%%%%%%%%%%%%%%%%%%%%%%%%%%%%%%%%
%%%%%%%%%%%%%%%%%%%%%%%%%%%%%%%%%%%%%%%%%%%%%%%%%%%%%%%%%%%%%%%%%%%%%%%% 

\section*{Table of contents}

\par\medskip
0. Introduction  
\par
1. Preliminaries on Ariki-Koike algebras and cyclotomic $q$-Schur
   algebras 
\par
2. The standard basis for $\mathcal{S}^{0}(\Lambda)$ 
\par
3. A relationship between $\mathcal{S}^{0}(\Lambda)$ and
   $\mathcal{S}(\Lambda)$ 
\par
4. A relationship between $\mathcal{S}^{\flat}(\bs{m}, n)$ and
   $\mathcal{S}^{0}(\Lambda)$
\par
5. An estimate for decomposition numbers

%%%%%%%%%%%%%%%%%%%%%%%%%%%%%%%%%%%%%%%%%%%%%%%%%%%%%%%%%%%%%%%%%%%%%%%%
%%%%%%%%%%%%%%%%%%%%%%%%%%%%%%%%%%%%%%%%%%%%%%%%%%%%%%%%%%%%%%%%%%%%%%%% 

\par\medskip

%%%%%%%%%%%%%%%%%%%%%%%%%%%%%%%%%%%%%%%%%%%%%%%%%%%%%%%%%%%%%%%%%%%%%%%%
%%%%%%%%%%%%%%%%%%%%%%%%%%%%%%%%%%%%%%%%%%%%%%%%%%%%%%%%%%%%%%%%%%%%%%%%
\section{Preliminaries on Ariki-Koike algebras and Cyclotomic
  $q$-Schur algebras} 

\subsection{}\label{def:the definition of
  Ariki-Koike algebras} 
Fix positive integers $r$ and $n$ and let $\mathfrak{S}_{n}$ be the
symmetric group of degree $n$. Let $R$ be an integral
domain with 1 and $q, Q_1,
\ldots , Q_r$ be elements in $R$,
with invertible $q$. The Ariki-Koike algebra associated to the complex 
reflection group $W_{n, r} = G( r, 1, n )$, is the associative unital
algebra 
$\mathscr{H} = \mathscr{H}_{n, r}$ over $R$ with generators $T_1,
\ldots, T_n$ subject to the following conditions, 
\begin{equation*}
\begin{array}{rll}
( T_{1} - Q_{1} ) \cdots ( T_{1} - Q_{r} ) &= 0, &  \\
( T_{i} -q )( T_{i} + q^{-1} ) &= 0 &( i \geq 2), \\
T_{1} T_{2} T_{1} T_{2} &= T_{2} T_{1} T_{2} T_{1},&  \\
T_{i} T_{j} &= T_{j} T_{i} & ( | i - j | \geq 2), \\
T_{i} T_{i+1} T_{i} &=T_{i+1} T_{i} T_{i+1} & ( 2 \le i \le n-1 ).
\end{array}
\end{equation*} 
It is known that $\mathscr{H}$ is a free $R$-module of rank
$n!r^{n}$. The subalgebra $\mathscr{H}( \mathfrak{S}_{n} )$ of
$\mathscr{H}$ generated by $T_{2}, \ldots, T_{n}$ is isomorphic to the
Iwahori-Hecke algebra $\mathscr{H}_{n}$ of the symmetric group
$\mathfrak{S}_{n}$. 

For $i = 2, \ldots , n$ let $s_i$ be the transposition $( i-1, i )$ in
$\mathfrak{S}_{n}$.Then $\{ s_2, \ldots, s_n \}$ generate
$\mathfrak{S}_{n}$. For $w \in \mathfrak{S}_{n}$, we set $T_{w} =
T_{{i}_{1}} \cdots T_{{i}_{k}}$ where $w = s_{i_1} \cdots s_{i_k}$ is
a reduced expression. Then
$T_{w}$ is independent of the choice of a reduced expression. We also
put $L_{k} = T_{k} \cdots T_{2}T_{1}T_{2} \cdots T_{k}$ for $k = 1, 2,
\ldots , n$. Note that all $L_1, \ldots , L_n$ commutes. Moreover,
these 
elements produce a basis of $\mathscr{H}$.

\addtocounter{thm}{1}
\begin{thm}[{\cite[Theorem 3.10]{AK}}]
The Ariki-Koike algebra $\mathscr{H}$ is free as an $R$-module with
basis $\{ L_{1}^{a_{1}} \cdots L_{n}^{a_{n}} T_{w} \mid w \in
\mathfrak{S}_{n},~ 0 \le a_{i} < r \text{ for } 1 \le i \le n \}$.
\end{thm}

Recall that a composition of $n$ is sequence $\sigma = ( \sigma_{1},
\sigma_{2}, \ldots )$ of non-negative integers such that $| \sigma | =
\sum_{i} \sigma_{i} = n$. $\sigma$ is a partition if in addition
$\sigma_{1} \geq \sigma_{2} \geq \cdots $. If $\sigma_{i} = 0$ for all
$i > k$ then we write $\sigma = ( \sigma_{1}, \ldots , \sigma_{k} )$.

An $r$-composition (or multicomposition) of $n$ is an $r$-tuple
$\lambda = ( \lambda^{( 1 )}, \ldots, \lambda^{( r )} )$ of
compositions with $\lambda^{( i )} = ( \lambda^{( i )}_1, \lambda^{( i
  )}_2, \ldots )$ such that $| \lambda^{( 1 )}
| + \cdots + | \lambda^{( r )} | = n$. An $r$-composition $\lambda$ is
an $r$-partition if each $\lambda^{( i )}$ is a partition. If
$\lambda$ is an $r$-partition of $n$ then we write $\lambda \vdash
n$. The diagram $[ \lambda ]$ of the $r$-composition $\lambda$ is the
set $[ \lambda ] = \{ ( i, j, s ) \mid 1 \le i \le \lambda_{j}^{( s )},
 1 \le s \le r \}$. The elements of $[ \lambda ]$ are called nodes. 
The set of $r$-compositions of $n$ is partially ordered by dominance,
i.e, if $\lambda$ and $\mu$ are two $r$-compositions then
$\lambda$ dominates $\mu$, and we write $\lambda \unrhd \mu$, if
\begin{equation*}
\underset{c = 1}{\overset{s-1}{\sum}} | \lambda^{( c )}| + \underset{j
  = 1}{\overset{i}{\sum}} | \lambda_{j}^{( s )}| \geq \underset{c =
  1}{\overset{s-1}{\sum}} | \mu^{( c )}| + \underset{j
  = 1}{\overset{i}{\sum}} | \mu_{j}^{( s )}|
\end{equation*}  
for $1 \le s \le r$ and for all $i \geq 1$. If $\lambda \unrhd \mu$
and $\lambda \neq \mu$ then we write $\lambda \rhd \mu$.

If $\lambda$ is an $r$-composition let $\mathfrak{S}_{\lambda}
=\mathfrak{S}_{\lambda^{( 1 )}} \times \cdots \times
\mathfrak{S}_{\lambda^{( r )}}$ be the corresponding Young subgroup of
$\mathfrak{S}_{n}$. Set
\begin{equation*}
x_{\lambda} = \underset{w \in \mathfrak{S}_{\lambda}}{\sum} q^{l( w )}
T_{w}, \qquad u^{+}_{\lambda} = \underset{s = 2}{\overset{r}{\prod}}
\underset{k = 1}{\overset{a_{s}}{\prod}} ( L_{k} - Q_{s} ),  
\end{equation*}  
where $a_s = | \lambda^{( 1 )} | + \cdots + | \lambda^{( s-1 )} |$ for
$2 \le s \le r$. If $s = 1$ then we set $a_{s} = 0$. Set $m_{\lambda} =
x_{\lambda} u^{+}_{\lambda} = u^{+}_{\lambda} x_{\lambda}$ and define
$M^{\lambda}$ to be the right ideal $M^{\lambda} = m_{\lambda}
\mathscr{H}$ of $\mathscr{H}$.

For any $r$-composition $\mu$, a $\mu$-tableau $\mathfrak{t} =(
\mathfrak{t}^{( 1 )}, \ldots , \mathfrak{t}^{( r )} )$ is a bijection
$t : [ \mu ] \rightarrow \{ 1, 2, \ldots , n \}$, where
$\mathfrak{t}^{( i )}$ is a tableau of $\text{Shape}( \mathfrak{t}^{(
  i )} ) = \mu^{( i )}$. We write $\text{Shape}(
  \mathfrak{t} ) = \mu $ if $\mathfrak{t}$ is a $\mu$-tableau. A
  $\mu$-tableau $\mathfrak{t}$ is called standard
  (resp. row standard) if all $\mathfrak{t}^{( i )}$ are standard
  (resp. row standard). 
Let 
$\text{Std}(\lambda)$ be the set of standard $\lambda$-tableaux. 

For each $r$-composition $\mu$, let $\mathfrak{t}^{\mu}$ be the
$\mu$-tableau with 
the numbers $1, 2, \ldots , n$ attached in order from left to right
along its rows and from top to bottom, and from $\mu^{( 1 )}$ to
$\mu^{( r )}$. If $\mathfrak{t}$ is any row standard $\mu$-tableau
let $d( \mathfrak{t} ) \in \mathfrak{S}_{n}$ be the unique permutation
such that $\mathfrak{t} = \mathfrak{t}^{\mu} d( \mathfrak{t}
)$. Furthermore, let $* : \mathscr{H} \rightarrow \mathscr{H}$ be the
anti-isomorphism given by $T_{i}^{*} = T_{i}$ for $i = 1, 2, \ldots ,
n$, and set $m_{\mathfrak{s} \mathfrak{t}} = T_{d( \mathfrak{s} )}^{*}
m_{\lambda} T_{d( \mathfrak{t} )}$. 
\begin{thm}[{\cite[Theorem 3.26; Murphy bases]{DJMa}}]
The Ariki-Koike algebra $\mathscr{H}$ is free as an $R$-module with
cellular basis $\{ m_{\mathfrak{s} \mathfrak{t}} \mid \mathfrak{s},
\mathfrak{t} \in {\rm{Std}}(\lambda) \text{ for some } \lambda
\vdash n \}$.
\end{thm}
Here, and below, whenever we write expressions involving a pair of
tableaux (such as $m_{\mathfrak{s} \mathfrak{t}}$ or $\varphi_{ST}$
and so on), we implicitly assume that the two tableaux are of the same
shape. 

\addtocounter{subsection}{2}
\subsection{}\label{sub:the definition of
  H^lambda} 
For each $r$-partition $\lambda$ let $\mathscr{H}^{\lambda}$ be the
$R$-submodule of $\mathscr{H}$ with basis $\{ m_{\mathfrak{u}
  \mathfrak{b}} \mid \mathfrak{u}, \mathfrak{b} \in
\text{Std}(\mu) \text{ for some } \mu \rhd \lambda \}$. Then
$\mathscr{H}^{\lambda}$ is a two-sided ideal of $\mathscr{H}$.

Let $S^{\lambda}$ be the Specht module (or cell module)
corresponding to the $r$-partition $\lambda$. Thus, $S^{\lambda}
= ( m_{\lambda} + \mathscr{H}^{\lambda} ) \mathscr{H}$, a submodule of
$\mathscr{H} / \mathscr{H}^{\lambda}$. For each $\mathfrak{t} \in
\text{Std}(\lambda)$ let $m_{\mathfrak{t}} = m_{\mathfrak{t}^{\lambda}
  \mathfrak{t}} + \mathscr{H}^{\lambda}$. Then $S^{\lambda}$ is free
as an $R$-module with basis $\{ m_{\mathfrak{t}} \mid \mathfrak{t} \in
\text{Std}(\lambda) \}$. 

Furthermore, there is an associative symmetric bilinear form on
$S^{\lambda}$ which is determined by 
\begin{equation*}
m_{\mathfrak{t}^{\lambda} \mathfrak{s}} m_{\mathfrak{t}
  \mathfrak{t}^{\lambda}}  \equiv
\langle m_{\mathfrak{s}}, m_{\mathfrak{t}} \rangle
m_{\mathfrak{t}^{\lambda} \mathfrak{t}^{\lambda}} \mod
  \mathscr{H}^{\lambda} 
\end{equation*}
for all $\mathfrak{s}, \mathfrak{t} \in
\text{Std}(\lambda)$. The radical $\text{rad}{S^{\lambda}}$ of this
form is again an $\mathscr{H}$-module, so $D^{\lambda} = S^{\lambda} / 
\text{rad}{S^{\lambda}}$ is an $\mathscr{H}$-module. When $R$ is a
field, $D^{\lambda}$ is either $0$ or absolutely irreducible and all
the simple $\mathscr{H}$-modules up to isomorphisms arise uniquely in
this way. 

We can now give a definition of the cyclotomic $q$-Schur algebras. A
set $\Lambda$ of $r$-compositions of $n$ is saturated if $\Lambda$ is
finite and whenever $\lambda$ is an $r$-partition such that $\lambda
\unrhd \mu$ for some $\mu \in \Lambda$ then $\lambda \in \Lambda$. If
$\Lambda$ is a saturated set of $r$-compositions, we denote by
$\Lambda^{+}$ be  
the set of $r$-partitions in $\Lambda$.

\addtocounter{thm}{1}
\begin{defi}
Suppose that $\Lambda$ is a saturated set of multicompositions of
$n$. The cyclotomic $q$-Schur algebra with weight poset $\Lambda$ is
the endomorphism algebra
\begin{equation*}
\mathcal{S}( \Lambda ) = \End_{\mathscr{H}}( M( \Lambda ) ), \qquad
\text{ where } M( \Lambda ) = \underset{ \lambda \in
  \Lambda}{\bigoplus} M^{\lambda}. 
\end{equation*}
\end{defi}

Let $\lambda$ be an $r$-partition and $\mu$ an $r$-composition. A
$\lambda$-Tableau of type $\mu$ is a map $T : [ \lambda ] \rightarrow
\{ ( i, s ) \mid i \geq 1, ~ 1 \le s \le r \}$ such that $\mu^{( s
  )}_{i} = \sharp \{ x \in [ \lambda ] \mid T( x ) = ( i, s ) \}$ for
all 
$i \geq 1$ and $1 \le s \le r$. We regard $T$ as an $r$-tuple
$T = ( T^{( 1 )}, \ldots , T^{( r )} )$, where $T^{( s )}$ is the
$\lambda^{( s )}$-tableau with $T^{( s )}(i, j) = T(i, j, s)$ for all 
$(i, j, s) \in [ \lambda ]$. In this way we identify the standard
tableaux above with the Tableaux of type $w = ( ( 0 ), \ldots , (
1^{n} ) )$. If $T$ is a Tableau of type $\mu$ then we write
$\text{Type}( T ) 
= \mu$.

Given two pairs $( i, s )$ and $( j, t )$ write $( i, s ) \preceq ( j,
t )$ if either $s < t$, or $s = t$ and $i \le j$. 

\begin{defi}
A Tableau $T$ is $($row$)$ semistandard if, for $1 \le t \le r$, the
entries in $T^{( t )}$ are 

$({\rm{i}})$ weakly increasing along the rows with respect to
$\preceq$, 

$({\rm{ii}})$ strictly increasing down columns,

$({\rm{iii}})$ $( i, s )$ appears in $T^{( t )}$ only if $ s \geq
t$. 
\end{defi}

Let $\mathcal{T}_0(\lambda, \mu)$ be the set of semistandard
$\lambda$-Tableaux of type $\mu$ and let
$\mathcal{T}_0(\lambda) = \mathcal{T}^{\Lambda}_0(\lambda) =
\bigcup_{\mu \in \Lambda}
\mathcal{T}_0(\lambda, \mu)$.

Notice that if $\mathcal{T}_0(\lambda, \mu)$ is non-empty, then
$\lambda \unrhd \mu$.

Suppose that $\mathfrak{t}$ is a standard $\lambda$-tableau and let
$\mu$ be an $r$-composition. Let $\mu( \mathfrak{t} )$ be the Tableau
obtained from $\mathfrak{t}$ by replacing each entry $j$ with $( i, k
)$ if $j$ appears in row $i$ of $( \mathfrak{t}^{\mu} )^{( k )}$. The 
tableau $\mu( \mathfrak{t} )$ is a $\lambda$-Tableau of type $\mu$. It
is not necessarily semistandard. 

If $S$ and $T$ are semistandard $\lambda$-Tableaux of type $\mu$ and
$\nu$ respectively,                                                                 
let

\begin{equation*}
\qquad m_{ST}
= 
\underset{ 
\substack{
\mathfrak{s}, \mathfrak{t}\in \text{Std}(\lambda) \\
\mu( \mathfrak{s} ) = S, ~ \nu( \mathfrak{t} ) = T
}
}{\sum}
q^{l( d(\mathfrak{s}) ) + l( d(\mathfrak{t}) )} m_{\mathfrak{s}
  \mathfrak{t}}. 
\end{equation*}
 
 For $S$ and $T$ as above we define a map $\varphi_{ST}$ on
 $M(\Lambda)$ 
 by $\varphi_{ST}( m_{\alpha} h ) =
 \delta_{\alpha \nu} m_{ST} h$, for all $h \in \mathscr{H}$ and all
 $\alpha \in \Lambda$. Here $\delta_{\alpha \nu}$ is the Kronecker
 delta, i.e, $\delta_{\alpha \nu} = 1$ if $\alpha = \nu$ and it is
 zero otherwise. Then $\varphi_{ST}$ is well-defined, and it belongs
 to $\mathcal{S}( \Lambda
 )$. Moreover, 

\begin{thm}[{\cite[Theorem 6.6]{DJMa}}]
The cyclotomic $q$-Schur algebra $\mathcal{S}( \Lambda )$ is free as
an $R$-module with cellular basis $\mathcal{C}( \Lambda ) = \{
\varphi_{ST} \mid S, T \in \mathcal{T}^{\Lambda}_0(\lambda) \text{ for
  some } \lambda \in \Lambda^{+} \}$. 
\end{thm} 

The basis $\{ \varphi_{ST} \}$ is called a semistandard basis of
$\mathcal{S}( \Lambda ) $. Since this basis is cellular, the map $* :
\mathcal{S}( \Lambda ) \rightarrow \mathcal{S}( \Lambda )$ which is
determined by $\varphi_{ST}^{*} = \varphi_{TS}$ is an
anti-automorphism 
of $\mathcal{S}( \Lambda )$. This involution is closely related to the
$*$-involution on $\mathscr{H}$. Explicitly, if $\varphi : M^{\nu}
\rightarrow M^{\mu}$ is an $\mathscr{H}$-module homomorphism then
$\varphi^{*} : M^{\mu} \rightarrow M^{\nu}$ is the homomorphism given
by $\varphi^{*}( m_{\mu} h ) = ( \varphi( m_{\nu} ) )^{*} h$, for all
$h \in \mathscr{H}$.

For each $r$-partition $\lambda \in \Lambda^{+}$, we define
$\mathcal{S}^{\vee \lambda} = \mathcal{S}^{\vee}( \Lambda )^{\lambda}$
as the $R$-span of $\varphi_{ST}$ such that $S, T \in
\mathcal{T}^{\Lambda}_0(\alpha)$ with $\alpha \rhd \lambda$, which is
a two-sided ideal of
$\mathcal{S}( \Lambda )$. We define the Weyl module $W^{\lambda}$ by
the right $\mathcal{S}( \Lambda )$-submodule of $\mathcal{S}( \Lambda
) / \mathcal{S}^{\vee}( \Lambda )^{\lambda}$ generated by the image
$\varphi_{\lambda} = \varphi_{T^{\lambda} T^{\lambda}} \in
\mathcal{S}( \Lambda )$  where $T^{\lambda} = \lambda(
\mathfrak{t}^{\lambda} )$. For each $T \in
\mathcal{T}^{\Lambda}_0(\lambda)$, let 
$\varphi_{T}$ be the image of $\varphi_{T^{\lambda} T}$ in
$W^{\lambda}$. Then the Weyl module $W^{\lambda}$ is $R$-free with
basis $\{ \varphi_{T} \mid T \in \mathcal{T}^{\Lambda}_0(\lambda)
\}$. As in the case of Specht modules there is an inner product on
$W^{\lambda}$ which is determined by 
\begin{equation*}
\varphi_{T^{\lambda}
  S} \varphi_{T T^{\lambda}} \equiv \langle \varphi_{S}, \varphi_{T}
  \rangle \varphi_{T^{\lambda} T^{\lambda}} \mod \mathcal{S}^{\vee
  \lambda}.
\end{equation*}

Let $\text{rad}{W^{\lambda}} = \{ x \in W^{\lambda} \mid \langle x, y
\rangle = 0 \text{ for all } y \in W^{\lambda} \}$. The quotient
module $L^{\lambda} = W^{\lambda} /
\text{rad}{W^{\lambda}}$ is absolutely irreducible and $\{ L^{\lambda}
\mid \lambda \in \Lambda^{+} \}$ is a complete set of non-isomorphic
irreducible $\mathcal{S}( \Lambda )$-modules.

\addtocounter{subsection}{3}
\subsection{}
For an $r$-composition $\mu$, we define the type $\alpha = \alpha( \mu
)$ of $\mu$ by $\alpha = ( n_1, \ldots , n_r )$ with $n_{i} = | \mu^{(
  i )}|$, and the size of $\mu$ by $n = \sum_{i = 1}^{r} n_{i}$. We
also define a sequence $\bold{a} = \bold{a}( \mu ) = ( a_1, \ldots ,
a_r 
)$. (Recall that $a_{i} = \sum_{k = 1}^{i - 1} | \mu^{( k )} | =
\sum_{k = 1}^{i - 1} n_{k} $.)

 We define a partial order $\geq$ on the set $\mathbb{Z}^{r}_{\geq 0}$
 by 
$\bold{a} \geq {\bold{a}}'$ for $\bold{a} = ( a_1, \ldots , a_r )$, 
${\bold{a}}' = ( {a}'_1, \ldots , {a}'_r ) \in \mathbb{Z}^{r}_{\geq
  0}$ 
if $a_{i} \geq {a}'_{i}$ for any $i$. We write $\bold{a} >
{\bold{a}}'$ if 
$\bold{a} \geq {\bold{a}}'$ and $\bold{a} \neq {\bold{a}}'$. It is
clear that
\begin{equation}\label{eq:A certain fact of r-compositions} 
\text{If } \lambda \unrhd \mu, \text{ then } \bold{a}( \lambda ) \geq 
\bold{a}( \mu ) \text{ for } r\text{-compositions } \lambda, \mu.
\end{equation}
Hence if $\mathcal{T}_0(\lambda, \mu)$ is non-empty, then
$\lambda 
\unrhd \mu$, and so we have $\bold{a}( \lambda ) \geq \bold{a}( \mu
)$. 

For any $r$-partition $\lambda$ and $r$-composition $\mu$, we define a 
subset $\mathcal{T}^{+}_0(\lambda, \mu)$ of $\mathcal{T}_0(\lambda,
\mu)$ by
\begin{equation*}
\mathcal{T}^{+}_0(\lambda, \mu) = \{ S \in \mathcal{T}_0(\lambda,
\mu) \mid \bold{a}( \lambda ) = \bold{a}( \mu ) \}.
\end{equation*} 
 Note that the condition $\bold{a}( \lambda ) = \bold{a}( \mu )$ is
 equivalent to $\alpha( \lambda ) = \alpha( \mu )$. Take $S \in
 \mathcal{T}^{+}_0(\lambda, \mu)$. Then one can check that $S \in
 \mathcal{T}^{+}_0(\lambda, \mu)$ if and only if each entry of $S^{( k
 )}$ is of the form $( i, k )$ for some $i$. Hence in this case $S^{(
 k )}$ can be identified with a semistandard $\lambda^{( k )}$-Tableau
 of type $\mu^{( k )}$ under the usual definition of the semistandard
 Tableaux for $1$-partitions $\lambda^{( k )}$ and $1$-compositions
 $\mu^{( k )}$. It follows that we have a bijection
\begin{equation*}
\mathcal{T}^{+}_0(\lambda, \mu) \simeq \mathcal{T}_0(\lambda^{( 1 )},
\mu^{( 1 )} ) \times \cdots \times \mathcal{T}_0(\lambda^{( r )}, 
\mu^{( r )} )
\end{equation*}
via $S \leftrightarrow ( S^{( 1 )}, \ldots , S^{( r )} )$. Moreover,
if $\mathfrak{s} \in \text{Std}( \lambda )$ is such that $\mu(
\mathfrak{s} ) = S$ with $S \in \mathcal{T}^{+}_0(\lambda, \mu)$, then
the entries of $i$-th component of $\mathfrak{s}$ consist of numbers
$a_{i} + 1, \ldots , a_{i + 1}$ for $\bold{a}( \lambda ) = ( a_{1},
\ldots , a_{r} )$. In particular, $d( \mathfrak{s} ) \in
\mathfrak{S}_{\alpha}$ for $\alpha = \alpha( \lambda )$. 

Fix an $r$-tuple $\bs{m} = (
m_1, \ldots , m_r )$ of non-negative integers.
Then, an $r$-composition $\mu = ( \mu^{( 1 )}, \ldots , \mu^{( r )} )$ 
with $\mu^{( i )} =( \mu^{( i )}_{1}, \ldots , \mu^{( i )}_{m_{i}} )
\in 
\mathbb{Z}_{\geq 0}^{m_{i}}$ is called an $(r,
\bs{m})$-composition, and $(r, \bs{m})$-partition is defined
similarly. We denote by  
$\widetilde{\mathcal{P}}_{n,r} = \widetilde{\mathcal{P}}_{n,r}( \bs{m}
)$ (resp. $\mathcal{P}_{n,r} =\mathcal{P}_{n,r}( \bs{m} )$) the set
of $(r, \bs{m})$-compositions (resp. $(r, \bs{m})$-partitions) of
size $n$. (Note that
$\mathcal{P}_{n,r}( \bs{m} )$ are naturally identified with each other
for any $\bs{m}$ such that $m_{i} \geq n$. However,
$\widetilde{\mathcal{P}}_{n,r}$ depends on the choice of $\bs{m}$.)
Finally, let 
\begin{equation*}
\begin{array}{ll}
\mathcal{C}^{0}(\Lambda) = \underset{
\mu, \nu \in \Lambda, ~ \lambda \in \Lambda^{+} 
}{\bigcup}
\{ \varphi_{ST} \in \mathcal{C}(\Lambda) 
\mid 
&S \in \mathcal{T}_0(\lambda, \mu), ~ T \in \mathcal{T}_0(\lambda,
\nu),\\ 
&\bold{a}( \lambda ) > \bold{a}( \mu ) \text{ if } \alpha( \mu ) \neq
\alpha( \nu ) \}
\end{array}
\end{equation*}
and we define $\mathcal{S}^{0}(\Lambda)$ as the $R$-submodule of 
$\mathcal{S}(\Lambda)$ with basis
$\mathcal{C}^{0}(\Lambda)$.   

\medskip

%%%%%%%%%%%%%%%%%%%%%%%%%%%%%%%%%%%%%%%%%%%%%%%%%%%%%%%%%%%%%%%%%%%%%%%%
%%%%%%%%%%%%%%%%%%%%%%%%%%%%%%%%%%%%%%%%%%%%%%%%%%%%%%%%%%%%%%%%%%%%%%%%

\section{The standard basis for $\mathcal{S}^{0}(\Lambda)$}\label{The
  standard basis for S^0(Lambda)}

\subsection{}
In \cite{SawS}, the $R$-submodule
$\mathcal{S}^{0}(\Lambda)$ was constructed as an $R$-subalgebra of
$\mathcal{S}(\Lambda)$ when $\Lambda = \widetilde{\mathcal{P}}_{n,r}$
subject to the conditions $(\ast)$, $(\ast \ast)$ in Introduction. 
But, as a fact, $\mathcal{S}^{0}(\Lambda)$ becomes an $R$-subalgebra
with no condition and it become
a standardly based algebra in the sense of \cite{DR1} which is an 
extension of the cellular algebra. The standardly based algebra has a
basis, namely 
standard basis, corresponding to the cellular basis in the case of the  
cellular algebra. The aim of this section is to prove that
$\mathcal{S}^{0}(\Lambda)$ is an $R$-subalgebra of
$\mathcal{S}(\Lambda)$ and $\mathcal{C}^{0}(\Lambda)$ is a standard
basis for 
$\mathcal{S}^{0}(\Lambda)$.  

First, we recall the definition of standardly based algebras in
\cite{DR1}. 

\addtocounter{thm}{1}
\begin{defi}[{\cite{DR1}}]\label{def:standardly based}
Assume that $R$ is a commutative ring with 1. Let $A$ be an
$R$-algebra and $(\Lambda,\ge)$ a poset. $A$ is called a
\emph{standardly 
based algebra on} $\Lambda$ $($or \emph{standardly based}$)$ if the
following conditions hold.

$(a)$ For any $\lambda \in \Lambda$, there are index sets
$I(\lambda)$,$J(\lambda)$ and subsets

\begin{equation*}
\mathscr{A}^{\lambda}=\{a^{\lambda}_{i,j} \mid (i,j) \in I(\lambda)
\times J(\lambda) \}
\end{equation*}

of $A$ such that the union $\mathscr{A} = \underset{{\lambda \in
    \Lambda}}{\bigcup} \mathscr{A}^{\lambda}$ is disjoint and forms an 
    $R$-basis for $A$.

$(b)$ For any $a \in A$, $a^{\lambda}_{i,j} \in \mathscr{A}$, we
have 

\begin{equation*}
\begin{array}{c}
a\cdot a^{\lambda}_{i,j}\equiv \underset{i' \in
  I(\lambda)}{\sum}f_{i',\lambda}(a,i)~a^{\lambda}_{i',j}\mod
  A(>\lambda)\\  
a^{\lambda}_{i,j} \cdot a  \equiv \underset{j' \in J(\lambda)}{\sum}
f_{\lambda,j'}(j,a) ~  a^{\lambda}_{i,j'}  \mod A(> \lambda),
\end{array}
\end{equation*}

where $A( > \lambda )$ is the $R$-submodule of $A$ spanned by
$\mathscr{A}^{\mu}$ with $\mu > \lambda$, and $f_{i',\lambda}(a,i)
$, $f_{\lambda,j'}(j,a)  \in R$ are
independent of $j$ and $i$, respectively. Such a basis $\mathscr{A}$ is
called a \emph{standard basis} for the algebra $A$.

\end{defi}

Note that the cellular algebra is a special case of the standardly
based algebras.
We prepare some notation.\\
Let
\begin{equation*}
\varOmega= ( \Lambda^{+} \times \{0, 1\} ) \setminus \{ (\lambda, 1)
     \mid 
     \mathcal{T}_0(\lambda,\mu) = \emptyset 
     \text{ for any } \mu \in \Lambda \text{ such that } 
      \mathbf a(\lambda) > \mathbf a(\mu)\}
\end{equation*}
and we define a partial order
$(\lambda_1,\varepsilon_1)\geq(\lambda_2,\varepsilon_2)$ on
$\varOmega$ by 
$(\lambda_1,\varepsilon_1)>(\lambda_2,\varepsilon_2)$ if
$\lambda_1\rhd\lambda_2$, or $\lambda_1=\lambda_2$ and
$\varepsilon_1>\varepsilon_2$. 
For a $( \lambda, \varepsilon ) \in \varOmega$, we define index sets 
$I(\lambda,\varepsilon)$, $J(\lambda,\varepsilon)$ by 
\begin{equation*}
\begin{array}{ll}

I(\lambda,\varepsilon)= 

\left\{\begin{array}{ll}
\mathcal{T}_0^+(\lambda) & \qquad \text{ if } \varepsilon=0, \\[3mm]
\underset{
\mu \in \Lambda, ~
     \bold{a}(\lambda)>
     \bold{a}(\mu)
}{\bigcup}
\mathcal{T}_0(\lambda,\mu) & \qquad \text{ if } \varepsilon=1,
\end{array}\right.

\\[10mm]

J(\lambda,\varepsilon)=

\left\{\begin{array}{ll}
\mathcal{T}_0^+(\lambda) & \qquad \text{ if } \varepsilon=0, \\[3mm]
\mathcal{T}_0(\lambda) & \qquad \text{ if } \varepsilon=1,
\end{array}\right.

\end{array}
\end{equation*}
where $\mathcal{T}_0^+(\lambda) = \bigcup_{\mu \in \Lambda}
\mathcal{T}_0^{+}(\lambda, \mu )$. Then $I(\lambda,\varepsilon)$ and
$J(\lambda, \varepsilon)$ are not
empty for all $( \lambda, \varepsilon ) \in \varOmega$. Assume that
$( \lambda, \varepsilon ) \in \varOmega$. We define a subset $
\mathcal{C}^0(\lambda,\varepsilon)$ of $ \mathcal{S}^0(\Lambda)$ by 

\begin{equation*}
\mathcal{C}^0(\lambda,\varepsilon)=\{ \varphi_{ST} \mid
(S,T) \in I(\lambda,\varepsilon)\times J(\lambda,\varepsilon) \}.
\end{equation*}
It is easy to see that 
\begin{equation}\label{eq:the property of C^0(m,n)}
\text{the union} \underset{{(\lambda,\varepsilon) \in
    \varOmega}}{\bigcup} \mathcal{C}^0(\lambda,\varepsilon) ~ \text{is    
    disjoint and is equal to the set} ~ \mathcal{C}^0(\Lambda).\\
\end{equation}

First, we need the following Lemma.
\begin{lem}[{\cite{M2}}]\label{lem:multiplication formula}
Suppose that $\lambda_i \in \Lambda^{+}$ and $\mu_i, \nu_i \in
\Lambda$ $(i=1,2)$. Assume that
$\varphi_{S_1T_1}$, $\varphi_{S_2T_2} \in
\mathcal{C}(\Lambda)$ where $S_i \in
\mathcal{T}_0(\lambda_i,\mu_i)$, $T_i \in
\mathcal{T}_0(\lambda_i,\nu_i)$. Then
\begin{equation*}
\varphi_{S_1T_1} \cdot \varphi_{S_2T_2} = \delta_{\nu_1
  \mu_2} \cdot \underset{\substack{
\lambda \in \Lambda^{+}, ~ \lambda \unrhd \lambda_1 \text{and}
\lambda_2, \\
S \in \mathcal{T}_0(\lambda,\mu_1), ~ T
  \in\mathcal{T}_0(\lambda,\nu_2)}}{\sum} r_{ST} \cdot \varphi_{ST}
\end{equation*}
where $r_{ST} \in R$, and $\delta_{\nu_1 \mu_2}$ is such that
  $\delta_{\nu_1 
  \mu_2}=1$ if $\nu_1=\mu_2$ and it is zero otherwise.
\end{lem}
\begin{proof}
This lemma was shown by Mathas \cite[(2.8)]{M2}. However, since this
fact 
itself is important for later
discussions, we give the proof here.

It is sufficient to consider the case where $\nu_1 = \mu_2$. For all
$\mu \in \Lambda$, by definition of
$\varphi_{ST}$, we may suppose that
$\varphi_{S_2T_2}(m_{\mu}) = \delta_{\mu
  \nu_2}m_{S_2T_2}= \delta_{\mu \nu_2}
m_{\mu_2}h$ with some $h \in \mathscr{H}$. Since
$m_{S_iT_i} \in M^{\nu_i \ast} \cap
M^{\mu_i}$ $(i=1,2)$, we have
$(\varphi_{S_1T_1} \cdot
\varphi_{S_2T_2})(m_{\mu}) = \delta_{\mu
  \nu_2}m_{S_1T_1}h \in M^{\nu_2 \ast} \cap
M^{\mu_1}$. By \cite[Proposition 6.3]{DJMa}, we deduce that 
\begin{equation*}
m_{S_1T_1}h = \underset{\substack{
\lambda \in \Lambda^{+},\\
S \in \mathcal{T}_0(\lambda,\mu_1), ~ T \in
\mathcal{T}_0(\lambda,\nu_2)}}{\sum} r_{ST} \cdot m_{ST}
\end{equation*}
where $r_{ST} \in R$. Therefore, 
\begin{equation*}
\varphi_{S_1T_1} \cdot \varphi_{S_2T_2}
=\underset{\substack{\lambda \in \Lambda^{+},\\
S \in \mathcal{T}_0(\lambda,\mu_1), ~ T \in
\mathcal{T}_0(\lambda,\nu_2)}}{\sum} r_{ST} \cdot \varphi_{ST}.
\end{equation*}
Note that the set $\{\varphi_{ST} \}$ is the cellular basis
for $\mathcal{S}(\Lambda)$ by \cite[Theorem 6.6]{DJMa}. Then by
the property of 
  cellular basis, the last sum is over $\lambda \in \Lambda^{+}$
  with $\lambda \unrhd \lambda_1$ and $\lambda \unrhd \lambda_2$, $S
  \in 
  \mathcal{T}_0(\lambda,\mu_1)$ and $T \in
  \mathcal{T}_0(\lambda,\nu_2)$. This proves the lemma.
\end{proof}

For all $\lambda \in \Lambda^{+}$, $\mu, \nu \in
\Lambda$ with $\mathcal{T}_0(\lambda, \mu),
\mathcal{T}_0(\lambda, \nu )\neq
\emptyset$ and every $S \in \mathcal{T}_0(\lambda, \mu)$, $T \in
\mathcal{T}_0(\lambda, \nu)$, since definitions of $\varphi_{ST}$ and 
$m_{ST}$ 
we 
find that
\begin{equation*}
\begin{array}{ll}
(\varphi_{ST^{\lambda}}
\varphi_{T^{\lambda}T})(m_{{\mu}'}) 
& = {\delta}_{{\mu}' \nu} \cdot
\varphi_{ST^{\lambda}}(m_{T^{\lambda} T}) 
= {\delta}_{{\mu}' \nu} \cdot
\varphi_{ST^{\lambda}}(m_{\lambda}) 
\underset{
\substack{
\mathfrak{t} \in \text{Std}(\lambda)\\
{\nu}(\mathfrak{t}) = T
} 
}{\sum} T_{d(\mathfrak{t})} \\
&= {\delta}_{{\mu}' \nu} \cdot m_{S T^{\lambda}}\underset{
\substack{
\mathfrak{t} \in \text{Std}(\lambda)\\
{\nu}(\mathfrak{t}) = T
}
}
{\sum} T_{d(\mathfrak{t})}
= {\delta}_{{\mu}' \nu} \cdot m_{ST} 
\end{array}
\end{equation*}
where ${\mu}' \in \Lambda$. Therefore,
\begin{equation}\label{eq:decomposition of varphiST}
\varphi_{ST^{\lambda}} \varphi_{T^{\lambda}T} =
\varphi_{ST} \quad
(\text{for } 
^{\forall}S \in \mathcal{T}_0(\lambda, \mu),
^{\forall}T \in \mathcal{T}_0(\lambda, \nu)).
\end{equation}

The next lemma is a sharper version of Lemma \ref{lem:multiplication
  formula}. 
\begin{lem}\label{lem:sharper multiplication formulae}
Suppose that $\lambda_i \in \Lambda^{+}$ and $\mu_i, \nu_i \in
\Lambda$ $(i=1,2)$, and let $ \nu_1 =  \mu_2$. Assume that
$\varphi_{S_1T_1}$, $\varphi_{S_2T_2} \in
\mathcal{C}^0(\Lambda)$ where $S_i \in
\mathcal{T}_0(\lambda_i,\mu_i)$, $T_i \in 
\mathcal{T}_0(\lambda_i,\nu_i)$. Then
\begin{equation*}\label{eq:sharper multiplication formulae}
\varphi_{S_1T_1} \cdot \varphi_{S_2T_2}
=
\begin{cases}
\underset{\varphi_{ST} \in
  \mathcal{C}^0(\lambda_1,0)}{\sum}
  r_{ST} \cdot\varphi_{ST}    
 + 
\underset{
\lambda \rhd \lambda_1 ~
}{\sum}
\underset{\varphi_{ST} \in
  \mathcal{C}^0(\lambda)}{\sum} 
& r_{ST} \cdot
\varphi_{ST}, \\

& \text{if} ~  \varphi_{S_1T_1} \in
  \mathcal{C}^0(\lambda_1,0),\\

& \\

\underset{
\lambda \unrhd \lambda_1 ~ }{\sum}
\underset{
\varphi_{ST} \in
  \mathcal{C}^0(\lambda,1)}{\sum} r_{ST} \cdot
\varphi_{ST}, &
\text{if} ~  \varphi_{S_1T_1} \in
  \mathcal{C}^0(\lambda_1,1),\\

\underset{
\lambda \unrhd \lambda_2  ~ }{\sum}
\underset{
\varphi_{ST} \in
  \mathcal{C}^0(\lambda)}{\sum} r_{ST} \cdot
\varphi_{ST}, &
\text{if} ~  \varphi_{S_2T_2} \in
  \mathcal{C}^0(\lambda_2,0),\\

\underset{
\lambda \unrhd \lambda_2 ~ }{\sum}
\underset{
\varphi_{ST} \in
  \mathcal{C}^0(\lambda,1)}{\sum} r_{ST} \cdot
\varphi_{ST}, &
\text{if} ~  \varphi_{S_2T_2} \in
  \mathcal{C}^0(\lambda_2,1), \\
\end{cases}
\end{equation*}
where $r_{ST} \in R$ and $\mathcal{C}^0(\lambda) =
\mathcal{C}^0(\lambda, 0) \cup \mathcal{C}^0(\lambda, 1)$ and all
$\lambda, S, T$ occurring in the above formulas are such that $\lambda 
\in 
\Lambda^{+}$,  
and the 
semistandard Tableaux $S, T$ with ${\rm{Type}}( S ) = \mu_{1}$ and
${\rm{Type}}( T ) = \nu_{2}$,
respectively. 

Hence the $\mathcal{S}^0(\Lambda)$ is a subalgebra of
$\mathcal{S}(\Lambda)$.

\end{lem}

\begin{proof}
By Lemma \ref{lem:multiplication formula}, $\varphi_{S_1T_1} \cdot
\varphi_{S_2T_2}$ is a linear combination of
$\varphi_{ST}$, where $S \in
\mathcal{T}_0(\lambda,\mu_1)$, $T \in\mathcal{T}_0(\lambda,\nu_2)$
with $\lambda \in 
\Lambda^{+}$, $\lambda \unrhd \lambda_1$ and $\lambda_2$. 

First assume that $\varphi_{S_1T_1} \in
\mathcal{C}^0(\lambda_1,0)$. Then for
the above $\varphi_{ST}$, we have $\bold{a}(\lambda) > \bold{a}(\mu)$ 
if $\alpha(\mu) \neq \alpha(\nu_{2})$ since $\varphi_{S_2T_2} \in
\mathcal{C}^0(\Lambda)$ and $\lambda
\unrhd \lambda_2$ (Note that $\mu \unrhd \nu \Rightarrow
\bold{a}(\mu) \geq \bold{a}(\nu)$ for $\mu, \nu \in
\Lambda$. cf. \eqref{eq:A certain fact of r-compositions}). So
$\varphi_{ST} \in 
\mathcal{C}^0(\Lambda)$, hence $\varphi_{ST} \in
\mathcal{C}^0(\lambda,0) \cup \mathcal{C}^0(\lambda,1)$. Moreover, the 
case $\varphi_{ST} \in
\mathcal{C}^0(\lambda_1,1)$ cannot happen. In fact, if it happens,
then 
$\alpha(\lambda_1) \neq \alpha(\mu_1)$ which contradicts our
assumption that $\varphi_{S_1T_1} \in
\mathcal{C}^0(\lambda_1,0)$. Therefore, the first equality holds.

The second equality is easy. In fact, assume that $\varphi_{S_1T_1} \in
\mathcal{C}^0(\lambda_1,1)$. Then $\varphi_{ST}$, as in the previous
case, are elements in the $\mathcal{C}^0(\lambda,1)$ since
$\lambda \unrhd \lambda_1$ and $\bold{a}(\lambda_1) >
\bold{a}(\mu_1)$. Hence we obtain the second one.

Next assume that $\varphi_{S_2T_2} \in
\mathcal{C}^0(\lambda_2,0)$. If $\alpha(\mu_1) \neq \alpha(\nu_2)$,
then the definition of
$\mathcal{C}^0(\lambda_2,0)$ and our assumption $\mu_2 =\nu_1$ implies 
that $\alpha(\mu_1) \neq \alpha(\nu_1)$. It follows that
$\bold{a}(\lambda_1) 
> \bold{a}(\mu_1)$ since $\varphi_{S_1T_1} \in
\mathcal{C}^0(\Lambda)$. Thus we have $\bold{a}(\lambda) >
\bold{a}(\mu_1)$ by $\lambda \unrhd \lambda_1$. This shows that
$\varphi_{ST}$ is in
$\mathcal{C}^0(\Lambda)$. Therefore $\varphi_{ST} \in
\mathcal{C}^0(\lambda,0) \cup \mathcal{C}^0(\lambda,1)$ and the third
equality holds. 

Finally, suppose that $\varphi_{S_2T_2} \in
\mathcal{C}^0(\lambda_2,1)$. If
$\alpha(\mu_1) = \alpha(\nu_1)$ then $\bold{a}(\lambda_2)
> \bold{a}(\mu_1)$ since $\varphi_{S_2T_2} \in
\mathcal{C}^0(\lambda_2,1)$ and $\mu_2 = \nu_1$ by our
assumption. Hence $\bold{a}(\lambda) > \bold{a}(\mu_1)$, since
$\lambda \unrhd \lambda_2$. On the other hand, if $\alpha(\mu_1) \neq 
\alpha(\nu_1)$ then also $\bold{a}(\lambda) > \bold{a}(\mu_1)$ and
$\varphi_{S_1T_1} \in
\mathcal{C}^0(\Lambda)$ and $\lambda \unrhd \lambda_1$. This argument 
means that $\varphi_{ST} \in \mathcal{C}^0(\lambda,1)$ and
hence the fourth equality holds. The lemma is proved.
\end{proof}

\addtocounter{subsection}{3}
\subsection{}\label{standardly based subsection}
Let $\mu \in \Lambda$. We define $\varphi_{\mu} \in
\mathcal{S}(\Lambda)$  as the identity map on $M^{\mu}$
and zero map on $M^{\kappa}$ with $\kappa \neq
\mu$. Moreover, let $1_{\mathcal{S}(\Lambda)}$ be the
identity element of $\mathcal{S}(\Lambda)$. Then, from
the definition,  we can
write $1_{\mathcal{S}(\Lambda)} = \sum_{\mu \in
  \Lambda} \varphi_{\mu}$. On the
other hand, since $\varphi_{\mu}(m_{\mu}) =
m_{\mu} \in M^{\mu \ast} \cap
M^{\mu}$, we have
\begin{equation*}
\varphi_{\mu}(m_{\mu}) = \underset{\substack{
S,T \in \mathcal{T}_0(\lambda,\mu)\\
\lambda \in \Lambda^{+}}}{\sum} r_{S,T} \cdot m_{ST}
\qquad (r_{S,T} \in R)
\end{equation*}
by \cite[Proposition 6.3]{DJMa}. This shows that 
\begin{equation}\label{eq:identity map on M^mu_natural}
\varphi_{\mu} = \underset{\substack{
S,T \in \mathcal{T}_0(\lambda,\mu)\\
\lambda \in \Lambda^{+}}}{\sum} r_{S,T} \cdot
\varphi_{ST}  \qquad (\text{ for any } \mu \in
\Lambda ).
\end{equation}
Thus all the $\varphi_{ST}$ in the right hand side are contained in 
$\mathcal{C}^0(\Lambda)$. Hence we have
\begin{equation}\label{id:identity element of S^0}
1_{\mathcal{S}(\Lambda)} = \underset{ \mu \in \Lambda}{\sum} \quad  
\underset{\substack{
S,T \in \mathcal{T}_0(\lambda,\mu)\\
\lambda \in \Lambda^{+}}}{\sum} r_{S, T} \cdot
\varphi_{ST} \in \mathcal{S}^0(\Lambda).
\end{equation}

For any $(\lambda, \varepsilon) \in \varOmega$, we define by
$\mathcal{S}_{0}^{\vee (\lambda, \varepsilon)} =
\mathcal{S}^0(\Lambda)(>(\lambda, \varepsilon))$ the $R$-submodule of     
$\mathcal{S}^0(\Lambda)$ spanned by
$\varphi_{UV}$ where $(U,V) \in I({\lambda}',{\varepsilon}')
\times J({\lambda}',{\varepsilon}')$ for some
$({\lambda}',{\varepsilon}') \in \varOmega$ with
$({\lambda}',{\varepsilon}') > 
(\lambda, \varepsilon)$. Note that $\mathcal{S}^0(\Lambda) \cap
\mathcal{S}^{\vee \lambda} =
\mathcal{S}_{0}^{\vee (\lambda,1)}$ for every $\lambda \in
\Lambda^{+}$. Similarly, we define
$\mathcal{S}^0(\Lambda)(\geq (\lambda, \varepsilon))$ as the
$R$-submodule spanned by $\varphi_{UV}$ with
$({\lambda}',{\varepsilon}') \geq (\lambda, \varepsilon)$.

We can now state.

\addtocounter{thm}{1}
\begin{thm}\label{thm:standardly based}
The subalgebra $\mathcal{S}^0(\Lambda)$ is standardly based on
$(\varOmega, \geq)$ with standard basis $\mathcal{C}^0(\Lambda)$, that  
is,

$({\rm{i}})$ The union $\underset{{(\lambda,\varepsilon) \in
    \varOmega}}{\bigcup} \mathcal{C}^0(\lambda,\varepsilon) =
    \mathcal{C}^0(\Lambda)$ is
    disjoint and forms an $R$-basis for $\mathcal{S}^0(\Lambda)$.

$({\rm{ii}})$ For any $\varphi \in \mathcal{S}^0(\Lambda)$,
$\varphi_{ST} \in \mathcal{C}^0(\lambda,\varepsilon)$, we have

\begin{equation}\label{eq:multiplication formulae of S^0(m,n)}
\begin{array}{c}
\varphi \cdot \varphi_{ST} \equiv \underset{S'
 \in
 I(\lambda,\varepsilon)}{\sum}f_{S',(\lambda,\varepsilon)}(\varphi,S) 
 \cdot \varphi_{S'T} \mod
 \mathcal{S}_{0}^{\vee (\lambda, \varepsilon)}\\
\varphi_{ST} \cdot \varphi \equiv \underset{T'
 \in
 J(\lambda,\varepsilon)}{\sum}f_{(\lambda,\varepsilon),T'}(T,\varphi) 
 \cdot \varphi_{ST'} \mod \mathcal{S}_{0}^{\vee (\lambda, \varepsilon)},
\end{array}
\end{equation}

where $\varphi_{S'T}, \varphi_{ST'} \in
\mathcal{C}^0(\Lambda)$ and $f_{S',(\lambda,\varepsilon)}(\varphi,S)$, 
$f_{(\lambda,\varepsilon),T'}(T,\varphi) \in R$ are
independent of $T$ and $S$, respectively.
\end{thm}
\begin{proof}
The first condition (i) is immediate from \eqref{eq:the property of
  C^0(m,n)}. We 
show (ii). Take $\varphi_{ST} \in
\mathcal{C}^0(\lambda,\varepsilon)$ and $\varphi \in
\mathcal{S}^0(\Lambda)$. Note that $\varphi_{ST}$ is an element in
  $\mathcal{C}(\Lambda)$ which is the cellular basis for
$\mathcal{S}(\Lambda)$ and $\varphi$ is an element in
  $\mathcal{S}(\Lambda)$. Hence, by the property of
cellular 
basis, $\varphi_{ST} \cdot \varphi$ can be written as    
\begin{equation*}
\varphi_{ST} \cdot \varphi = \underset{T' \in
 \mathcal{T}_0(\lambda) }{\sum}r_{T'} \cdot \varphi_{ST'} +
 \underset{\substack{{\lambda}' \in \Lambda^{+}, ~ {\lambda}' \rhd 
 \lambda\\ 
U,V \in \mathcal{T}_0({\lambda}')}}{\sum}r_{UV} \cdot \varphi_{UV}
\end{equation*}
with $r_{T'}$, $r_{UV} \in R$, where $r_{T'}$ does not depend on
$S$. Then by rewriting $\varphi$ as
a linear combination of the basis elements in
$\mathcal{C}^0(\Lambda)$ and by combining the formulas in Lemma
\ref{lem:sharper multiplication formulae}, we obtain the second
equality. By a similar argument, applying the third and fourth
formulas in Lemma \ref{lem:sharper multiplication formulae} instead,
the 
first equality also holds. The theorem follows.
\end{proof}

\addtocounter{subsection}{1}
\subsection{}(cf. {\cite{DR1}}).\label{def:the definition of
  standardly (full-)based algebra}
Let $A$ be a standardly based algebra on $\Lambda$ as given in
\ref{def:standardly based}. For any
$\lambda \in \Lambda$, let $f_{\lambda}: J(\lambda) \times I(\lambda) 
\rightarrow R$ be a function, whose value $f_{\lambda}(j,i')$ at
$(j,i') \in J(\lambda) \times I(\lambda)$ is defined by 

\begin{equation*}
a_{ij}^{\lambda} a_{i'j'}^{\lambda} \equiv f_{\lambda}(j,i')
a_{ij'}^{\lambda} \mod A(>\lambda).
\end{equation*}
The function $f_{\lambda}$ induces a bilinear form ${\beta}_{\lambda}:
A^{\lambda} \times A^{\lambda} \rightarrow R$ such that
${\beta}_{\lambda}(a_{ij}^{\lambda},a_{i'j'}^{\lambda}) =
f_{\lambda}(j',i)$, where $A^{\lambda}$ is the free $R$-submodule of
$A$ spanned by $a_{ij}^{\lambda}$ for all $(i,j) \in I(\lambda) \times 
J(\lambda)$. We say $A$ is a standardly \emph{full}-based algebra if
$ \text{Im}({\beta}_{\lambda})=R$ for all $\lambda \in \Lambda$. 
 
The following result has been proved by Du and Rui
\cite[(3.2.1), (4.2.7)]{DR2}. 

\addtocounter{thm}{1}
\begin{thm}[{\cite{DR2}}]\label{th:quasi-hereditary}

Suppose that $A$ is a standardly based algebra. 

$(\rm{i})$ Let $R$ be a commutative Noetherian ring. If $A$ is a
standardly  
full-based algebra, then $A$ is a quasi-hereditary algebra over $R$ in   
the sense of \cite{CPS}.

$(\rm{ii})$ If $R$ is a commutative local Noetherian ring, then $A$ is
split 
quasi-hereditary if and only if $A$ is a standardly full-based 
algebra. 
\end{thm}

\addtocounter{subsection}{1}
\subsection{}\label{sub:dagger tableau}
We shall show that the $\mathcal{S}^0(\Lambda)$ turns out to be a
standardly full-based algebra under a certain condition. For
$(\lambda,1)  
\in \varOmega$, put $\lambda = (\lambda^{(1)}, \ldots
,\lambda^{(r)})$, $ \lambda^{(i)} = (\lambda^{(i)}_1, \ldots
,\lambda^{(i)}_{n_i})$ and take the smallest positive integer $l$ such    
that $| \lambda^{(l)}| \neq 0 $ $(1 \leq l \leq r)$. Note that $l \neq
 r$ since $(\lambda, 1) \in \varOmega$. Then, one can
 define an $r$-partition $\lambda^{\dagger}$ associated to $\lambda$
 as follows, 
\begin{equation*}
\begin{array}{l}
\lambda^{\dagger} = (\lambda^{\dagger (1)}, \ldots ,\lambda^{\dagger
  (r)}), \\[3mm]

\lambda^{\dagger (l)} = \underset{|\lambda^{(l)}|-1 ~
  \text{times}}{(\underbrace{1,1, \ldots ,1})},  ~~  \lambda^{\dagger
  (l+1)} = \underset{|\lambda^{(l+1)}|+1 ~
  \text{times}}{(\underbrace{1,1, \ldots ,1})}, \\[10mm]

\lambda^{\dagger (j)} = \underset{|\lambda^{(j)}| ~
  \text{times}}{(\underbrace{1,1, \ldots ,1})} \quad  (j \neq l,
  l+1).   \qquad   \\[10mm]

(\text{If } |\lambda^{(j)}| = 0  ~~ \text{then }
  \lambda^{\dagger 
  (j)} ~ \text{is the empty partition.})
\end{array}
\end{equation*}    
We remark that $\lambda^{\dagger}$ is defined only when $(\lambda, 1)
\in 
\varOmega$, and it satisfies the property
\begin{equation}\label{geq:the property of lambda dagger}
\bold{a}(\lambda)> \bold{a}(\lambda^{\dagger}). 
\end{equation}
 
Moreover, let $T^{\lambda \dagger}$ be a semistandard
$\lambda$-Tableau of type $\lambda^{\dagger}$ in which the entries are 
laid in increasing order from ``left-upper'' to right along the
rows. (The ``left-upper'' means that the left has priority over the
upper.) For example, if
$\lambda=((3,2,1),(2,2))$ then $\lambda^{\dagger} = ( ( 1^5 ), ( 1^5 )
)$ and
\begin{center}
$T^{\lambda \dagger} = \biggl($ 
\begin{tabular}{ccc}
\begin{minipage}{13mm}
\setlength{\unitlength}{1.2mm}
\begin{picture}(21,6)
\put(0,0){\line(1,0){5}}
\put(0,2){\line(1,0){10}}
\put(0,4){\line(1,0){15}}
\put(0,6){\line(1,0){15}}

\put(0,0){\line(0,1){6}}
\put(5,0){\line(0,1){6}}
\put(10,2){\line(0,1){4}}
\put(15,4){\line(0,1){2}}

\put(0.5,3){\makebox(4,4){$\scriptscriptstyle{(1,1)}$}}
\put(5.5,3){\makebox(4,4){$\scriptscriptstyle{(2,1)}$}}
\put(10.5,3){\makebox(4,4){$\scriptscriptstyle{(3,1)}$}}
\put(0.5,1){\makebox(4,4){$\scriptscriptstyle{(4,1)}$}}
\put(5.5,1){\makebox(4,4){$\scriptscriptstyle{(5,1)}$}}
\put(0.5,-1){\makebox(4,4){$\scriptscriptstyle{(1,2)}$}}

\end{picture} 
\end{minipage}
& \quad , &
\begin{minipage}{13mm}
\setlength{\unitlength}{1.2mm}
\begin{picture}(21,6)
\put(0,2){\line(1,0){10}}
\put(0,4){\line(1,0){10}}
\put(0,6){\line(1,0){10}}

\put(0,2){\line(0,1){4}}
\put(5,2){\line(0,1){4}}
\put(10,2){\line(0,1){4}}

\put(0.5,3){\makebox(4,4){$\scriptscriptstyle{(2,2)}$}}
\put(5.5,3){\makebox(4,4){$\scriptscriptstyle{(3,2)}$}}
\put(0.5,1){\makebox(4,4){$\scriptscriptstyle{(4,2)}$}}
\put(5.5,1){\makebox(4,4){$\scriptscriptstyle{(5,2)}$}}

\end{picture} 
\end{minipage}
\end{tabular} 
$\biggl)$.
\end{center}
Thus, if $\mathfrak{s} \in \text{Std}(\lambda)$ and
$\lambda^{\dagger}(\mathfrak{s}) = T^{\lambda \dagger}$ then
$\mathfrak{s} = \mathfrak{t}^{\lambda}$. 

Finally, let  
\begin{equation*}
P_n(q,Q_1, \ldots ,Q_r) = \prod_{i=1}^{n} (1+q+ \cdots +q^{i-1}) \cdot    
\prod_{j=1}^{r-1}  \prod_{-n < k < n} (q^{2k} Q_j - Q_{j+1}).
\end{equation*}

  We have the following result.

\addtocounter{thm}{1}
\begin{prop}
Assume that $R$ is a commutative Noetherian ring. Suppose that
$P_n(q,Q_1, 
\ldots ,Q_r) \in R$ is
invertible, and $\lambda^{\dagger} \in \Lambda^{+}$ for any $(
\lambda, 1 ) \in \varOmega$. Then $\mathcal{S}^0(\Lambda)$
is the standardly full-based, and
hence is quasi-hereditary over $R$.
\end{prop}

\begin{proof}
For all $(\lambda, \varepsilon) \in \varOmega$, it is enough to show
that  
there exists some $T \in I(\lambda, \varepsilon) \cap
J(\lambda, \varepsilon)$ such that $(\varphi_{TT})^2 \equiv
r \cdot \varphi_{TT} \mod
\mathcal{S}_{0}^{\vee (\lambda,\varepsilon)}$ with $r \in R$
invertible. If $(\lambda, 0) \in \varOmega$ then one can take $T = 
T^{\lambda}$. Since
$\varphi_{T^{\lambda}T^{\lambda}}$ is the identity map on
$M^{\lambda}$, this case is immediate.

Assume that $(\lambda, 1) \in \varOmega$. First, we claim that   
\begin{equation}\label{eq:applying L to m}
m_{\lambda} \cdot (L_{|\lambda^{(l)}|} - Q_{l+1})
\equiv (q^{2k}Q_{l} - Q_{l+1}) \cdot m_{\lambda}
\mod \mathscr{H}^{\lambda}
\end{equation}
where $l$ is the integer attached to $\lambda$ as in \ref{sub:dagger
  tableau}, and $k$ is
some integer such that $|k| < n$ and $\mathscr{H}^{\lambda}$ is as in 
  \ref{sub:the definition of H^lambda}. In fact,
applying James-Mathas' result \cite[Proposition 3.7]{JM} for
$\mathfrak{t}^{\lambda}$ and $|\lambda^{(l)}|$, we see that
$m_{\lambda} \cdot L_{|\lambda^{(l)}|} \equiv q^{2k}Q_{l}
\cdot m_{\lambda} \mod \mathscr{H}^{\lambda}$ and hence
\eqref{eq:applying L to m} holds. 

Now, by \eqref{geq:the property
  of lambda dagger} and by our assumption, one can choose the element 
$\varphi_{T^{\lambda \dagger} T^{\lambda \dagger}} \in
\mathcal{C}^0(\lambda, 1)$. For any ${\mu}' \in
\Lambda$, we have $\varphi_{T^{\lambda \dagger} T^{\lambda
  \dagger}}(m_{{\mu}'}) =
{\delta}_{{\mu}' {\lambda}^{\dagger}}
\cdot m_{{\mathfrak{t}}^{\lambda} {\mathfrak{t}}^{\lambda}}
={\delta}_{{\mu}' {\lambda}^{\dagger}} \cdot m_{\lambda}$. Hence
\begin{equation}\label{eq:calculation 1}
\begin{array}{ll}
(\varphi_{T^{\lambda \dagger} T^{\lambda \dagger}})^2(m_{{\mu}'}) 

&= {\delta}_{{\mu}'
    {\lambda}^{\dagger}} (\varphi_{T^{\lambda \dagger}
  T^{\lambda \dagger}})(m_{\lambda}) 
={\delta}_{{\mu}' {\lambda}^{\dagger}}
    (\varphi_{T^{\lambda \dagger}
  T^{\lambda \dagger}})(u^{+}_{\lambda}
    x_{\lambda} )\\

&= {\delta}_{{\mu}' {\lambda}^{\dagger}}
    (\varphi_{T^{\lambda \dagger} T^{\lambda
  \dagger}})(m_{{\lambda}^{\dagger}} (L_{|{\lambda}^{(l)}|} -
    Q_{l+1}) x_{\lambda}) \\

&= {\delta}_{{\mu}' {\lambda}^{\dagger}} \cdot m_{\lambda}
    (L_{|{\lambda}^{(l)}|} - Q_{l+1}) x_{\lambda} 
\end{array}
\end{equation}
where the second and the third equality follows from the definition of 
$m_{\lambda}$ and $m_{{\lambda}^{\dagger}}$,
respectively. Note that $ m_{\lambda} x_{\lambda} =
P(q) m_{\lambda}$ where $P(q)$
is  some products of Poincar\'e
polynomials. Therefore, by \eqref{eq:applying L to m}, we have 
\begin{equation*}
(\varphi_{T^{\lambda \dagger} T^{\lambda
    \dagger}})^2(m_{{\mu}'})=
\left\{\begin{array}{ll}
P(q) \cdot (q^{2k}Q_{l} - Q_{l+1}) m_{\lambda} + h  & \text{ if 
}{\mu}' = {\lambda}^{\dagger} \\ 
0 & \text{ otherwise },\\
\end{array}\right.
\end{equation*}
where $h \in \mathscr{H}^{\lambda}$. 
Now, by definition, $({\varphi}_{T^{\lambda \dagger}
  T^{\lambda \dagger}})^2 (m_{{\lambda}^{\dagger}}) \in
  M^{{\lambda}^{\dagger} \ast} \cap
M^{{\lambda}^{\dagger}}$. 
So, by \cite[Proposition 6.3]{DJMa},
$({\varphi}_{T^{\lambda \dagger} T^{\lambda \dagger}})^2
  (m_{{\lambda}^{\dagger}})
= \sum r_{UV} m_{UV}$ where $r_{UV} \in R$ and the sum is over
$U, V \in \mathcal{T}_0(\alpha, {\lambda}^{\dagger})$ for some $\alpha 
\in \Lambda^{+}$. Hence, by \eqref{eq:calculation 1} and by a property 
of cellular basis, we deduce that
\begin{equation*}
\begin{array}{ll}
P(q) \cdot (q^{2k}Q_{l} - Q_{l+1}) m_{\lambda}
& \equiv ({\varphi}_{T^{\lambda \dagger} T^{\lambda \dagger}})^2
    (m_{{\lambda}^{\dagger}}) \mod
    \mathscr{H}^{\lambda} \\

&= m_{\lambda} (L_{|{\lambda}^{(l)}|} -
Q_{l+1}) x_{\lambda} 
= m_{T^{\lambda \dagger} T^{\lambda \dagger}} (L_{|{\lambda}^{(l)}|} - 
Q_{l+1}) x_{\lambda}\\ 
&= \sum_{T' \in \mathcal{T}_0(\lambda,
  {\lambda}^{\dagger})} r_{T'} m_{T^{\lambda \dagger} T'} +
\sum_{U', V'} r_{U'V'}
m_{U'V'}  
\end{array}
\end{equation*}
where $r_{T'}, r_{U'V'} \in R$ and the last sum is over $U', V' \in 
  \mathcal{T}_0(\alpha, {\lambda}^{\dagger})$ for $\alpha \in
  \Lambda^{+}$ such that $\alpha \rhd \lambda$. Comparing the
  coefficient of $m_{T^{\lambda \dagger} T'}$ on both sides reveals 
  that 
  $r_{T^{\lambda \dagger}} = P(q) \cdot (q^{2k}Q_{l} - Q_{l+1})$
  and $r_{T'} = 0$ unless $T' = T^{\lambda \dagger}$. Thus,
  $({\varphi}_{T^{\lambda \dagger} T^{\lambda \dagger}})^2
    (m_{{\lambda}^{\dagger}}) = P(q) \cdot (q^{2k}Q_{l} -
  Q_{l+1}) ({\varphi}_{T^{\lambda \dagger} T^{\lambda \dagger}})
  (m_{{\lambda}^{\dagger}}) + \sum_{U', V'} r_{U'V'}
{\varphi}_{U'V'} (m_{{\lambda}^{\dagger}})$. Consequently,
  $(\varphi_{T^{\lambda \dagger} T^{\lambda \dagger}})^2
   \equiv P(q) \cdot
  (q^{2k}Q_{l} - Q_{l+1}) \cdot
  \varphi_{T^{\lambda \dagger} T^{\lambda \dagger}} \mod
  \mathcal{S}_{0}^{\vee (\lambda, 1)}$. By our assumption, $P(q)
  \cdot (q^{2k}Q_{l} - Q_{l+1})$ is invertible and the proposition 
  follows. 
\end{proof}

\medskip

%%%%%%%%%%%%%%%%%%%%%%%%%%%%%%%%%%%%%%%%%%%%%%%%%%%%%%%%%%%%%%%%%%%%%%%%
%%%%%%%%%%%%%%%%%%%%%%%%%%%%%%%%%%%%%%%%%%%%%%%%%%%%%%%%%%%%%%%%%%%%%%%%

\section{A relationship between $\mathcal{S}^0(\Lambda)$ and
  $\mathcal{S}(\Lambda)$}

\subsection{}
We shall describe a relationship between the original cyclotomic
$q$-Schur 
algebra $\mathcal{S}(\Lambda)$ and its $R$-subalgebra
$\mathcal{S}^0(\Lambda)$.

Because $\mathcal{S}(\Lambda)$ is equipped with the involution $*$,
we can define subsets $\mathcal{S}^0(\Lambda)^*$,
$\mathcal{C}^0(\Lambda)^*$ of
$\mathcal{S}(\Lambda)$ by 
\begin{equation*}
\mathcal{S}^0(\Lambda)^* = \{ \varphi^* \mid \varphi \in
\mathcal{S}^0(\Lambda) \}, \quad \mathcal{C}^0(\Lambda)^* = \{
\varphi^{*}_{ST} \mid \varphi_{ST} \in \mathcal{C}^0(\Lambda) \},
\end{equation*}   
respectively. Then $\mathcal{S}^0(\Lambda)^*$ is an $R$-subalgebra of
$\mathcal{S}(\Lambda)$. Moreover, observe that
$\mathcal{S}^0(\Lambda)^*$ turns out to be a standardly based
algebra on $\mathcal{C}^0(\Lambda)^*$ by applying the involution
$*$ to Theorem \ref{thm:standardly based}.

\addtocounter{thm}{1}
\begin{prop}
$
\mathcal{S}(\Lambda) = \mathcal{S}^0(\Lambda) \cdot
\mathcal{S}^0(\Lambda)^* = \sum_{\lambda \in \Lambda^{+}}
\mathcal{S}^0(\Lambda) \varphi_{T^{\lambda} T^{\lambda}}
\mathcal{S}^0(\Lambda)^*
$.
\end{prop}

\begin{proof}
For all $\lambda \in \Lambda^{+}$, $\mu, \nu \in
\Lambda$ with $\mathcal{T}_0(\lambda, \mu)$, $\mathcal{T}_0(\lambda,
\nu) \neq \emptyset$ and for any $S \in \mathcal{T}_0(\lambda, \mu)$, $T
\in  
\mathcal{T}_0(\lambda, \nu)$, we certainly have
$\varphi_{ST^{\lambda}} \in 
\mathcal{C}^0(\Lambda)$ and hence $(\varphi_{TT^{\lambda}})^* =
\varphi_{T^{\lambda}T} \in
{\mathcal{C}^0(\Lambda)}^*$. Moreover, by using
\eqref{eq:decomposition of varphiST}, we see that all the basis
$\mathcal{C}(\Lambda)$ for
$\mathcal{S}(\Lambda)$ is contained in $\mathcal{S}^0(\Lambda) \cdot 
\mathcal{S}^0(\Lambda)^*$. Hence $\mathcal{S}(\Lambda) =
\mathcal{S}^0(\Lambda) \cdot
\mathcal{S}^0(\Lambda)^*$. Finally, for any $S \in
\mathcal{T}_0(\lambda, \mu)$, $T \in \mathcal{T}_0(\lambda, \nu )$, we   
see that $\varphi_{ST^{\lambda}}
\varphi_{T^{\lambda}T} = \varphi_{ST^{\lambda}}
\varphi_{T^{\lambda} T^{\lambda}}
\varphi_{T^{\lambda}T} \in \mathcal{S}^0(\Lambda)
\varphi_{T^{\lambda} T^{\lambda}}
\mathcal{S}^0(\Lambda)^*$, so the last equality follows.
\end{proof}

Next we introduce the Weyl module for
$\mathcal{S}^0(\Lambda)$. By \eqref{eq:multiplication
  formulae of S^0(m,n)} in Theorem
\ref{thm:standardly based}, it is easy to see that $R$-modules
$\mathcal{S}^0(\Lambda)(\geq (\lambda, \varepsilon))$ and
$\mathcal{S}_{0}^{\vee (\lambda, \varepsilon)} =
\mathcal{S}^0(\Lambda)(> (\lambda, \varepsilon))$ are two-sided 
ideals of $\mathcal{S}^0(\Lambda)$. Fix a $(\lambda, \varepsilon) \in 
\varOmega$. For $S \in I(\lambda, \varepsilon)$, we define the Weyl
module $Z^{(\lambda, \varepsilon)}_{S}$ for
$\mathcal{S}^0(\Lambda)$ by the $R$-submodule of $\{
\mathcal{S}^0(\Lambda)(\geq (\lambda, \varepsilon)) \} / \{
\mathcal{S}^0(\Lambda)(> (\lambda, \varepsilon)) \}$ with basis $\{
\varphi_{ST} + \mathcal{S}_{0}^{\vee (\lambda, \varepsilon)} \mid T
\in J(\lambda, \varepsilon)\}$. Moreover, by \eqref{eq:multiplication
  formulae of S^0(m,n)}, we see that $Z^{(\lambda, \varepsilon)}_{S}$
is the right 
$\mathcal{S}^0(\Lambda)$-module and the action of
$\mathcal{S}^0(\Lambda)$ on $Z^{(\lambda, \varepsilon)}_{S}$ is
independent of the choice of $S$, i.e, $Z^{(\lambda,
  \varepsilon)}_{S_1} \simeq Z^{(\lambda,\varepsilon)}_{S_2}$ for
all $S_1, S_2 \in I(\lambda, \varepsilon)$.
 However, since $T^{\lambda}$ is not an element in 
$I(\lambda, 1)$ for $(\lambda, 1) \in \varOmega$, one should
pay attention that
there is no ``canonical''-Weyl module for the case $(\lambda,
1)$. (That is, we can not define $Z^{(\lambda, 1)}_{T^{\lambda}}$.) 
For the convenience sake let $Z^{(\lambda,
  0)} = Z^{(\lambda,
  0)}_{T^{\lambda}}$ and put
$\varphi^{0}_{T} = \varphi_{T^{\lambda} T} +
\mathcal{S}_{0}^{\vee (\lambda, \varepsilon)}$ for any $T \in 
J(\lambda, 0) = \mathcal{T}_0^{+}(\lambda)$. We can now prove the
following.  
\begin{thm}\label{thm:tensor theorem}
Let $\lambda \in \Lambda^{+}$. Then there exists an isomorphism of
$\mathcal{S}(\Lambda)$-modules  
\begin{equation*}
Z^{(\lambda, 0)} \otimes_{\mathcal{S}^0(\Lambda)}
\mathcal{S}(\Lambda) \simeq W^{\lambda}
\end{equation*}
which maps $\varphi^{0}_{T^{\lambda}} \psi \otimes \varphi$ to
$\varphi_{T^{\lambda}} \psi \varphi$.
\end{thm}
\begin{proof}
First we note that $Z^{(\lambda, 0)} = {\varphi}^{0}_{T^{\lambda}} 
  {\mathcal{S}}^0(\Lambda)$. In fact, since $\mathcal{S}_{0}^{\vee
  (\lambda, 0)}$ is a two-sided ideal of $\mathcal{S}(\Lambda)$, for 
  any $T \in 
  \mathcal{T}^{+}_0(\lambda, \mu)$ with $\mu \in
  \Lambda$, we have ${\varphi}^{0}_{T^{\lambda}} \cdot
  \varphi_{T^{\lambda} T} = {\varphi}^{0}_{T}$ by
  \eqref{eq:decomposition of varphiST}. Now we can write
  $Z^{(\lambda, 0)} \otimes_{\mathcal{S}^0(\Lambda)}
\mathcal{S}(\Lambda) = ({\varphi}^{0}_{T^{\lambda}}
  \otimes_{\mathcal{S}^0(\Lambda)} 1) \cdot
  \mathcal{S}(\Lambda)$. On the other hand, we can also 
  write $W^{\lambda} =
  {\varphi}_{T^{\lambda}} \cdot
  \mathcal{S}(\Lambda)$ (note that
  ${\varphi}_{T^{\lambda}} =
  {\varphi}_{T^{\lambda} T^{\lambda}} +
  \mathcal{S}^{\vee \lambda}$). Hence it is enough
  to find two $\mathcal{S}(\Lambda)$-homomorphisms $f:
  {\varphi}_{T^{\lambda}} \cdot
  \mathcal{S}(\Lambda) \rightarrow ({\varphi}^{0}_{T^{\lambda}}
  \otimes_{\mathcal{S}^0(\Lambda)} 1) \cdot
  \mathcal{S}(\Lambda)$ and $g: ({\varphi}^{0}_{T^{\lambda}}
  \otimes_{\mathcal{S}^0(\Lambda)} 1) \cdot
  \mathcal{S}(\Lambda) \rightarrow
  {\varphi}_{T^{\lambda}} \cdot
  \mathcal{S}(\Lambda)$ such that
  $f({\varphi}_{T^{\lambda}}) = {\varphi}^{0}_{T^{\lambda}}
  \otimes_{\mathcal{S}^0(\Lambda)} 1$ and
  $g({\varphi}^{0}_{T^{\lambda}} \otimes_{\mathcal{S}^0(\Lambda)} 1) =  
  {\varphi}_{T^{\lambda}}$. 
We first define $f$. Consider the
  $\mathcal{S}(\Lambda)$-homomorphism $\widetilde{f}:
  {\varphi}_{T^{\lambda} T^{\lambda}} \cdot
  \mathcal{S}(\Lambda) \rightarrow ({\varphi}^{0}_{T^{\lambda}}
  \otimes_{\mathcal{S}^0(\Lambda)} 1) \cdot
  \mathcal{S}(\Lambda)$ satisfying
  $\widetilde{f}({\varphi}_{T^{\lambda} T^{\lambda}}) =
  {\varphi}^{0}_{T^{\lambda}} \otimes_{\mathcal{S}^0(\Lambda)}
  1$. Note that 
  $\widetilde{f}$ is a well defined homomorphism. To see this we take  
  an element ${\varphi} \in
  \mathcal{S}(\Lambda)$ such that
  ${\varphi}_{T^{\lambda} T^{\lambda}} {\varphi}
  = 0$. Write ${\varphi} = \sum r_{{\lambda}', ST} \cdot
  {\varphi}^{{\lambda}'}_{ST}$ for some $r_{{\lambda}', ST} \in
  R$ where ${\varphi}^{{\lambda}'}_{ST} =
  {\varphi}_{ST}$ with $S,T \in \mathcal{T}_0({\lambda}')$
  and the sum is taken over $r$-partitions ${\lambda}' \in
  \Lambda^{+}$ and semistandard tableaux $S,T \in
  \mathcal{T}_0({\lambda}')$. Then, since
  ${\varphi}_{T^{\lambda} T^{\lambda}}$ is the identity
  map on $M^{\lambda}$ and is zero on
  $M^{\kappa}$ for $\lambda \neq \kappa \in
  \Lambda$ and $\{ {\varphi}^{{\lambda}'}_{ST}
  \}$ is the basis for $\mathcal{S}(\Lambda)$,
  $r_{{\lambda}', ST} = 0$ whenever
  $S \in \mathcal{T}_0({\lambda}', \lambda), T \in
  \mathcal{T}_0({\lambda}', \nu)$ with ${\lambda}' \in
  \Lambda^{+}$ and $\nu \in
  \Lambda$. It follows that ${\varphi}^{\lambda'}_{S T^{\lambda'}} \in 
  \mathcal{S}^{0}(\Lambda)$ for all $S \in \mathcal{T}_0({\lambda}')$
  such 
  that $r_{\lambda', S T} \neq 0$. Hence, by using
  \eqref{eq:decomposition 
  of varphiST}, we have  
\begin{equation*}
\begin{array}{ll}
{\varphi}^{0}_{T^{\lambda}}
\otimes_{\mathcal{S}^0(\Lambda)}
{\varphi} 

&= {\varphi}^{0}_{T^{\lambda}}
\otimes_{\mathcal{S}^0(\Lambda)} ( \sum
r_{{\lambda}', ST} \cdot {\varphi}^{{\lambda}'}_{ST} ) \\[3mm]

&={\varphi}^{0}_{T^{\lambda}}
\otimes_{\mathcal{S}^0(\Lambda)} ( \sum
r_{{\lambda}', ST} \cdot {\varphi}^{{\lambda}'}_{S
  T^{{\lambda}'}} {\varphi}^{{\lambda}'}_{T^{{\lambda}'} T} )  \\[3mm]  

&= \sum
r_{{\lambda}', ST} ( {\varphi}^{0}_{T^{\lambda}}
\otimes_{\mathcal{S}^0(\Lambda)}
{\varphi}^{{\lambda}'}_{S T^{{\lambda}'}} \cdot
{\varphi}^{{\lambda}'}_{T^{{\lambda}'} T} ) \\[3mm]

&=\sum
r_{{\lambda}', ST} ( {\varphi}^{0}_{T^{\lambda}}
\cdot {\varphi}^{{\lambda}'}_{S T^{{\lambda}'}}
\otimes_{\mathcal{S}^0(\Lambda)}
{\varphi}^{{\lambda}'}_{T^{{\lambda}'} T} )   
\end{array}
\end{equation*}
where all sums are taken over ${\lambda}' \in
  \Lambda^{+}$ and $S,T \in
  \mathcal{T}_0({\lambda}')$. Now we have ${\varphi}^{0}_{T^{\lambda}}  
{\varphi}^{{\lambda}'}_{S T^{{\lambda}'}} = 0$ unless $S$ is the
  tableaux of type $\lambda$, and $r_{{\lambda}', ST} = 0$ if $S$ is
  of type $\lambda$. It follows that ${\varphi}^{0}_{T^{\lambda}}
\otimes_{\mathcal{S}^0(\Lambda)}
  {\varphi} = 0$ and so $\widetilde{f}$ is well-defined.
 
Take ${\varphi}^{\alpha}_{ST} \in
  \mathcal{S}(\Lambda)$ where $S, T \in
  \mathcal{T}_0(\alpha)$ with $\alpha \in \Lambda^{+}$ and
  suppose that $0 \neq {\varphi}_{T^{\lambda} T^{\lambda}}
  {\varphi}^{\alpha}_{ST} \in
  \mathcal{S}^{\vee \lambda}$. Then $S \in
  \mathcal{T}_0(\alpha, \lambda)$ and so ${\varphi}^{\alpha}_{ST} =
  {\varphi}_{T^{\lambda} T^{\lambda}}
  {\varphi}^{\alpha}_{ST}
  \in
  \mathcal{S}^{\vee \lambda}$ since
  ${\varphi}_{T^{\lambda} T^{\lambda}}$ is the identity map
  on $M^{\lambda}$. Thus, $\alpha \rhd \lambda$.
Furthermore, again by using  \eqref{eq:decomposition
  of varphiST} and the property that ${\varphi}_{S T^{\alpha}} \in
  \mathcal{S}^0(\Lambda)$, we have ${\varphi}^{0}_{T^{\lambda}}
\otimes_{\mathcal{S}^0(\Lambda)}
  {\varphi}^{\alpha}_{ST} 
= {\varphi}^{0}_{T^{\lambda}} \cdot {\varphi}^{\alpha}_{S T^{\alpha}} 
\otimes_{\mathcal{S}^0(\Lambda)}
  {\varphi}^{\alpha}_{T^{\alpha} T}$. Since
  ${\varphi}^{\alpha}_{S
  T^{\alpha}}$ is an element in the standard basis
  $\mathcal{C}^0(\Lambda)$ for 
  $\mathcal{S}^0(\Lambda)$ and $\alpha \rhd \lambda$ we
  have ${\varphi}^{}_{T^{\lambda} T^{\lambda}}
  {\varphi}^{\alpha}_{S T^{\alpha}} \in \mathcal{S}_{0}^{\vee
  (\lambda, 1)}$. Thus, ${\varphi}^{0}_{T^{\lambda}} \cdot
  {\varphi}^{\alpha}_{S T^{\alpha}} = 0$ and then
  ${\varphi}^{0}_{T^{\lambda}}
\otimes_{\mathcal{S}^0(\Lambda)}
  {\varphi}^{\alpha}_{ST} =0$. 
It follows that $\widetilde{f}(
  \mathcal{S}^{\vee \lambda} \cap
  {\varphi}_{T^{\lambda} T^{\lambda}} \cdot
  \mathcal{S}(\Lambda) ) = 0$. So
  $\widetilde{f}$ induces an
  $\mathcal{S}(\Lambda)$-homomorphism $f$ from
  ${\varphi}_{T^{\lambda}} \cdot
  \mathcal{S}(\Lambda)$ to $({\varphi}^{0}_{T^{\lambda}}
  \otimes_{\mathcal{S}^0(\Lambda)} 1) \cdot
  \mathcal{S}(\Lambda)$ such that
  $f({\varphi}_{T^{\lambda}}) = {\varphi}^{0}_{T^{\lambda}}
  \otimes_{\mathcal{S}^0(\Lambda)} 1$. 

Next we define $g$. To do this we define the
  $\mathcal{S}^0(\Lambda)$-balanced map $\widetilde{g}$ from the
  direct 
  product module ${\varphi}^{0}_{T^{\lambda}}
  \mathcal{S}^0(\Lambda) \times \mathcal{S}(\Lambda)$
  (for the right $\mathcal{S}^0(\Lambda)$-module
  ${\varphi}^{0}_{T^{\lambda}} 
  \mathcal{S}^0(\Lambda)$ and the left $\mathcal{S}^0(\Lambda)$-module   
  $\mathcal{S}(\Lambda)$) to ${\varphi}_{T^{\lambda}}
  \mathcal{S}(\Lambda)$ such that
  $\widetilde{g}( ( {\varphi}^{0}_{T^{\lambda}}
  {\varphi}', {\varphi} ) ) =
  ({\varphi}_{T^{\lambda}} {\varphi}' )
  {\varphi}$ for any ${\varphi}' \in
  \mathcal{S}^0(\Lambda)$ and ${\varphi} \in
  \mathcal{S}(\Lambda)$. Since the first two formulas in Lemma
  \ref{lem:sharper multiplication formulae} assert that
  ${\varphi}_{T^{\lambda} T^{\lambda}}
  {\varphi}' \in
  \mathcal{S}^{\vee \lambda}$ if ${\varphi}_{T^{\lambda} T^{\lambda}}
  {\varphi}' \in \mathcal{S}_{0}^{\vee (\lambda, 0)}$, we see that
  $\widetilde{g}$ is a well-defined $\mathcal{S}^0(\Lambda)$-balance 
  map. Then by using the universality property of the tensor product,
  $g$ is defined.  
\end{proof}

The following lemma is an analogy of Lemma 6.7 (i) in \cite{SawS}. 

\begin{lem}\label{lem:injective homomorphism f_lambda from
  Z^(lambda,0) to W^lambda_natural}
Suppose that $\lambda \in \Lambda^{+}$. We regard
$W^{\lambda}$ as right $\mathcal{S}^0(\Lambda)$-module via
the restriction. Then there exists an injective
$\mathcal{S}^0(\Lambda)$-module homomorphism $f_{\lambda} :
Z^{(\lambda, 0)} \rightarrow W^{\lambda}$ such that
  $f_{\lambda}( {\varphi}^{0}_{T}) =
  {\varphi}_{T}$ for $T \in J(\lambda, 0) =
  \mathcal{T}^{+}_0(\lambda)$.
\end{lem}

\begin{proof}
First, We can write $Z^{(\lambda, 0)} = \varphi_{T^{\lambda}
  T^{\lambda}} \mathcal{S}^{0}(\Lambda) / \{ \varphi_{T^{\lambda}
  T^{\lambda}} \mathcal{S}^{0}(\Lambda) \cap \mathcal{S}_{0}^{\vee
  (\lambda, 0)} \}$ and $W^{\lambda} = \varphi_{T^{\lambda}
  T^{\lambda}} \mathcal{S}(\Lambda) / \{ \varphi_{T^{\lambda}
  T^{\lambda}} \mathcal{S}(\Lambda) \cap \mathcal{S}^{\vee
  \lambda} \}$ and, moreover, there exists an inclusion map
\begin{equation*}
\varphi_{T^{\lambda}
  T^{\lambda}} \mathcal{S}^{0}(\Lambda) / \{ \varphi_{T^{\lambda}
  T^{\lambda}} \mathcal{S}^{0}(\Lambda) \cap \mathcal{S}^{\vee
  \lambda} \} \hookrightarrow \varphi_{T^{\lambda}
  T^{\lambda}} \mathcal{S}(\Lambda) / \{ \varphi_{T^{\lambda}
  T^{\lambda}} \mathcal{S}(\Lambda) \cap \mathcal{S}^{\vee
  \lambda} \} = W^{\lambda}
\end{equation*}
since $\mathcal{S}^{0}(\Lambda) \subset 
  \mathcal{S}(\Lambda)$. Then we have a natural surjective
  $\mathcal{S}^{0}(\Lambda)$-module 
  homomorphism 
\begin{equation*}
f : \varphi_{T^{\lambda}
  T^{\lambda}} \mathcal{S}^{0}(\Lambda) / \{ \varphi_{T^{\lambda}
  T^{\lambda}} \mathcal{S}^{0}(\Lambda) \cap \mathcal{S}_{0}^{\vee
  (\lambda, 0)} \}  
\rightarrow \varphi_{T^{\lambda}
  T^{\lambda}} \mathcal{S}^{0}(\Lambda) / \{ \varphi_{T^{\lambda}
  T^{\lambda}} \mathcal{S}^{0}(\Lambda) \cap \mathcal{S}^{\vee
  \lambda} \}
\end{equation*} 
since $\mathcal{S}_{0}^{\vee
  (\lambda, 0)} \subset \mathcal{S}^{\vee
  \lambda}$. But since $\mathcal{S}^{0}(\Lambda) \cap
  \mathcal{S}^{\vee 
  \lambda} = \mathcal{S}_{0}^{\vee (\lambda, 1)} \subset
  \mathcal{S}_{0}^{\vee (\lambda, 0)}$, we have $(
  \varphi_{T^{\lambda}  
  T^{\lambda}} \mathcal{S}^{0}(\Lambda) \cap \mathcal{S}^{\vee
  \lambda} ) \subset \mathcal{S}_{0}^{\vee (\lambda, 0)}$. Hence $(  
  \varphi_{T^{\lambda}  
  T^{\lambda}} \mathcal{S}^{0}(\Lambda) \cap \mathcal{S}^{\vee
  \lambda} ) \subset ( 
  \varphi_{T^{\lambda}  
  T^{\lambda}} \mathcal{S}^{0}(\Lambda) \cap \mathcal{S}_{0}^{\vee
  (\lambda, 0)} )$ and, therefore, $(  
  \varphi_{T^{\lambda}  
  T^{\lambda}} \mathcal{S}^{0}(\Lambda) \cap \mathcal{S}^{\vee
  \lambda} ) = ( 
  \varphi_{T^{\lambda}  
  T^{\lambda}} \mathcal{S}^{0}(\Lambda) \cap \mathcal{S}_{0}^{\vee
  (\lambda, 0)} )$. Consequently, $f$ is an isomorphism. The statement 
  for the basis now clear.
\end{proof}

Take a $(\lambda, 0) \in \varOmega$. Then $I(\lambda, 0) =
J(\lambda, 0) = \mathcal{T}^{+}_0(\lambda)$. Hence, by
\eqref{eq:multiplication formulae of S^0(m,n)} in Theorem
\ref{thm:standardly based}, we have the following. 

\addtocounter{subsection}{3}
\subsection{}
Suppose that $S, T \in \mathcal{T}^{+}_0(\lambda)$. Then there exists 
an element $r_{ST} \in R$ such that for any $U, V \in
\mathcal{T}^{+}_0(\lambda)$
\begin{equation*}
{\varphi}_{US} \cdot {\varphi}_{TV} \equiv
r_{ST} \cdot {\varphi}_{UV} \mod \mathcal{S}_{0}^{\vee (\lambda, 0)}. 
\end{equation*}

We define a bilinear form $\langle ~ ,
~ 
\rangle_{0} : Z^{(\lambda, 0)} \times Z^{(\lambda, 0)} \rightarrow R$
by $\langle {\varphi}^{0}_{S}, {\varphi}^{0}_{T} \rangle_{0} =
r_{ST}$. Hence we have
\begin{equation}\label{eq:the bilinear form on Z^(lambda,0)}
\langle {\varphi}^{0}_{S}, {\varphi}^{0}_{T} \rangle_{0} \cdot
  {\varphi}_{UV} \equiv 
  {\varphi}_{US} \cdot {\varphi}_{TV} \mod
  \mathcal{S}_{0}^{\vee (\lambda, 0)},
\end{equation}
where $U$ and $V$ are any elements of $\mathcal{T}^{+}_0(\lambda)$. It 
is easy to see that 
\begin{equation}\label{eq:the relationship between <,>_0 to
    <,>_{natural}} 
\langle {\varphi}^{0}_{S}, {\varphi}^{0}_{T} \rangle_{0} = \langle
  {\varphi}_{S}, 
  {\varphi}_{T} \rangle \qquad \text{ for every }
  S, T \in \mathcal{T}^{+}_0(\lambda).
\end{equation}
Let $\text{rad}Z^{(\lambda, 0)} = \{ x \in Z^{(\lambda, 0)} \mid
\langle x, y \rangle_{0} = 0 \text{ for all } y \in Z^{(\lambda, 0)}\}$.

\addtocounter{thm}{1}
\begin{lem}
${\rm{rad}}Z^{(\lambda, 0)}$ is an
$\mathcal{S}^{0}(\Lambda)$-submodule of $Z^{(\lambda, 0)}$.
\end{lem}

\begin{proof}
Take $x \in \text{rad}Z^{(\lambda, 0)}$ and $y \in Z^{(\lambda, 0)}$, 
  and write $x = \sum_{T_1 \in
  \mathcal{T}^{+}_0(\lambda)} r_{T_1}^{( x )} \cdot
  {\varphi}^{0}_{T_1}$ and $y = \sum_{T_2 \in
  \mathcal{T}^{+}_0(\lambda)} r_{T_2}^{( y )} \cdot
  {\varphi}^{0}_{T_2} \in Z^{(\lambda, 0)}$ with $r_{T_1}^{( x )}$,
  $r_{T_2}^{( y )} \in 
  R$. Then, for any ${\varphi} \in
  \mathcal{S}^{0}(\Lambda)$, we have $\langle x {\varphi}, y
  \rangle_{0} = \sum_{T_1, T_2 \in \mathcal{T}^{+}_0(\lambda)}
  r_{T_1}^{( x )} r_{T_2}^{( y )} \langle {\varphi}^{0}_{T_1}
  {\varphi}, {\varphi}^{0}_{T_2} \rangle_{0}$. We
  claim that 
\begin{equation}\label{eq:claim 1}
\langle {\varphi}^{0}_{T_1}
  {\varphi}, {\varphi}^{0}_{T_2} \rangle_{0} \cdot
  {\varphi}_{S_1 S_2} \equiv ( {\varphi}_{S_1
  T_1} {\varphi} ) {\varphi}_{T_2 S_2} \mod
  \mathcal{S}_{0}^{\vee (\lambda, 0)}
\end{equation}
 for any $S_1, S_2 \in
  \mathcal{T}^{+}_0(\lambda)$. In fact, since
  ${\varphi}_{S_1 T_1} \in \mathcal{C}^0(\lambda, 0)$ and
  ${\varphi} \in \mathcal{S}^{0}(\Lambda)$, Theorem
  \ref{thm:standardly based} implies that 
\begin{equation*}
{\varphi}_{S_1 
  T_1} \cdot {\varphi} \equiv \sum_{T \in \mathcal{T}^{+}_0(\lambda)} 
  f_{(\lambda, 0), T}(T_1, 
  {\varphi} ) \cdot {\varphi}_{S_1 T} \mod 
  \mathcal{S}_{0}^{\vee (\lambda, 0)} 
\end{equation*}
 where $f_{(\lambda, 0), T}(T_1,
  {\varphi} ) \in R$. Hence, 
\begin{equation*}
\begin{array}{ll}  
( {\varphi}_{S_1  
  T_1} \cdot {\varphi} ) \cdot {\varphi}_{T_2 
  S_2}  
&\equiv \underset{T \in \mathcal{T}^{+}_0(\lambda)}{\sum} f_{(\lambda,  
  0), T}(T_1,  
  {\varphi} ) \cdot {\varphi}_{S_1 T}\cdot  
  {\varphi}_{T_2 S_2} \\[6mm]
&\equiv \underset{T \in \mathcal{T}^{+}_0(\lambda)}{\sum} f_{(\lambda, 
  0), T}(T_1, 
  {\varphi} ) \cdot \langle {\varphi}^{0}_{T},
  {\varphi}^{0}_{T_2} \rangle_{0} \cdot
  {\varphi}_{S_1 S_2} \\[6mm]
&\equiv \langle \underset{T \in \mathcal{T}^{+}_0(\lambda)}{\sum}
  f_{(\lambda, 0), T}(T_1, {\varphi} ) \cdot
  {\varphi}^{0}_{T}, {\varphi}^{0}_{T_2} \rangle_{0}
  \cdot {\varphi}_{S_1 S_2}.
\end{array}
\end{equation*}
On the other hand, by the definition of ${\varphi}^{0}_{T_1}$
and by Theorem \ref{thm:standardly based}, we have 
\begin{equation*}
\begin{array}{ll}
{\varphi}^{0}_{T_1} \cdot {\varphi} 
&= ( {\varphi}_{T^{\lambda} T_1} +
\mathcal{S}_{0}^{\vee (\lambda, 0)} ) \cdot {\varphi}

= {\varphi}_{T^{\lambda} T_1} \cdot {\varphi} +
\mathcal{S}_{0}^{\vee (\lambda, 0)} \\[3mm]

&=\underset{T \in \mathcal{T}^{+}_0(\lambda)}{\sum} f_{(\lambda, 0), 
  T}(T_1, {\varphi} ) \cdot
{\varphi}_{T^{\lambda} T} +
\mathcal{S}_{0}^{\vee (\lambda, 0)} 

= \underset{T \in \mathcal{T}^{+}_0(\lambda)}{\sum} f_{(\lambda, 0), 
  T}(T_1, {\varphi} ) \cdot {\varphi}^{0}_{T}.
\end{array}
\end{equation*}
Note that
$\mathcal{S}_{0}^{\vee (\lambda, 0)}$ is a two-sided ideal of
$\mathcal{S}^{0}(\Lambda)$ in the second equality. \eqref{eq:claim 1}
follows from this. By Theorem \ref{thm:standardly based}, one can
write $\varphi \cdot \varphi_{T_2 S_2} = \sum_{S \in
  \mathcal{T}^{+}_0(\lambda)} f_{S, (\lambda, 0)}(\varphi, T_2)
\varphi_{S S_2}$. Then by \eqref{eq:claim 1}, we have, for $T_2$, $S_2
\in \mathcal{T}^{+}_0(\lambda)$,  
\begin{equation*}
\begin{array}{ll} 
\langle {\varphi}^{0}_{T_1} 
  {\varphi}, {\varphi}^{0}_{T_2} \rangle_{0} \cdot 
  {\varphi}_{S_1 S_2}  
&\equiv ( {\varphi}_{S_1 
    T_1} \cdot {\varphi} ) \cdot {\varphi}_{T_2 
    S_2} 
= {\varphi}_{S_1
    T_1} \cdot ( {\varphi} \cdot {\varphi}_{T_2
    S_2} ) \\

&\equiv {\varphi}_{S_1
    T_1} \cdot ( \underset{S \in \mathcal{T}^{+}_0(\lambda)}{\sum}
  f_{S, (\lambda, 0) } ( {\varphi}, T_2 ) \cdot
  {\varphi}_{S S_2} ) \\

&\equiv  \underset{S \in \mathcal{T}^{+}_0(\lambda)}{\sum}
  f_{S, (\lambda, 0) } ( {\varphi}, T_2 ) \cdot \langle
  {\varphi}^{0}_{T_1}, {\varphi}^{0}_{S} \rangle_{0}
  \cdot {\varphi}_{S_1 S_2}
\end{array}
\end{equation*}
The last equality follows from \eqref{eq:the bilinear form on
Z^(lambda,0)}. Since ${\varphi}_{S_1 S_2}$ are free $R$-basis of
$\mathcal{S}^{0}(\Lambda)$ and ${\varphi}_{S_1 S_2} \not\in
\mathcal{S}_{0}^{\vee (\lambda, 0)}$, we have 
\begin{equation*}
\langle {\varphi}^{0}_{T_1}
  {\varphi}, {\varphi}^{0}_{T_2} \rangle_{0} =
  \sum_{S \in \mathcal{T}^{+}_0(\lambda)} f_{S, (\lambda, 0) } (
  {\varphi}, T_2 ) \cdot \langle {\varphi}^{0}_{T_1},
  {\varphi}^{0}_{S} \rangle_{0}.
\end{equation*}
Hence
\begin{equation*}
\begin{array}{ll}
\langle x {\varphi}, y \rangle_{0} 
&= \underset{T_1, T_2 \in \mathcal{T}^{+}_0(\lambda)}{\sum} r_{T_1}^{( 
  x )} \cdot r_{T_2}^{( y )} \cdot \langle {\varphi}^{0}_{T_1}
{\varphi}, {\varphi}^{0}_{T_2}
\rangle_{0} \\[6mm]

&= \underset{T_1, T_2 \in \mathcal{T}^{+}_0(\lambda)}{\sum} r_{T_1}^{(  
  x )} r_{T_2}^{( y )} \cdot ( \underset{S \in
  \mathcal{T}^{+}_0(\lambda)}{\sum} f_{S, 
  (\lambda, 0) } ( {\varphi}, T_2 ) \cdot \langle
{\varphi}^{0}_{T_1}, {\varphi}^{0}_{S} \rangle_{0} )
\\[6mm]

&= \underset{T_2 \in \mathcal{T}^{+}_0(\lambda)}{\sum} r_{T_2}^{( y )}
\cdot ( \underset{S \in \mathcal{T}^{+}_0(\lambda)}{\sum} f_{S,
  (\lambda, 0) } ( {\varphi}, T_2 ) \cdot ( \underset{T_1
  \in \mathcal{T}^{+}_0(\lambda)}{\sum} r_{T_1}^{( x )} \cdot \langle
{\varphi}^{0}_{T_1}, {\varphi}^{0}_{S} \rangle_{0} )
~ ) \\[6mm]

&= \underset{T_2 \in \mathcal{T}^{+}_0(\lambda)}{\sum} r_{T_2}^{( y )} 
\cdot ( \underset{S \in \mathcal{T}^{+}_0(\lambda)}{\sum} f_{S,
  (\lambda, 0) } ( {\varphi}, T_2 ) \cdot \langle x,
{\varphi}^{0}_{S}
\rangle_{0} ) \\[6mm]

&= \langle x, \underset{T_2 \in \mathcal{T}^{+}_0(\lambda)}{\sum}
r_{T_2}^{( y )} \cdot ( \underset{S \in
  \mathcal{T}^{+}_0(\lambda)}{\sum} f_{S, 
  (\lambda, 0) } ( {\varphi}, T_2 ) \cdot
{\varphi}^{0}_{S} ) \rangle_{0}.
\end{array}
\end{equation*}
Since $\underset{T_2 \in \mathcal{T}^{+}_0(\lambda)}{\sum}
r_{T_2}^{( y )} \cdot ( \underset{S \in
  \mathcal{T}^{+}_0(\lambda)}{\sum} f_{S, 
  (\lambda, 0) } ( {\varphi}, T_2 ) \cdot
{\varphi}^{0}_{S} ) \in Z^{(\lambda, 0)}$ and $x \in
\text{rad}Z^{(\lambda, 0)}$, the last formula is equal to
$0$. This shows that $x {\varphi} \in \text{rad}Z^{(\lambda,
  0)}$. It is immediate that if $x_1, x_2 \in \text{rad}Z^{(\lambda,
  0)}$ then $x_1 + x_2 \in \text{rad}Z^{(\lambda, 0)}$. Thus, the
lemma follows.  
\end{proof}

We put $L_{0}^{\lambda} = Z^{(\lambda, 0)} /
\text{rad}Z^{(\lambda, 0)}$. Then we have the following.  

\begin{prop}\label{prop:the property of L^lambda and radZ^(lambda,0)}
Suppose that $R$ is a field, and $\lambda \in \Lambda^{+}$. Then 

$({\rm{i}})$ $L_{0}^{\lambda} \neq 0$ and

$({\rm{ii}})$ ${\rm{rad}}Z^{(\lambda, 0)}$ is the unique maximal
submodule of $Z^{(\lambda, 0)}$ and $L_{0}^{\lambda}$ is absolutely
irreducible. Moreover, the
Jacobson radical of $Z^{(\lambda, 0)}$ is equal to
$\text{rad}Z^{(\lambda, 0)}$.
\end{prop}

\begin{proof}
For any $\lambda \in \Lambda^{+}$, we have
  ${\varphi}^{0}_{T^{\lambda}} \in 
  Z^{(\lambda, 0)}$ and ${\varphi}^{0}_{T^{\lambda} T^{\lambda}}
  {\varphi}^{0}_{T^{\lambda} T^{\lambda}} ={\varphi}^{0}_{T^{\lambda}
  T^{\lambda}}$. Hence $\langle {\varphi}^{0}_{T^{\lambda}},
  {\varphi}^{0}_{T^{\lambda}} \rangle_{0} = 1 \neq 0$. Thus,
  ${\varphi}^{0}_{T^{\lambda}} \not\in \text{rad}Z^{(\lambda, 0)}$,
  and we have $L_{0}^{\lambda} 
  \neq 0$. This proves (i). Let $x \in Z^{(\lambda, 0)} \setminus
  \text{rad}Z^{(\lambda,0)}$. Then $\langle x, y \rangle_{0} \neq 0$
  for some $y \in Z^{(\lambda, 0)}$. We write $y = \sum_{S \in
  \mathcal{T}^{+}_0(\lambda)} r_S \cdot {\varphi}^{0}_{S}$. For each 
  $T \in \mathcal{T}^{+}_0(\lambda)$ set $y_{T} = \sum_{S \in
  \mathcal{T}^{+}_0(\lambda)} r_S \cdot {\varphi}_{S T}$,
  an element of $\mathcal{S}^{0}(\Lambda)$. Moreover, write $x = 
  \sum_{U \in 
  \mathcal{T}^{+}_0(\lambda)} r'_U \cdot {\varphi}^{0}_{U}$. We note
  that ${\varphi}^{0}_{U} \cdot {\varphi}_{ST} = \langle
  {\varphi}^{0}_{U}, {\varphi}^{0}_{S} \rangle_{0} \cdot
  {\varphi}^{0}_{T}$ 
  in $Z^{(\lambda, 0)}$ by \eqref{eq:the bilinear form on
  Z^(lambda,0)}. Thus we 
  have 
\begin{equation*}
\begin{array}{ll}
x y_{T} 
&= \underset{S, U \in
  \mathcal{T}^{+}_0(\lambda)}{\sum} r_{S} r'_{U} \cdot
{\varphi}^{0}_{U} 
  {\varphi}_{ST} 
= \underset{S, U \in
  \mathcal{T}^{+}_0(\lambda)}{\sum} r_{S} r'_{U} \cdot \langle
  {\varphi}^{0}_{U}, {\varphi}^{0}_{S} \rangle_{0} \cdot 
  {\varphi}^{0}_{T} \\[6mm]

&= \langle \underset{U \in
  \mathcal{T}^{+}_0(\lambda)}{\sum} r'_{U} {\varphi}^{0}_{U}, ~
  \underset{S \in  
  \mathcal{T}^{+}_0(\lambda)}{\sum} r_{S} {\varphi}^{0}_{S}
  \rangle_{0}  
  \cdot {\varphi}^{0}_{T} = \langle x, y 
  \rangle_{0} \cdot {\varphi}^{0}_{T}.
\end{array}
\end{equation*}
Now $R$ is a field and
  $\langle x, y \rangle_{0} \neq 0$. Hence we have $x \cdot
( 1 /  \langle x, 
  y \rangle_{0} ) \cdot 
 y_{T} = {\varphi}^{0}_{T}$ for any $T \in
  \mathcal{T}^{+}_0(\lambda)$. Consequently, $x$ generates
  $Z^{(\lambda, 0)}$ as an $\mathcal{S}^{0}(\Lambda)$-module. This
  argument can be applied to any
  element of $Z^{(\lambda, 0)}$ which does not belong to the
  radical. It follows that $\text{rad}Z^{(\lambda, 0)}$ is the unique
  maximal proper submodule
  of $Z^{(\lambda, 0)}$ and $L_{0}^{\lambda}$ is irreducible. The same
  argument shows that $L_{0}^{\lambda}$ is irreducible for any
  extension field of $R$, so $L_{0}^{\lambda}$ is absolutely
  irreducible. This proves (ii).
\end{proof}

\medskip

%%%%%%%%%%%%%%%%%%%%%%%%%%%%%%%%%%%%%%%%%%%%%%%%%%%%%%%%%%%%%%%%%%%%%%%%
%%%%%%%%%%%%%%%%%%%%%%%%%%%%%%%%%%%%%%%%%%%%%%%%%%%%%%%%%%%%%%%%%%%%%%%%

\section{A relationship between $\mathcal{S}^{\flat}(\bs{m},n)$ and
  $\mathcal{S}^0(\Lambda)$}\label{sec:A relationship between
  S^flat(m,n) and S^0(Lambda)} 

First, we recall the definition of modified Ariki-Koike algebras and
their 
cyclotomic $q$-Schur algebras (\cite{SawS}).

\subsection{}\label{sec:modified Ariki-Koike algebra}
From now on, throughout this paper, we consider the following
condition on parameters $Q_1, \ldots ,
Q_{r}$ in $R$ whenever we consider the modified Ariki-Koike algebras
(and their cyclotomic $q$-Schur algebras).

\begin{equation}\label{eq:The condition of the Modified Ariki-Koike}
Q_{i} - Q_{j} \text{ are invertible in } R \text{ for any } i \neq
j. \qquad \qquad \qquad \qquad \qquad \qquad \qquad \qquad
\end{equation}

Let $A$ be a square matrix of degree $r$ whose $i$-$j$ entry is given
by $Q_{j}^{i-1}$ for $1 \le i, j \le r$. Thus $A$ is the Vandermonde
matrix, and $\Delta = \det{A} = \prod_{i > j} ( Q_i -Q_j )$ is
invertible by \eqref{eq:The condition of the Modified Ariki-Koike}. We 
express the inverse of $A$ as $A^{-1} = \Delta^{-1} B$ with $B = (
h_{ij} )$, and define a polynomial $F_{i}( X ) \in R[ X ]$, for $1 \le
i \le r$, by $F_{i}( X ) =\sum_{ 1 \le j \le r } h_{ij} X^{j-1}$.

The modified Ariki-Koike algebra $\mathscr{H}^{\flat} =
 \mathscr{H}^{\flat}_{n, r}$ is an associative algebra over $R$ with
 generators $T_2, \cdots , T_n$ and $\xi_{1}, \ldots , \xi_{n}$ and
 relations  
\begin{equation}\label{eq:Defining relations of H^flat}
\begin{array}{ll}
( T_i - q )( T_i + q^{-1} ) 
= 0 
& ( 2 \le i \le n), \\

(\xi_{i} - Q_{1} ) \cdots ( \xi_{i} - Q_{r} )
= 0
& ( 1 \le i \le n), \\

T_{i} T_{i+1} T_{i}
= T_{i+1} T_i T_{i+1}
& ( 2 \le i \le n), \\

T_{i} T_{j}
= T_{j} T_{i}
& ( | i - j | \geq 2 ), \\

\xi_{i} \xi_{j}
= \xi_{j} \xi_{i}
& ( 1 \le i, j \le n ), \\

T_{j} \xi_{j}
= \xi_{j-1} T_j + \Delta^{-2}
\underset{c_1 < c_2}{\sum}
(Q_{c_2} - Q_{c_1} )( q - q^{-1} ) F_{c_1}( \xi_{j-1} ) F_{c_2}(
\xi_{j} ), 
& \\

T_{j} \xi_{j-1} 
= \xi_{j} T_{j} - \Delta^{-2}
\underset{c_1 < c_2}{\sum}
(Q_{c_2} - Q_{c_1} )( q - q^{-1} ) F_{c_1}( \xi_{j-1} ) F_{c_2}(
\xi_{j} ),
& \\

T_{j} \xi_{k}
=\xi_{j} T_{j}
& ( k \neq j-1, j ).
\end{array}
\end{equation}

It is known that if $R = \mathbb{Q}( \overline{q}, \overline{Q}_{1},
\ldots 
, \overline{Q}_{r} )$, the field of rational functions with variables
$\overline{q}, \overline{Q}_{1}, \ldots
, \overline{Q}_{r}$, $\mathscr{H}^{\flat}$ is isomorphic to
$\mathscr{H}$, and it gives an alternate presentation of $\mathscr{H}$
apart from \ref{def:the definition of Ariki-Koike algebras}.

The subalgebra $\mathscr{H}^{\flat}( \mathfrak{S}_{n} )$ of
$\mathscr{H}^{\flat}$ generated by $T_2, \ldots , T_n$ is isomorphic
to $\mathscr{H}_{n}$, hence it can be naturally identified with the
corresponding subalgebra $\mathscr{H}( \mathfrak{S}_{n} )$ of
$\mathscr{H}$. Moreover, it is known by \cite{S} that the set $\{
\xi_{1}^{c_1} \cdots \xi_{n}^{c_n} T_{w} \mid w \in
\mathfrak{S}_{n}, ~ 0 \le c_{i} < r \text{ for } 1 \le i \le n \}$
gives rise to a basis of $\mathscr{H}^{\flat}$.

Let $V = \bigoplus_{i = 1}^{r} V_{i}$ be a free $R$-module, with rank 
$V_{i} = m_{i}$. We put $m = \sum m_{i}$. It is known by \cite{SakS}
that we can define a right $\mathscr{H}$-module structure on
$V^{\otimes n}$. We denote this representation by ${\rho} :
\mathscr{H} \rightarrow \End{V^{\otimes n}}$. Note that this
construction works without the condition \eqref{eq:The
  condition of the Modified Ariki-Koike}. Also it is shown in
\cite{S} that, under the assumption \eqref{eq:The
  condition of the Modified Ariki-Koike}, a right action
of 
$\mathscr{H}^{\flat}$ on $V^{\otimes n}$ can be defined. We denote
this 
representation by ${\rho}^{\flat} : \mathscr{H}^{\flat} \rightarrow
\End{V^{\otimes n}}$. By \cite[Lemma 3.5]{S}, we know that
$\text{Im}{\rho} \subset \text{Im}{{\rho}^{\flat}}$.

We consider the condition
\begin{equation}\label{con:m_i>=n}
m_{i} \geq n \text{ ~ for ~ } i = 1, \cdots, r. \qquad \qquad \qquad
\qquad \qquad \qquad \qquad \qquad \qquad \qquad \qquad
\end{equation}

\addtocounter{thm}{1}
\begin{lem}[{\cite[Lemma 1.5]{SawS}}]
Under the conditions \eqref{eq:The condition of the Modified
  Ariki-Koike}, \eqref{con:m_i>=n}, there exists
an $R$-algebra homomorphism $\rho_{0} : \mathscr{H} \rightarrow
\mathscr{H}^{\flat}$ such that $\rho_{0}$ induces the identity on
$\mathscr{H}_{n}$. $($Here we regard $\mathscr{H}_{n} \subset
\mathscr{H}$, $\mathscr{H}_{n} \subset \mathscr{H}^{\flat}$ under the    
previous identifications.$)$ If ${\rm{Im}}{{\rho}^{\flat}} =
  {\rm{Im}}{\rho}$ and $R$ is a field, then $\mathscr{H} \simeq
\mathscr{H}^{\flat}$. 
\end{lem}

From now on, throughout the paper, we fix an $r$-tuple $\bs{m} = (
m_1, \ldots , m_r )$ of non-negative integers and always assume the
condition \eqref{con:m_i>=n} whenever we consider
$\mathscr{H}^{\flat}$.

Any $\mu \in \widetilde{\mathcal{P}}_{n,r}( \bs{m} )$ may be regarded
as an 
element in $\mathcal{P}_{n,1}$ (i.e, $1$-composition) of $n$ by
arranging the entries of $\mu = ( \mu_{j}^{( i )} )$ in order
\begin{equation*}
\mu_{1}^{( 1 )}, \ldots , \mu_{m_1}^{( 1 )}, \mu_{1}^{( 2 )}, \ldots ,
\mu_{m_2}^{( 2 )}, \ldots , \mu_{1}^{( r )}, \ldots , \mu_{m_r}^{( r
  )}, 
\end{equation*}    
which we denote by $\{ \mu \}$.

For $\alpha = ( n_1, \ldots , n_r ) \in \mathbb{Z}_{\geq 0}$ such that 
$\sum n_i = n$, we define $c( \alpha )$ by
\begin{equation*}
c( \alpha ) = ( \underset{n_1\text{-times}}{\underbrace{r, \ldots , 
    r}}, 
\underset{n_2\text{-times}}{\underbrace{r-1, \ldots , r-1}}, \ldots ,    
\underset{n_r\text{-times}}{\underbrace{1, \ldots , 1}} )
\end{equation*}
and let $c( \alpha ) = ( c_1, \ldots , c_n)$. We define $F_{\alpha}
\in 
\mathscr{H}^{\flat}$ by $F_{\alpha} = \Delta^{-n} F_{c_1}( \xi_{1} ) 
F_{c_2}( \xi_{2} ) \cdots
F_{c_n}( \xi_{n} )$. 
For any $\mu \in \widetilde{\mathcal{P}}_{n,r}$, put $m^{\flat}_{\mu}
= 
F_{\alpha( \mu )} \cdot m_{ \{ \mu \} }$ where $m_{\{ \mu \}} =
\sum_{ w \in \mathfrak{S}_{ \{ \mu \} } } q^{l( w )} T_{w} ~ ( =
x_{\mu}) \in \mathscr{H}_{n}$.

We define
an $R$-linear anti-automorphism $h \rightarrow h^{*}$ on
$\mathscr{H}^{\flat}$ by the condition that $*$ fixes the
generators $T_{i} ~ ( 2 \le i \le n )$ and $\xi_{j} ~ ( 1 \le j \le n 
)$. As discussed in \cite[2.7]{SawS}, this condition induces a
well-defined 
anti-automorphism on $\mathscr{H}^{\flat}$. Moreover, by Lemma 2.9 in
\cite{SawS}, we know that $(
m^{\flat}_{\mu} )^{*} = m^{\flat}_{\mu}$.
 For $\mathfrak{s}, \mathfrak{t} \in \text{Std}( \lambda )$ with
$\lambda \in \mathcal{P}_{n,r}$, we define an element
$m^{\flat}_{\mathfrak{s} \mathfrak{t}} \in \mathscr{H}^{\flat}$ by
$m^{\flat}_{\mathfrak{s} \mathfrak{t}} = T^{*}_{d( \mathfrak{s} )}
m^{\flat}_{\mu} T_{d( \mathfrak{t} )}$. By the above fact, we have
$( m^{\flat}_{\mathfrak{s} \mathfrak{t}} )^{*} =
m^{\flat}_{\mathfrak{t} \mathfrak{s}}$.

\begin{thm}[{\cite[Theorem 2.18]{SawS}}]
The modified Ariki-Koike algebra $\mathscr{H}^{\flat}$ is free as an
$R$-module with cellular basis $\{ m^{\flat}_{\mathfrak{s}
  \mathfrak{t}} \mid \mathfrak{s}, \mathfrak{t} \in \text{Std}(
\lambda ) \text{ for some } \lambda \in \mathcal{P}_{n,r} \}$.
\end{thm}

Put $M_{\flat}^{\mu} = m^{\flat}_{\mu} \mathscr{H}^{\flat}$ for $\mu
\in \widetilde{\mathcal{P}}_{n,r}$. We define a cyclotomic $q$-Schur
algebra 
$\mathcal{S}^{\flat}(\bs{m},n)$ as follows.

\begin{defi}
The cyclotomic $q$-Schur algebra for $\mathscr{H}^{\flat}$
with weight poset  
$\widetilde{\mathcal{P}}_{n,r}$ is the
endomorphism algebra
\begin{equation*}
\mathcal{S}^{\flat}(\bs{m},n) = \End_{\mathscr{H}^{\flat}}( M^{\flat}(  
\widetilde{\mathcal{P}}_{n,r} ) ), \qquad \text{ where } M^{\flat}(
\widetilde{\mathcal{P}}_{n,r} ) = \underset{\mu \in
  \widetilde{\mathcal{P}}_{n,r}}{\bigoplus} M_{\flat}^{\mu}. 
\end{equation*}
\end{defi}

For an $r$-tuples $\alpha \in \widetilde{\mathcal{P}}_{n,1}$, let
$M_{\flat}^{\alpha} = \bigoplus_{\mu ; \alpha( \mu ) = \alpha}
M_{\flat}^{\mu}$. Then by Proposition 5.2 (i) in \cite{SawS}, we 
have $\mathcal{S}^{\flat}(\bs{m},n) \simeq \bigoplus_{\alpha \in
  \widetilde{\mathcal{P}}_{n,1}} \End_{\mathscr{H}^{\flat}}
M_{\flat}^{\alpha}$ as $R$-algebras.

\begin{thm}[{\cite[Theorem 5.5]{SawS}}]\label{weak separation
    condition} 
Let $\mathcal{S}^{\flat}(\bs{m},n)$ be the cyclotomic $q$-Schur
algebra associated to the modified Ariki-Koike algebra
$\mathscr{H}^{\flat}$ and $\mathcal{S}( m_{i}, n_{i} )$ be the
$q$-Schur algebra associated to the Iwahori-Hecke algebra
$\mathscr{H}_{n_{i}}$. Then there
exists an isomorphism of $R$-algebras
\begin{equation*}
\mathcal{S}^{\flat}(\bs{m},n) \simeq 
\underset{
\substack{
( n_1, \ldots , n_r ) \\
n = n_1 + \cdots + n_{r}
}
}{\bigoplus} 
\mathcal{S}( m_{1}, n_{1} ) \otimes \cdots \otimes
\mathcal{S}( m_{r}, n_{r} ).
\end{equation*}
\end{thm}

Let $\mu, \nu \in \widetilde{\mathcal{P}}_{n,r}$ and $\lambda \in
\mathcal{P}_{n,r}$. We assume that $\alpha( \mu ) = \alpha( \nu ) =
\alpha( \lambda )$. For $S \in \mathcal{T}^{+}_0(\lambda, \mu)$ and $T  
\in \mathcal{T}^{+}_0(\lambda, \nu)$, put
\begin{equation*}
m^{\flat}_{ST} = \underset{
\substack{
\mathfrak{s}, \mathfrak{t} \in \text{Std}( \lambda ) \\
\mu( \mathfrak{s} ) = S, ~ \nu( \mathfrak{t} ) = T 
}
}
{\sum}
q^{l( d( \mathfrak{s} ) ) + l( d( \mathfrak{t} ) )}
m^{\flat}_{\mathfrak{s} 
  \mathfrak{t}}. 
\end{equation*}

Moreover, for $S \in \mathcal{T}^{+}_0(\lambda, \mu)$ and $T \in
\mathcal{T}^{+}_0(\lambda, \nu)$, one can define
$\varphi^{\flat}_{ST} \in \mathcal{S}^{\flat}(\bs{m},n)$ by
$\varphi^{\flat}_{ST}( m^{\flat}_{\alpha} h ) = \delta_{\alpha \nu}
m^{\flat}_{ST} h$, for all $h \in \mathscr{H}^{\flat}$ and all $\alpha 
\in \widetilde{\mathcal{P}}_{n,r}$.

\begin{thm}[{\cite[Theorem 5.9]{SawS}}]
The cyclotomic $q$-Schur algebra $\mathcal{S}^{\flat}(\bs{m},n)$ is
free as an $R$-module with cellular basis
$\mathcal{C}^{\flat}(\bs{m},n) = \{ \varphi^{\flat}_{ST} \mid S, T \in 
\mathcal{T}^{+}_0(\lambda), \text{ for some } \lambda \in
\mathcal{P}_{n,r} \}$. 
\end{thm}

\addtocounter{subsection}{5}
\subsection{}
Let $\mathcal{S}^0(\Lambda)$ be as in Section \ref{The standard basis
  for S^0(Lambda)}. We describe a relationship between the algebra
$\mathcal{S}^0(\Lambda)$ and the cyclotomic $q$-Schur algebra
$\mathcal{S}^{\flat}(\bs{m},n)$ in the case where $\Lambda =
  \widetilde{\mathcal{P}}_{n,r}$. But in the moment, we shall consider 
  an arbitrary $\Lambda$ as in Section \ref{The standard basis
  for S^0(Lambda)}. 

First, let $\mathcal{C}^{00}(\Lambda) = \{
{\varphi}_{ST} \mid (S, T) \in I(\lambda, 1) \times
J(\lambda, 1),  ~  \lambda \in \Lambda^{+} \} \subset
\mathcal{C}^{0}(\Lambda)$ and $\mathcal{S}^{00}(\Lambda)$ be the
$R$-span of ${\varphi}_{ST} \in
\mathcal{C}^{00}(\Lambda)$, which is an $R$-submodule of
$\mathcal{S}^{0}(\Lambda)$. We note that, $\mathcal{S}^{00}(\Lambda)$ 
is a two-sided ideal of $\mathcal{S}^{0}(\Lambda)$ by the second
and fourth formula in Lemma \ref{lem:sharper multiplication formulae}.
Thus one can define the quotient algebra
$\overline{\mathcal{S}^{0}}(\Lambda) = \mathcal{S}^0(\Lambda) /
\mathcal{S}^{00}(\Lambda)$. We write $\overline{x} = x +
\mathcal{S}^{00}(\Lambda)$ $( x \in \mathcal{S}^0(\Lambda) )$. It is 
easy to see that $\overline{\mathcal{S}^{0}}(\Lambda)$ has a free
$R$-basis $\{ \overline{{\varphi}}_{ST} \mid S \in
I(\lambda, 0), ~ T \in J(\lambda, 0), ~ \lambda \in \Lambda^{+}
\}$. Note that the condition $(S,T) \in I(\lambda, 0) \times
J(\lambda, 
0)$ is nothing but $S,T \in \mathcal{T}^{+}_0(\lambda)$. For
$\lambda \in \Lambda^{+}$, let $\overline{\mathcal{S}_{0}}^{\vee
  \lambda} = \overline{\mathcal{S}_{0}}^{\vee}(\Lambda)^{\lambda}$ be   
the $R$-submodule of 
$\overline{\mathcal{S}^{0}}(\Lambda)$ spanned by
$\overline{{\varphi}}_{ST}$ with $S, T \in
\mathcal{T}^{+}_0(\alpha)$ for various $\alpha \in \Lambda^{+}$
such that $\alpha \rhd \lambda$. We show the following.

\addtocounter{thm}{1}
\begin{thm}
The algebra $\overline{\mathcal{S}^{0}}(\Lambda)$ has a free basis
\begin{equation*}
\overline{\mathcal{C}^{0}}(\Lambda) = \{
\overline{{\varphi}}_{ST} \mid S,T \in
\mathcal{T}^{+}_0(\lambda), ~ \lambda \in \Lambda^{+} \}
\end{equation*}
satisfying the following properties.

$({\rm{i}})$ The $R$-linear map $\ast :
\overline{\mathcal{S}^{0}}(\Lambda) \rightarrow
\overline{\mathcal{S}^{0}}(\Lambda)$ determined by
$\overline{{\varphi}}^{\ast}_{ST} =
\overline{{\varphi}}_{TS}$, for all $S,T \in
\mathcal{T}^{+}_0(\lambda)$ and all $\lambda \in \Lambda^{+}$,
is an anti-automorphism of $\overline{\mathcal{S}^{0}}(\Lambda)$. 

$({\rm{ii}})$ Let $T \in \mathcal{T}^{+}_0(\lambda)$. Then for all
$\overline{{\varphi}} \in
\overline{\mathcal{S}^{0}}(\Lambda)$, and any $V \in
\mathcal{T}^{+}_0(\lambda)$, there exists $r_{V} \in 
R$ such that      
\begin{equation*}
\overline{{\varphi}}_{ST} \cdot
\overline{{\varphi}} \equiv \underset{V \in
  \mathcal{T}^{+}_0(\lambda)}{\sum} r_{V}
\overline{{\varphi}}_{SV} \mod \mathcal{S}_{\overline0}^{\vee
  \lambda}
\end{equation*}
for any $S \in \mathcal{T}^{+}_0(\lambda)$, where $r_{V}$ is
independent of the choice of $T$.

In particular, $\overline{\mathcal{C}^{0}}(\Lambda)$ is a cellular  
basis of $\overline{\mathcal{S}^{0}}(\Lambda)$.
\end{thm}

\begin{proof}
To show (i) we first take a small detour. Let $f :
\mathcal{S}^{0}(\Lambda) \rightarrow
\overline{\mathcal{S}^{0}}(\Lambda)$ be the $R$-linear 
map given by 
\begin{equation*}
f({\varphi}_{ST})=
\left\{\begin{array}{ll} 
\overline{{\varphi}}_{TS} & \text{ if }
{\varphi}_{ST} \in \mathcal{C}^0(\lambda, 0) 
\text{ for some } \lambda \in \Lambda^{+},\\
0 & \text{ otherwise. } \\
\end{array}\right.
\end{equation*} 
We claim that $f$ is a surjective algebra
anti-homomorphism. It is clear from the definition that $f$ is  
surjective.
Then, to prove the claim it suffices to show the equations
\begin{equation}\label{eq:(i)}
f({\varphi}_{S_1T_1} {\varphi}_{S_2T_2}) 
= f({\varphi}_{S_2T_2}) \cdot
f({\varphi}_{S_1T_1})  
\end{equation}
for all $S_i, T_i \in \mathcal{T}^{+}_0(\lambda_i)$, $\lambda_i \in 
\Lambda^{+}$ ($i=1,2$) and 
\begin{equation}\label{eq:(ii)}    
f(1_{\mathcal{S}^{0}(\Lambda)}) 
= 1_{\overline{\mathcal{S}^{0}}(\Lambda)}.
\end{equation}
We show \eqref{eq:(ii)}. By
applying $f$ to 
\eqref{id:identity element of S^0}, we have
\begin{equation*}
f( 1_{\mathcal{S}^{0}(\Lambda)} ) 
= \sum_{\lambda \in \Lambda^{+}}
\sum_{S, T \in \mathcal{T}^{+}_0(\lambda)}
r_{T, S} \cdot \overline{{\varphi}_{ST}} 
= \overline{ (
1_{\mathcal{S}^{0}(\Lambda)} )^{*} } 
= \overline{1}_{\mathcal{S}^{0}(\Lambda)}
= 1_{\overline{\mathcal{S}^{0}}(\Lambda)}  
\end{equation*}
where $*$ is the involution of $\mathcal{S}(\Lambda)$. Hence we obtain      
\eqref{eq:(ii)}. We show \eqref{eq:(i)}. There are three cases to
consider. First, if ${\varphi}_{S_1T_1} \in
\mathscr{C}^0(\lambda, 1)$ with some $\lambda \in \Lambda^{+}$
then ${\varphi}_{S_1T_1} \cdot {\varphi}_{S_2T_2} \in
\mathcal{S}^{00}(\Lambda)$ by Lemma \ref{lem:sharper multiplication
  formulae} and hence $f({\varphi}_{S_1T_1}
{\varphi}_{S_2T_2}) = 0 =
f({\varphi}_{S_1T_1})$, so we are done in this
case. Similarly, it holds if ${\varphi}_{S_2T_2} \in
\mathcal{C}^0(\lambda, 1)$ 
with some 
$\lambda \in \Lambda^{+}$. So, we only need to check the case where
${\varphi}_{S_iT_i} \in \mathcal{C}^0(\lambda_i, 0)$ with
some $\lambda_i \in \Lambda^{+} ~ (i=1,2)$. Let $\mu_i, \nu_i
\in \Lambda$ with
$\alpha(\lambda_i) = \alpha(\mu_i) = \alpha(\nu_i)$ and $S_i \in
\mathcal{T}_0(\lambda_i,\mu_i), T_i \in
\mathcal{T}_0(\lambda_i,\nu_i) ~ (i=1,2)$. If $\nu_1 \neq \mu_2$ then 
${\varphi}_{S_1 T_1}
  \cdot {\varphi}_{S_2 T_2} = 0 = {\varphi}_{T_2S_2} \cdot
  {\varphi}_{T_1S_1}$. Hence we have
\begin{equation*}
f({\varphi}_{S_1T_1} \cdot {\varphi}_{S_2T_2})
= 0 = \overline{{\varphi}_{T_2S_2} \cdot
  {\varphi}_{T_1S_1}} =
\overline{{\varphi}}_{T_2S_2} \cdot
\overline{{\varphi}}_{T_1S_1} =
f({\varphi}_{S_2T_2}) \cdot f({\varphi}_{S_1T_1}).
\end{equation*} 
So we may assume that $\nu_1 = \mu_2$. It is
easy to see from our assumptions that $\bold{a}(\mu_{1})
=\bold{a}(\mu_{2})$. Hence, for $\lambda \in \Lambda^{+}$ and $S
\in \mathcal{T}_0(\lambda, \mu_1)$, $T \in \mathcal{T}_0(\lambda,
\nu_2)$, if ${\varphi}_{ST} \in \mathcal{C}^0(\lambda, 1)$
then ${\varphi}_{TS} \in
\mathcal{C}^0(\lambda, 1)$. Therefore, by the first formula in
Lemma \ref{lem:sharper multiplication formulae}, we have
\begin{equation}\label{eq:(1)}
{\varphi}_{S_1T_1} \cdot {\varphi}_{S_2T_2} =
\underset{\substack{
\lambda \unrhd \lambda_1, ~ S, T\\
{\varphi}_{ST} \in \mathcal{C}^0(\lambda, 0)}}{\sum} r_{ST}
\cdot 
\varphi_{ST} + \underset{\substack{
\lambda \unrhd \lambda_1, ~ S, T\\
{\varphi}_{ST} \in \mathcal{C}^0(\lambda, 1) }}{\sum} r_{ST} \cdot   
\varphi_{ST}
\end{equation}
where $ r_{ST} \in R$ and both sums are taken over $r$-partitions
$\lambda 
\unrhd \lambda_1$ and semistandard Tableaux $S \in
\mathcal{T}_0(\lambda, 
\mu_1)$, $T \in \mathcal{T}_0(\lambda,
\nu_2)$. This gives us two formulas. First, we have 
\begin{equation}\label{eq:(2)}
f({\varphi}_{S_1T_1} \cdot
{\varphi}_{S_2T_2}) = \underset{\substack{
\lambda \unrhd \lambda_1, ~ S, T ~ \\ 
{\varphi}_{ST} \in \mathcal{C}^0(\lambda,0)}}
{\sum} r_{ST} \cdot 
\overline{{\varphi}}_{TS}.
\end{equation} 
Secondly, by applying
the involution $\ast$ of $\mathcal{S}(\Lambda)$ on both
sides,
\begin{equation}\label{eq:(3)}
\overline{{\varphi}_{T_2S_2} \cdot
{\varphi}_{T_1S_1}} = \underset{\substack{
\lambda \unrhd \lambda_1, ~ S, T ~ \\ 
{\varphi}_{ST} \in \mathcal{C}^0(\lambda, 0)}}
{\sum} r_{ST} \cdot
\overline{{\varphi}}_{TS}
\end{equation} 
since in the second sum in \eqref{eq:(1)}, we have
 ${\varphi}_{ST}^{\ast} =  
 {\varphi}_{TS} \in \mathcal{C}^0(\lambda, 1)$. Comparing
 \eqref{eq:(2)} and \eqref{eq:(3)}, we have
\begin{equation*}
f({\varphi}_{S_1T_1}
 \cdot  
{\varphi}_{S_2T_2}) = \overline{{\varphi}_{T_2S_2} \cdot
{\varphi}_{T_1S_1}} =
\overline{{\varphi}}_{T_2S_2} \cdot
\overline{{\varphi}}_{T_1S_1} =
f({\varphi}_{S_2T_2}) \cdot
f({\varphi}_{S_1T_1}).
\end{equation*}
This proves our claim. By definition of
$f$, we have $\ker f = \mathcal{S}^{00}(\Lambda)$. Hence $f$ induces
an 
algebra anti-automorphism on $\mathcal{S}^{\overline0}(\Lambda) =
\mathcal{S}^{0}(\Lambda) / \mathcal{S}^{00}(\Lambda)$, which maps
$\overline{{\varphi}}_{ST}$ to
$\overline{{\varphi}}_{TS}$ for any $S,T \in
\mathcal{T}^{+}_0(\lambda)$ with $\lambda \in \Lambda^{+}$. Thus
(i) holds. Then (ii) follows from the second formula of
\eqref{eq:multiplication formulae of S^0(m,n)} by applying the natural 
map $\mathcal{S}^{0}(\Lambda) \rightarrow
\overline{\mathcal{S}^{0}}(\Lambda)$. 
\end{proof}

In the case where $\mathcal{S}^{\flat}(\bs{m},n)$ is defined,
$\overline{\mathcal{S}^{0}}(\Lambda)$ can be identified with
$\mathcal{S}^{\flat}(\bs{m},n)$, i.e, we have the following
proposition. 

\begin{prop}\label{prop:the relationship between S^overline0 and S^0}  
Let $\Lambda = \widetilde{\mathcal{P}}_{n,r}$ and assume that
\eqref{eq:The condition of the Modified Ariki-Koike} and
\eqref{con:m_i>=n} holds. Then there exists an algebra isomorphism
${\flat}: 
\overline{\mathcal{S}^{0}}(\Lambda) \rightarrow
\mathcal{S}^{\flat}(\bs{m},n)$ 
satisfying the following. For $\overline{{\varphi}}_{ST}
\in \overline{\mathcal{C}^{0}}(\Lambda)$ such that $S, T \in
\mathcal{T}^{+}_0(\lambda)$ and $\lambda \in \Lambda^{+}$, we have $(    
\overline{{\varphi}}_{ST} )^{\flat} = {\varphi}^{\flat}_{ST}$.
\end{prop}

\begin{proof}
By assumption, $\mathcal{S}^{\flat}(\bs{m},n)$ is defined. Now
\cite[Proposition 6.4 (ii)]{SawS} says that there exists a surjective
algebra  
homomorphism $\widehat{f} : \mathcal{S}^{0}(\Lambda) \rightarrow
\mathcal{S}^{\flat}(\bs{m},n)$ such that
\begin{equation*}
\widehat{f}({\varphi}_{ST}) = 
\left
\{
\begin{array}{ll}
{\varphi}^{\flat}_{ST} 
&
\text{ if } S,T \in \mathcal{T}^{+}_0(\lambda), \\
0
&
\text{ otherwise } \\
\end{array}\right.
\end{equation*}
for ${\varphi}_{ST} \in \mathcal{C}^{0}(\Lambda)$. Then $\ker
\widehat{f} = 
\mathcal{S}^{00}(\Lambda)$ and so $\widehat{f}$ induces an algebra
isomorphism $\overline{\mathcal{S}^{0}}(\Lambda) \rightarrow
\mathcal{S}^{\flat}(\bs{m},n)$, which maps
$\overline{{\varphi}}_{ST}$ to ${\varphi}^{\flat}_{ST}$ for any 
$S,T \in \mathcal{T}^{+}_0(\lambda)$ and $\lambda \in
\Lambda^{+} = \mathcal{P}_{n,r}$. We obtain the result.
\end{proof}

We now return to the general setting, and consider
$\overline{\mathcal{S}^{0}}(\Lambda)$ for arbitrary $\Lambda$. The
above proposition says that the $\overline{\mathcal{S}^{0}}(\Lambda)$
is     
a natural "cover" of the $\mathcal{S}^{\flat}(\bs{m},n)$.

For $\lambda \in \Lambda^{+}$,
$\overline{{\varphi}}_{\lambda} =
\overline{{\varphi}}_{T^{\lambda} T^{\lambda}}$ is an
element in $\overline{\mathcal{S}^{0}}(\Lambda)$. Hence, by the
cellular theory \cite{GL}, one can define a Weyl module
$\overline{Z}^{\lambda}$ of $\overline{\mathcal{S}^{0}}(\Lambda)$ as 
the right $\overline{\mathcal{S}^{0}}(\Lambda)$-submodule of
$\overline{\mathcal{S}^{0}}(\Lambda) /
\overline{\mathcal{S}_{0}}^{\vee 
  \lambda}$ spanned by the image of
$\overline{{\varphi}}_{\lambda}$. We denote by
$\overline{{\varphi}}_{T}$ the image of
$\overline{{\varphi}}_{T^{\lambda} T}$ in
$\overline{\mathcal{S}^{0}}(\Lambda) /
\overline{\mathcal{S}_{0}}^{\vee 
  \lambda}$. Then the set $\{
\overline{{\varphi}}_{T} \mid T \in
\mathcal{T}^{+}_0(\lambda) \}$ is a free $R$-basis of
$\overline{Z}^{\lambda}$. Define a bilinear form $\langle ~ , ~
\rangle_{\overline{0}}$ on $\overline{Z}^{\lambda}$ by requiring that  
\begin{equation*}
\overline{{\varphi}}_{T^{\lambda}
  S}\overline{{\varphi}}_{T T^{\lambda}} \equiv \langle
  \overline{{\varphi}}_{S},
  \overline{{\varphi}}_{T} \rangle_{\overline{0}}
  \cdot \overline{{\varphi}}_{\lambda} \mod
  \overline{\mathcal{S}_{0}}^{\vee \lambda}
\end{equation*} 
for all $S,T \in \mathcal{T}^{+}_0(\lambda)$. Let
$\overline{L}^{\lambda} = \overline{Z}^{\lambda} / \text{rad} 
\overline{Z}^{\lambda}$, where $\text{rad}
\overline{Z}^{\lambda} = \{x\in \overline{Z}^{\lambda} \mid \langle x,
y 
\rangle_{\overline{0}} = 0
\text{ for all } y \in \overline{Z}^{\lambda} \}$. In the case where  
$R$ is a field, by a general theory
of cellular algebras, the set $\{ 
\overline{L}^{\lambda} \mid \lambda \in \Lambda^{+},
~ \overline{L}^{\lambda} \neq 0 \}$ gives a complete set of
non-isomorphic irreducible
$\overline{\mathcal{S}^{0}}(\Lambda)$-modules. Furthermore, we have
the following result.

\begin{prop}\label{prop:complete set of non-isomorphic irreducible
    S^overline0(m,n) modules}
Suppose that $R$ is a field. Then $\overline{L}^{\lambda} \neq 0$ for 
any $\lambda \in \Lambda^{+}$. Hence, $\{
\overline{L}^{\lambda} \mid \lambda \in \Lambda^{+} \}$ is a
complete set of non-isomorphic irreducible
$\overline{\mathcal{S}^{0}}(\Lambda)$-modules. Therefore,
$\overline{\mathcal{S}^{0}}(\Lambda)$ is quasi-hereditary. 
\end{prop}

\begin{proof}
We have $\overline{{\varphi}}_{T^{\lambda}
  T^{\lambda}}\overline{{\varphi}}_{T^{\lambda} T^{\lambda}} \equiv 
  \langle 
  \overline{{\varphi}}_{T^{\lambda}},
  \overline{{\varphi}}_{T^{\lambda}} \rangle_{\overline{0}} 
  \cdot \overline{{\varphi}}_{\lambda} \mod
  \overline{\mathcal{S}_{0}}^{\vee \lambda}$. But since
  $\overline{{\varphi}}_{T^{\lambda} 
  T^{\lambda}}\overline{{\varphi}}_{T^{\lambda}
  T^{\lambda}} = \overline{{\varphi}}_{T^{\lambda}
  T^{\lambda}}$, we see that $\langle
  \overline{{\varphi}}_{T^{\lambda}}, 
  \overline{{\varphi}}_{T^{\lambda}} \rangle_{\overline{0}}
  = 1$, and so $\overline{L}^{\lambda}$ is non-zero. In particular,
  $\overline{\mathcal{S}^{0}}(\Lambda)$ is a standardly full-based
  algebra (see Definition in \ref{def:the definition of
  standardly (full-)based algebra}). Hence Theorem
  \ref{th:quasi-hereditary} gives us that
  $\overline{\mathcal{S}^{0}}(\Lambda)$ is quasi-hereditary. The
  proposition is proved. 
\end{proof}

\begin{lem}\label{lem:Z(lambda,0) is isomorphic to overlineZlambda}
Suppose that $\lambda \in \Lambda^{+}$. Then there exists an
$\mathcal{S}^{0}(\Lambda)$-module isomorphism $h_{\lambda} :
Z^{(\lambda,0)} \rightarrow \overline{Z}^{\lambda}$ such that
$h_{\lambda}({\varphi}^{0}_{T}) =
\overline{{\varphi}}_{T}$ for any $T \in
\mathcal{T}^{+}_0(\lambda)$, where $\overline{Z}^{\lambda}$ is
regarded as an $\mathcal{S}^{0}(\Lambda)$-module via the canonical map  
$\pi : \mathcal{S}^{0}(\Lambda)
\rightarrow \overline{\mathcal{S}^{0}}(\Lambda)$.
\end{lem}

\begin{proof}
Let $\pi : \mathcal{S}^{0}(\Lambda) \rightarrow
\overline{\mathcal{S}^{0}}(\Lambda)$ be the natural surjection. Then  
$\pi$ maps $\mathcal{S}_{0}^{\vee \lambda}$ to
$\overline{\mathcal{S}_{0}}^{\vee 
  \lambda}$, hence it induces an $\mathcal{S}^{0}(\Lambda)$-module
surjective 
homomorphism $\overline{\pi} : \mathcal{S}^{0}(\Lambda) /
\mathcal{S}_{0}^{\vee \lambda} 
\rightarrow \overline{\mathcal{S}_{0}}(\Lambda) /
\overline{\mathcal{S}_{0}}^{\vee  
  \lambda}$. Here $Z^{(\lambda, 0)} = {\varphi}^{0}_{T^{\lambda}}
\mathcal{S}^{0}(\Lambda)$ and $\overline{Z}^{\lambda} =
\overline{{\varphi}}_{\lambda} \overline{\mathcal{S}_{0}}(\Lambda) =
\overline{{\varphi}}_{\lambda} \mathcal{S}_{0}(\Lambda)$, and
$\overline{\pi}( {\varphi}^{0}_{T^{\lambda}} ) =
\overline{{\varphi}}_{\lambda}$. Hence $\overline{\pi}$ induces a
surjective map $h_{\lambda} : Z^{(\lambda, 0)} \rightarrow
\overline{Z}^{\lambda}$. Since it is easy to check that
$\overline{\pi}( {\varphi}^{0}_{T} ) = \overline{{\varphi}}_{T}$ for
any $T \in \mathcal{T}^{+}_0(\lambda)$, the map $h_{\lambda}$ is an
isomorphism.   
\end{proof}

Combining Lemma \ref{lem:Z(lambda,0) is isomorphic to overlineZlambda} 
  and Lemma \ref{lem:injective homomorphism
  f_lambda from
  Z^(lambda,0) to W^lambda_natural}, we have the following.

\begin{lem}
We regard $\overline{Z}^{\lambda}$ and $W^{\lambda}$ as
right $\mathcal{S}^{0}(\Lambda)$-modules via the above and via the 
restriction, respectively. Then the map $\overline{f}_{\lambda} =
f_{\lambda} \circ h_{\lambda}^{-1} : \overline{Z}^{\lambda}
\rightarrow W^{\lambda}$, which maps
$\overline{{\varphi}}_{T}$ to ${\varphi}_{T}$
for any $T \in \mathcal{T}^{+}_0(\lambda)$, is an injective
$\mathcal{S}^{0}(\Lambda)$-module homomorphism. 
\end{lem}

We need two easy lemmas.

\begin{lem}\label{lem:irreducible S^0(m,n)-modules for composition
    factors of Z^(lambda,0)}
Suppose that $R$ is a field. For every $\lambda \in
\Lambda^{+}$, we regard $\overline{L}^{\lambda}$ as right
$\mathcal{S}^{0}(\Lambda)$-modules via the map $\pi$. Then
$\mathcal{S}^{0}(\Lambda)$-module 
$\overline{L}^{\lambda}$ is irreducible. Furthermore, $\{
 \overline{L}^{\alpha} \mid \alpha \in \Lambda^{+}, ~ \lambda \unrhd 
 \alpha \}$ is a
complete set of pairwise inequivalent irreducible
$\mathcal{S}^{0}(\Lambda)$-modules occurring in the composition factors 
of $Z^{(\lambda, 0)}$.   
\end{lem}

\begin{proof}
By Proposition \ref{prop:complete set of non-isomorphic irreducible 
    S^overline0(m,n) modules}, we have $\overline{L}^{\lambda} \neq
    0$. Since $\overline{L}^{\lambda}$ is irreducible as
    $\overline{\mathcal{S}^{0}}(\Lambda)$-module, it is irreducible as 
    $\mathcal{S}^{0}(\Lambda)$-module. It is clear that $\{
    \overline{L}^{\alpha} \mid \alpha \in \Lambda^{+}, ~ \lambda
    \unrhd \alpha \}$ is a complete set of irreducible
    $\mathcal{S}^{0}(\Lambda)$-modules occurring in $Z^{(\lambda, 0)}$.   
\end{proof}

\begin{lem}\label{lem:L^lambda isomorphic to overlineL^lambda}
Suppose that $R$ is a field. Then, for all $\lambda \in
\Lambda^{+}$, we have $L_{0}^{\lambda} \simeq \overline{L}^{\lambda}$ 
as $\mathcal{S}^{0}(\Lambda)$-modules.
\end{lem}

\begin{proof}
(The following proof is almost similar to \cite[Lemma 6.8]{SawS}.)
Recall the bilinear form $\langle ~ , ~ \rangle_{0}$ (resp. $\langle ~  
, ~ \rangle_{\overline{0}}$) on $Z^{(\lambda, 0)}$ (resp. on
$\overline{Z}^{\lambda}$) and the $\mathcal{S}^{0}(\Lambda)$-module
isomorphism $h_{\lambda}$ in Lemma
\ref{lem:Z(lambda,0) is isomorphic to overlineZlambda}. We note that
$h_{\lambda}^{-1} : \overline{Z}^{\lambda} \rightarrow Z^{(\lambda,
  0)}$ maps $\text{rad}\overline{Z}^{\lambda}$ onto
$\text{rad}Z^{(\lambda, 0)}$. In fact, it is easy to see that $\langle
{\varphi}^{0}_{S}, {\varphi}^{0}_{T} \rangle_{0} =
\langle \overline{\varphi}_{S},
\overline{\varphi}_{T} \rangle_{\overline{0}}$ for any $S,
T \in \mathcal{T}^{+}_0(\lambda)$. Hence
$h_{\lambda}^{-1}(\text{rad}\overline{Z}^{\lambda}) =
\text{rad}Z^{(\lambda, 0)}$, and $h_{\lambda}^{-1}$ induces an
$\mathcal{S}^{0}(\Lambda)$-module isomorphism $\overline{L}^{\lambda}
\rightarrow L_{0}^{\lambda}$.
\end{proof}

 The following result connects the decomposition numbers in 
$\overline Z^{\lambda}$ and in $Z^{(\lambda,0)}$. 

\begin{thm}\label{thm:the property of L^lambda and
    [overlineZ^lambda:overlineL^mu]=[Z^(lambda,0):L_0^mu]}
Suppose that $R$ is a field.  Then
\begin{enumerate}
\item
 $\{ L_{0}^{\alpha} \mid \alpha \in
\Lambda^{+}, ~ \lambda \unrhd \alpha \}$ is a complete set of pairwise  
inequivalent 
irreducible $\mathcal{S}^{0}(\Lambda)$-modules occurring in the
composition 
factors of the 
$\mathcal{S}^0(\Lambda)$-module $Z^{(\lambda,0)}$.

\medskip

\item
 For $\lambda$, $\mu \in \Lambda^+$, we have 
\begin{equation*}
[\overline Z^{\lambda} : \overline L^{\mu}] = [ Z^{(\lambda, 0)} :   
L_{0}^{\mu} ].
\end{equation*}
\item
For $\lambda$, $\mu \in \Lambda^+$ such that 
$\alpha(\lambda) \ne \alpha(\mu)$, we have
\begin{equation*}
[\overline Z^{\lambda} : \overline L^{\mu}] = 0.
\end{equation*} 
\end{enumerate}
\end{thm}

\begin{proof}
The first and the second statement follows from Lemma
    \ref{lem:irreducible 
    S^0(m,n)-modules for composition
    factors of Z^(lambda,0)} and Lemma \ref{lem:L^lambda isomorphic to 
    overlineL^lambda}.
We show the third statement.

Take $\lambda \in \Lambda^{+}$ and assume that $\alpha(\lambda) \neq
\alpha(\mu)$ for $\mu \in \Lambda^{+}$. Let 
\begin{equation*}
0 = M_{1} \subsetneqq M_{2} \subsetneqq \cdots \subsetneqq M_{k} =
\overline{Z}^{\lambda} 
\end{equation*}
be a composition series of $\overline{Z}^{\lambda}$ and assume that
$M_{i+1} / M_{i} \simeq \overline{L}^{\mu}$ for some integer $i$ with
$1 \le i \le k-1$. Note that, as in the proof of Proposition
\ref{prop:complete set of non-isomorphic irreducible
    S^overline0(m,n) modules},
$\overline{\varphi}_{T^{\mu}} + \text{rad}\overline{Z}^{\mu}$ is a
non-zero element of $\overline{L}^{\mu}$. Then, there exists $c_{\mu}
\in M_{i+1} / M_{i}$ corresponding to $\overline{\varphi}_{T^{\mu}} +
\text{rad}\overline{Z}^{\mu}$. Take a representative
$\widetilde{c}_{\mu} \in M_{i+1} \subset \overline{Z}^{\lambda}$ of
$c_{\mu}$. Take ${\varphi}_{\mu} \in
\mathcal{S}^{0}(\Lambda)$ (cf. \eqref{eq:identity map on
  M^mu_natural}). Since $M_{i+1}$ is an
$\overline{\mathcal{S}^{0}}(\Lambda)$-module of
$\overline{Z}^{\lambda}$, we have $\widetilde{c}_{\mu} \cdot
\overline{{\varphi}_{\mu}}  
\in M_{i+1}$. Moreover, $\widetilde{c}_{\mu} \cdot
\overline{{\varphi}_{\mu}}$ 
corresponds to  
\begin{equation*}
( \overline{\varphi}_{T^{\mu}} + \text{rad}\overline{Z}^{\mu} ) \cdot
\overline{{\varphi}_{\mu}} = \overline{\varphi}_{T^{\mu}} +
\text{rad}\overline{Z}^{\mu} \neq 0
\end{equation*} 
Therefore, $\widetilde{c}_{\mu} \cdot \overline{{\varphi}_{\mu}} \neq
0$. On the 
other 
hand, for any $T \in \mathcal{T}^{+}_0(\lambda)$, the type of $T$ is
not equal to $\mu$ since 
$\alpha(\lambda) \neq
\alpha(\mu)$. Hence, by the definition
of ${\varphi}_{\mu}$, we have
${\varphi}_{T^{\lambda} T} \cdot {\varphi}_{\mu} = 0$
for any $T \in \mathcal{T}^{+}_0(\lambda)$. This implies that $x \cdot  
\overline{{\varphi}_{\mu}} = 0$ for any $x \in
\overline{Z}^{\lambda}$, a contradiction. The theorem follows.  
\end{proof}

\medskip

%%%%%%%%%%%%%%%%%%%%%%%%%%%%%%%%%%%%%%%%%%%%%%%%%%%%%%%%%%%%%%%%%%%%%%%%
%%%%%%%%%%%%%%%%%%%%%%%%%%%%%%%%%%%%%%%%%%%%%%%%%%%%%%%%%%%%%%%%%%%%%%%% 

\section{An estimate for decomposition numbers}

We are now ready to estimate the decomposition numbers for the
cyclotomic $q$-Schur algebras.

\subsection{}
We keep the notation in Section \ref{sec:A relationship between
  S^flat(m,n) and S^0(Lambda)}, and consider the general $\Lambda$. By
  abuse of the notation, $\mathcal{S}^{0},
\overline{\mathcal{S}^{0}}, \mathcal{S}$ denotes
$\mathcal{S}^{0}(\Lambda),
\overline{\mathcal{S}^{0}}(\Lambda), \mathcal{S}(\Lambda)$
respectively, and $M \otimes \mathcal{S}$ denotes $M
\otimes_{\mathcal{S}^{0}} \mathcal{S}$ for any right
$\mathcal{S}^{0}$-module $M$. Suppose that $\lambda \in
\Lambda^{+}$ and $\psi \in \mathcal{S}^{0},
{\varphi} \in \mathcal{S}$. Let $g_{\lambda} =
g : Z^{(\lambda,0)} \otimes \mathcal{S} \rightarrow
W^{\lambda}$ be the isomorphism given in Theorem \ref{thm:tensor
  theorem}. Then for any $T \in
\mathcal{T}^{+}_0(\lambda)$, we have 
\begin{equation*}
{\varphi}^{0}_{T} \otimes {\varphi} =
{\varphi}^{0}_{T^{\lambda}} \cdot
{\varphi}_{T^{\lambda} T} \otimes {\varphi} =
{\varphi}^{0}_{T^{\lambda}} \otimes
{\varphi}_{T^{\lambda} T} {\varphi}
\overset{g}{\longmapsto}
{\varphi}_{T^{\lambda}}{\varphi}_{T^{\lambda} T}
{\varphi} = {\varphi}_{T} {\varphi}
\end{equation*} 
where the second equality follows from ${\varphi}_{T^{\lambda} T} \in
\mathcal{S}^{0}$ since $T
\in \mathcal{T}^{+}_0(\lambda)$. Thus, we see that

\addtocounter{thm}{1}
\begin{rem}\label{rem:Remark of the tensor theorem}
The $\mathcal{S}$-module isomorphism
\begin{equation*}
g_{\lambda} : Z^{(\lambda,0)} \otimes \mathcal{S}
\rightarrow W^{\lambda}
\end{equation*}
in Theorem \ref{thm:tensor theorem} maps ${\varphi}^{0}_{T}
\otimes {\varphi}$ to 
${\varphi}_{T} {\varphi}$ for any $T \in
\mathcal{T}^{+}_0(\lambda)$.
\end{rem}

\begin{lem}\label{lem:the key Lemma1}
Suppose that $\lambda \in \Lambda$. Let $M_1, M_2$ be non-zero
$\mathcal{S}^{0}$-submodules of $Z^{(\lambda, 0)}$ with $M_1
\subsetneqq M_2$, and let $\iota : M_1 \hookrightarrow
M_2$,~~~ ${\iota}_1 : M_1 \hookrightarrow Z^{(\lambda, 0)}$,~~~
${\iota}_2 : M_2 
\hookrightarrow Z^{(\lambda, 0)}$ be inclusion maps. Then

$({\rm{i}})$ the $\mathcal{S}$-module homomorphism ~ ~ $\iota
\otimes 
{\rm{id}}_{\mathcal{S}} : M_1 \otimes
\mathcal{S} \rightarrow M_2 \otimes \mathcal{S}$ is a non-zero map.  

$({\rm{ii}})$ $(\iota_{1} \otimes {\rm{id}}_{\mathcal{S}})(M_1 \otimes   
\mathcal{S})$ is a proper $\mathcal{S}$-submodule of $(\iota_{2}
\otimes 
{\rm{id}}_{\mathcal{S}})(M_2 \otimes
\mathcal{S})$.
\end{lem}

\begin{proof}
Put $I = \iota \otimes
\text{id}_{\mathcal{S}}$, $I_i = \iota_{i} \otimes
\text{id}_{\mathcal{S}}$ $(i=1,2)$. Take $0 \neq x \in
M_1$. Then, by $M_1 \subset Z^{(\lambda, 0)}$, we may write $x =
\sum_{T \in \mathcal{T}^{+}_0(\lambda)} r_{T}
{\varphi}^{0}_{T}$ with $r_{T} \in R$. Then we have $(g \circ I ) ( x
\otimes 1) = \sum_{T \in
  \mathcal{T}^{+}_0(\lambda)} r_{T} {\varphi}_{T}$ by Remark
\ref{rem:Remark of the tensor theorem}. Now suppose that $I( x \otimes
1) = 0$. Then $( g_{\lambda} \circ I )( x \otimes 1 ) = 0$ and so
$r_{T} =0$ for any $T \in \mathcal{T}^{+}_0(\lambda)$ since $\{
{\varphi}_{T} \mid T \in \mathcal{T}^{+}_0(\lambda)\}$ is a
subset of the basis of $W^{\lambda}$. So $x
= 0$, which contradicts the choice of $x$. Hence $I$ is non-zero. This
proves (i). 

Next we show (ii). Clearly, $I_1(M_1 \otimes \mathcal{S}) \subset
I_2(M_2 \otimes \mathcal{S})$ so it suffices to see that $I_1(M_1
\otimes \mathcal{S}) \neq I_2(M_2 \otimes \mathcal{S})$. By our
assumption, 
there exists a non-zero element $x \in M_2 \setminus M_1$. Consider $x
\otimes 1 \in M_2 \otimes \mathcal{S}$ and, by way
of contradiction, suppose that $I_2(x \otimes 1) \in I_1(M_1
\otimes \mathcal{S})$. Hence there exists $y_i \in
M_1$ and ${\varphi}_{i} \in \mathcal{S}$ such
that 
\begin{equation}\label{eq:the image of x tensor 1}
(g_{\lambda} \circ I_2)(x \otimes 1) = (g_{\lambda} \circ
I_1)(\sum_{i} y_i \otimes {\varphi}_{i}).
\end{equation}
Since $x, y_i \in Z^{(\lambda, 0)}$, one can write $x = \sum_{T \in 
  \mathcal{T}^{+}_0(\lambda)} r_{T} {\varphi}^{0}_{T}$
and $y_i = \sum_{T \in \mathcal{T}^{+}_0(\lambda)} r_{T}^{i}
{\varphi}^{0}_{T}$ with $r_{T}, r_{T}^{i} \in
R$. By substituting them into \eqref{eq:the image of x tensor 1}, we
  have 
\begin{equation}\label{eq:substituting formula (1)}
\sum_{T \in
  \mathcal{T}^{+}_0(\lambda)} r_{T} {\varphi}_{T} =
\sum_{i} \sum_{T \in \mathcal{T}^{+}_0(\lambda)} r_{T}^{i}
{\varphi}_{T} {\varphi}_{i}.
\end{equation}
Hence, by using the definition of ${\varphi}_{T}$, we have
\begin{equation}\label{eq:substituting formula (2)}
\sum_{T \in
  \mathcal{T}^{+}_0(\lambda)} r_{T} {\varphi}_{T^{\lambda}
  T} - \sum_{i} \sum_{T \in \mathcal{T}^{+}_0(\lambda)} r_{T}^{i}
{\varphi}_{T^{\lambda} T} {\varphi}_{i} =
\psi
\end{equation}
 for some $\psi \in
\mathcal{S}^{\vee \lambda}$. Take $\nu \in
\Lambda$ with $\alpha(\nu) = \alpha(\lambda)$, and multiply
${\varphi}_{\nu}$ on the both side of \eqref{eq:substituting formula
  (2)} from the
right, we obtain
\begin{equation}\label{eq:multiply varphi_nu from the right}
\underset{T \in \mathcal{T}_0(\lambda, \nu)}{\sum} r_{T}
{\varphi}_{T^{\lambda} T} - \underset{i}{\sum} \underset{T \in
  \mathcal{T}^{+}_0(\lambda)}{\sum} r_{T}^{i} 
{\varphi}_{T^{\lambda} T} {\varphi}_{i}
{\varphi}_{\nu} = \psi {\varphi}_{\nu}.
\end{equation} 
(Recall that ${\varphi}_{\nu}$ is the identity map on
$M^{\nu}$ and zero on $M^{\kappa}$ with $\nu
\neq \kappa \in \Lambda^{+}$ see \ref{standardly based subsection}.)
Now, by using Lemma
\ref{lem:multiplication formula} and \eqref{eq:identity map on
  M^mu_natural}, we see that ${\varphi}_{T^{\lambda} T} {\varphi}_{i} 
{\varphi}_{\nu}$ is an $R$-linear combination of elements
of the form ${\varphi}_{UV}$ where $U \in
\mathcal{T}_0({\lambda}', \lambda)$ and $V \in
\mathcal{T}_0({\lambda}', \nu)$ 
with ${\lambda}' \in \Lambda^{+}$. Since $\alpha(\nu) =
\alpha(\lambda)$, this implies that ${\varphi}_{T^{\lambda} T}
{\varphi}_{i} 
{\varphi}_{\nu} \in \mathcal{S}^{0}$. In addition, clearly
${\varphi}_{T^{\lambda} T} \in \mathcal{S}^{0}$ for any $T \in
\mathcal{T}_0(\lambda, \nu)$. Hence, the left side of
\eqref{eq:multiply varphi_nu from the right} is contained in
$\mathcal{S}^{0}$. On the other hand, since
$\psi \in
\mathcal{S}^{\vee \lambda}$ and
$\mathcal{S}^{\vee \lambda}$ is a two-sided
ideal of $\mathcal{S}$, $\psi
  {\varphi}_{\nu}$ is an element in
  $\mathcal{S}^{\vee \lambda}$. Comparing both
  sides of \eqref{eq:multiply varphi_nu from the right}, we have $\psi  
  {\varphi}_{\nu} \in \mathcal{S}^{0} \cap
  \mathcal{S}^{\vee \lambda} =
  \mathcal{S}_{0}^{\vee (\lambda, 1)} \subset
  \mathcal{S}_{0}^{\vee (\lambda, 0)}$. Therefore, we have 
\begin{equation}\label{eq:multiply varphi_nu and working modulo
    S_0^vee(lambda,0)} 
\underset{T \in \mathcal{T}_0(\lambda, \nu)}{\sum} r_{T}
{\varphi}_{T^{\lambda} T} 
\equiv \underset{i}{\sum} \underset{T \in 
  \mathcal{T}^{+}_0(\lambda)}{\sum} r_{T}^{i} 
{\varphi}_{T^{\lambda} T} {\varphi}_{i}
{\varphi}_{\nu} \mod \mathcal{S}_{0}^{\vee (\lambda,0)} 
\end{equation}
for any $\nu \in \Lambda$ with $\alpha(\nu) =
\alpha(\lambda)$. We note that ${\varphi}_{T^{\lambda} T}
{\varphi}_{i} 
{\varphi}_{\nu} =
{\varphi}_{T^{\lambda} T^{\lambda}} \cdot
{\varphi}_{T^{\lambda} T} {\varphi}_{i}
{\varphi}_{\nu}$ with ${\varphi}_{T^{\lambda} T} {\varphi}_{i}
{\varphi}_{\nu} \in \mathcal{S}^{0}$. Then \eqref{eq:multiply
  varphi_nu and working modulo
    S_0^vee(lambda,0)} implies that 
\begin{equation}
\sum_{T \in
  \mathcal{T}_0(\lambda, \nu)} r_{T}
{\varphi}^{0}_{T} = \sum_{i} \sum_{T \in \mathcal{T}^{+}_0(\lambda)} 
r_{T}^{i} {\varphi}^{0}_{T^{\lambda}} \cdot
{\varphi}_{T^{\lambda} T} {\varphi}_{i}
{\varphi}_{\nu}.
\end{equation}
It follows that
\begin{equation}\label{eq:summing up formulas}
\underset{T \in \mathcal{T}_0^{+}(\lambda)}{\sum} r_{T} 
{\varphi}^{0}_{T} = 
\underset{
\substack{
\nu \in \Lambda\\
\alpha(\nu) = \alpha(\lambda)
}
}
{\sum}
\underset{T \in \mathcal{T}_0(\lambda, \nu)}{\sum} 
r_{T} {\varphi}^{0}_{T}
= 
\underset{
\substack{
\nu \in \Lambda\\
\alpha(\nu) = \alpha(\lambda)
}
}
{\sum}
\underset{i}{\sum} 
\underset{T \in \mathcal{T}^{+}_0(\lambda)}{\sum} 
r_{T}^{i} {\varphi}^{0}_{T^{\lambda}}
\cdot {\varphi}_{T^{\lambda} T} {\varphi}_{i}
{\varphi}_{\nu}.
\end{equation}
Put ${\varphi}_{\Sigma} = \underset{\mu \in
  \Lambda, ~~  \alpha(\mu) = \alpha(\lambda)}{\sum}
  {\varphi}_{\mu}$. Then ${\varphi}_{T^{\lambda}
  T} = {\varphi}_{T^{\lambda}
  T} {\varphi}_{\Sigma}$ for all $T \in
  \mathcal{T}^{+}_0(\lambda)$ and, by using Lemma
  \ref{lem:multiplication formula} again,
  ${\varphi}_{\Sigma} {\varphi}_{i}
  {\varphi}_{\nu}$ is an $R$-linear combination of elements
  of the form ${\varphi}_{UV}$ where $U \in
  \mathcal{T}_{0}({\lambda}', \mu)$ and $V \in
  \mathcal{T}_{0}({\lambda}', \nu)$ with $\alpha(\mu) =
  \alpha(\lambda) =\alpha(\nu)$ for various ${\lambda}' \in
  \Lambda^{+}$ and $\mu \in
  \Lambda$. Therefore
  ${\varphi}_{\Sigma} {\varphi}_{i}
  {\varphi}_{\nu} \in \mathcal{S}^{0}$. We also know that
  ${\varphi}_{T^{\lambda} T} \in \mathcal{S}^{0}$ for $T
  \in \mathcal{T}^{+}_0(\lambda)$. Combining these facts we can
  compute 
\begin{equation*}
{\varphi}^{0}_{T^{\lambda}}
\cdot {\varphi}_{T^{\lambda} T}
{\varphi}_{i}
{\varphi}_{\nu}
={\varphi}^{0}_{T^{\lambda}}
\cdot {\varphi}_{T^{\lambda} T}
{\varphi}_{\Sigma} {\varphi}_{i}
{\varphi}_{\nu}
=( {\varphi}^{0}_{T^{\lambda}}
\cdot {\varphi}_{T^{\lambda} T} )
\cdot {\varphi}_{\Sigma}
{\varphi}_{i}
{\varphi}_{\nu} 
={\varphi}^{0}_{T}
\cdot {\varphi}_{\Sigma}
{\varphi}_{i}
{\varphi}_{\nu}.
\end{equation*}
Substituting this into \eqref{eq:summing up formulas}, we have 
\begin{equation}
\begin{array}{l}
x 
= \underset{T \in \mathcal{T}^{+}_0(\lambda)}{\sum} r_{T}
{\varphi}^{0}_{T} = \underset{
\substack{
\nu \in \Lambda \\
\alpha(\nu) =\alpha(\lambda)
}
}{\sum}
\underset{i}{\sum}
\underset{T \in \mathcal{T}^{+}_0(\lambda)}{\sum}
r_{T}^{i} {\varphi}^{0}_{T}
\cdot {\varphi}_{\Sigma}
{\varphi}_{i}
{\varphi}_{\nu} 

= \underset{
\substack{
\nu \in \Lambda \\
\alpha(\nu) =\alpha(\lambda)
}
}{\sum}
\underset{i}{\sum}
y_{i} \cdot {\varphi}_{\Sigma}
{\varphi}_{i}
{\varphi}_{\nu}.
\end{array}
\end{equation}
Since ${\varphi}_{\Sigma}
{\varphi}_{i}
{\varphi}_{\nu} \in \mathcal{S}^{0}$ and $y_{i} \in M_1$ we have
$y_{i} \cdot {\varphi}_{\Sigma}
{\varphi}_{i}
{\varphi}_{\nu} \in M_1$. Hence $x \in M_1$, a
contradiction. Therefore $I_2( x {\otimes} 1 ) \not\in I_1( M_1
{\otimes} \mathcal{S} )$, and this proves (ii). The
Lemma is proved. 
\end{proof}

\begin{lem}\label{lem:the unique maximal proper right
    widetildeS-submodule of L^lambda tensor widetildeS}
Suppose that $R$ is a field. For every $\lambda \in
\Lambda^{+}$, let ~ ${\iota} :
{\rm{rad}}Z^{(\lambda, 0)} \rightarrow Z^{(\lambda, 0)}$ be the
inclusion map. Then there exists an $\mathcal{S}$-submodule
$N^{\lambda}$ of
$L_{0}^{\lambda} \otimes_{\mathcal{S}^{0}} \mathcal{S}$ such
that 
\begin{equation}\label{eq:unique maximal submodule of L_0^lambda
    otimes S}
( L_{0}^{\lambda} \otimes_{\mathcal{S}^{0}} \mathcal{S} ) /
N^{\lambda} \simeq L^{\lambda}
\end{equation}
as $\mathcal{S}$-modules. Moreover,
$N^{\lambda}$ is the unique maximal submodule of $L_{0}^{\lambda} 
\otimes_{\mathcal{S}^{0}} \mathcal{S}$.
\end{lem}

\begin{proof}
By the definition of $L_{0}^{\lambda}$, for each $\lambda \in
\Lambda^{+}$, there exists a short exact sequence of
$\mathcal{S}^{0}$-modules  

\begin{equation*}
\begin{CD}
0 @>>> \text{rad}Z^{(\lambda, 0)} @>>> Z^{(\lambda,
  0)} @>>> L_{0}^{\lambda} @>>> 0.
\end{CD}
\end{equation*}
This sequence induces an exact sequence of right
$\mathcal{S}$-module 
\begin{equation*}
\begin{CD}
\text{rad}Z^{(\lambda, 0)} \otimes \mathcal{S}
@>I>> Z^{(\lambda,
  0)} \otimes \mathcal{S} @>>> L_{0}^{\lambda} \otimes
\mathcal{S} @>>> 0,
\end{CD}
\end{equation*}
with $I = \iota \otimes \text{id}_{\mathcal{S}}$. Further, note that
$\text{rad}W^{\lambda}$ is the unique
maximal submodule of
$W^{\lambda}$ and $L^{\lambda} \neq 0$ for every
$\lambda \in \Lambda^{+}$. Hence under the isomorphism in
Theorem \ref{thm:tensor theorem}, there exists an
$\mathcal{S}$-submodule $M \simeq 
\text{rad}W^{\lambda}$ of $Z^{(\lambda, 0)} \otimes
\mathcal{S}$ such that $\text{Im}( I ) \subset M \subsetneqq
Z^{(\lambda, 
  0)} \otimes S$.
Then we have an isomorphism $f : ( Z^{(\lambda, 0)} \otimes S ) /  
\text{Im}( I ) \xrightarrow{\sim} L_{0}^{\lambda} \otimes
\mathcal{S}$. If we put $N^{\lambda}  
= f( M / \text{Im}( I ) )$, we obtain \eqref{eq:unique maximal
  submodule of L_0^lambda
    otimes S}. 
The uniqueness of
$N^{\lambda}$ follows from the uniqueness of
$\text{rad}W^{\lambda}$ and we are done.
\end{proof}

\begin{lem}\label{lem:the key Lemma2}
Suppose that $R$ is a field. For any $\lambda \in
\Lambda^{+}$, we regard $L^{\lambda}$ as an $\mathcal{S}^{0}$-module
via the restriction. Then
$\mathcal{S}^{0}$-module $L^{\lambda}$ contains
$L_{0}^{\lambda}$ as an irreducible submodule.
\end{lem}

\begin{proof}
This Lemma can be proved in a similar way as \cite[Lemma
6.8]{SawS}. However, this fact itself is important for later
discussions, we give here the proof. 

Recall the bilinear form $\langle ~ , ~ \rangle$ (resp. $\langle ~ , ~
\rangle_{0}$) on $W^{\lambda}$ (resp. on $Z^{(\lambda, 0)}$) and the
injective $\mathcal{S}^{0}$-module homomorphism $f_{\lambda} :
Z^{(\lambda, 0)} \rightarrow W^{\lambda}$ in Lemma
\ref{lem:injective homomorphism f_lambda from
  Z^(lambda,0) to W^lambda_natural}. We note that $f_{\lambda}$ maps
$\text{rad}Z^{(\lambda, 0)}$ into
$\text{rad}W^{\lambda}$. In fact, by \eqref{eq:the
  relationship between <,>_0 to <,>_{natural}}, $\langle
{\varphi}_{S}, {\varphi}_{T} \rangle
= \langle
{\varphi}^{0}_{S}, {\varphi}^{0}_{T} \rangle_{0}$
for any $S, T \in \mathcal{T}^{+}_0(\lambda)$. On the other hand, if
$S \in \mathcal{T}^{+}_0(\lambda, \mu), T \in \mathcal{T}_0(\lambda,
\nu) \setminus \mathcal{T}^{+}_0(\lambda, \nu)$, then $\alpha(\mu)
\neq \alpha(\nu)$, and so $\mu \neq \nu$. It follows that
${\varphi}_{T^{\lambda} S} {\varphi}_{T
  T^{\lambda}} = 0$, and $\langle {\varphi}_{S},
{\varphi}_{T} \rangle = 0$. Hence
$f_{\lambda}(\text{rad}Z^{(\lambda, 0)}) \subset
\text{rad}W^{\lambda}$. Now $f_{\lambda}$ induces an
$\mathcal{S}^{0}$-module homomorphism $L_{0}^{\lambda} \rightarrow
L^{\lambda}$, which is clearly non-trivial. Since $L_{0}^{\lambda}$ 
is an irreducible $\mathcal{S}^{0}$-module, $L^{\lambda}$ contains an
irreducible module isomorphic to $L_{0}^{\lambda}$. 
\end{proof}

We can now state our main result.

\begin{thm}\label{thm:On a composition multiplicity for the cyclotomic
  q-Schur algebra 1}
Suppose that $R$ is a field. Then, for all $\lambda, \mu \in
\Lambda^{+}$, we have 
\begin{equation*}
[ \overline{Z}^{\lambda} : \overline{L}^{\mu} ] = [ Z^{(\lambda, 0)} :
L_{0}^{\mu} ] \le [ W^{\lambda} : L^{\mu} ].
\end{equation*}
\end{thm} 

\begin{proof}
Since the first equality was shown in Theorem \ref{thm:the property of
    L^lambda and
    [overlineZ^lambda:overlineL^mu]=[Z^(lambda,0):L_0^mu]}, it is 
    enough to see that $[
    Z^{(\lambda, 0)} : L_{0}^{\mu} ] \le [ W^{\lambda} :
    L^{\mu} ]$. Let
\begin{equation*}
0 = M_0 \varsubsetneqq M_1 \varsubsetneqq \cdots \varsubsetneqq M_{l-1}
\varsubsetneqq M_{l} = Z^{(\lambda, 0)} 
\end{equation*} 
be a composition series of $Z^{(\lambda, 0)}$, and for $k = 0, \ldots
, l-1$ assume that $M_{k+1} / M_{k} \simeq
L_{0}^{\mu_{k}}$ with some $\mu_{k} \in
\Lambda^{+}$. Let
${\iota}_{k} : M_{k} \hookrightarrow M_{k+1}$ be the inclusion map and
put $I_{k} = {\iota}_{k} \otimes \text{id}_{\mathcal{S}} :
M_{k} \otimes \mathcal{S} \rightarrow M_{k+1} \otimes
\mathcal{S}$. Note that
$I_{k}$ is an $\mathcal{S}$-module
homomorphism. Then we have an exact sequence $0 \rightarrow M_{k}
\rightarrow M_{k+1}  
\rightarrow L_{0}^{\mu_{k}} \rightarrow 0$ of
$\mathcal{S}^{0}$-modules. Then, we obtain an $\mathcal{S}$-module
exact sequence  
\begin{equation}\label{eq:S-module exact sequence}
\begin{CD}
M_{k} \otimes \mathcal{S} @>{I_{k}}>> M_{k+1}
\otimes \mathcal{S} @>>>
L_{0}^{\mu_{k}} \otimes \mathcal{S} @>>> 0. 
\end{CD}
\end{equation}
Let ${\iota}'_{k} : M_{k}
\hookrightarrow Z^{(\lambda, 0)}$ be the inclusion map and put $M_{k}
= ( {\iota}'_{k} \otimes \text{id}_{\mathcal{S}} )( M_{k} \otimes
\mathcal{S} 
)$. Then 
by Lemma \ref{lem:the
  key Lemma1} (i), (ii), we have a filtration of $\mathcal{S}$-module
in $W^{\lambda}$ 
\begin{equation*}
0 = \mathcal{M}_{0} \varsubsetneqq
\mathcal{M}_{1} \varsubsetneqq \mathcal{M}_{2}
\varsubsetneqq \cdots \varsubsetneqq \mathcal{M}_{l-1}
\varsubsetneqq \mathcal{M}_{l} = W^{\lambda}.
\end{equation*}
We shall compute $\mathcal{M}_{k+1} / \mathcal{M}_{k}$ for $0 \le 
k < l$. Let $f : M_{k+1} \otimes \mathcal{S}
\rightarrow ( {\iota}'_{k+1} \otimes
\text{id}_{\mathcal{S}} )( M_{k+1} \otimes
\mathcal{S} ) / ( {\iota}'_{k} \otimes
\text{id}_{\mathcal{S}} )( M_{k} \otimes
\mathcal{S} ) \simeq \mathcal{M}_{k+1} /
\mathcal{M}_{k}$ be the natural map. Clearly, $f$ is a surjective
$\mathcal{S}$-module homomorphism. Moreover, it is easy to see that
$\ker{f} \supset I_{k}(
M_{k} \otimes \mathcal{S} )$. Hence $f$ induces an
$\mathcal{S}$-module surjective map $\widetilde{f} : ( M_{k+1} \otimes  
\mathcal{S} ) / 
I_{k}( M_{k} \otimes \mathcal{S} ) \rightarrow
\mathcal{M}_{k+1} / \mathcal{M}_{k}$. Since
$\mathcal{M}_{k} \varsubsetneqq
\mathcal{M}_{k+1}$, $\widetilde{f}$ is a non-zero
map. Note that $( M_{k+1} \otimes \mathcal{S} ) / I_{k}(
M_{k} \otimes \mathcal{S} ) \simeq L_{0}^{{\mu}_{k}} \otimes
\mathcal{S}$ by \eqref{eq:S-module exact sequence}. Hence, there
exists a certain $\mathcal{S}$-submodule $\mathcal{N}_{k}$ of
$L_{0}^{{\mu}_{k}} \otimes \mathcal{S}$ such that $(
L_{0}^{{\mu}_{k}} \otimes \mathcal{S} ) / \mathcal{N}_{k}
\simeq \mathcal{M}_{k+1} / \mathcal{M}_{k}$. However, by Lemma 
\ref{lem:the unique maximal proper right widetildeS-submodule of
  L^lambda tensor widetildeS}, there exists a unique maximal submodule 
$N^{{\mu}_{k}}$ of
$L_{0}^{{\mu}_{k}} \otimes \mathcal{S}$ such that $(
L_{0}^{{\mu}_{k}} \otimes \mathcal{S} ) /
N^{{\mu}_{k}} \simeq L^{{\mu}_{k}}$. Since $N^{{\mu}_{k}} \supset
\mathcal{N}_{k}$, we have a surjective map
\begin{equation*}
\mathcal{M}_{k+1} / \mathcal{M}_{k} \simeq ( L_{0}^{{\mu}_{k}} \otimes 
\mathcal{S} ) / 
\mathcal{N}_{k} \rightarrow ( L_{0}^{{\mu}_{k}}
\otimes 
\mathcal{S} ) / N^{{\mu}_{k}} \simeq L^{{\mu}_{k}}.
\end{equation*} 
This proves that $[ Z^{(\lambda, 0)} : L_{0}^{\mu} ] \le [ W^{\lambda}
: L^{\mu} ]$.
\end{proof}

In the following special case, we have a more precise formula.

\begin{thm}\label{thm:On a composition multiplicity for the cyclotomic
  q-Schur algebra 2}
Suppose that $R$ is a field. Then, for all $\lambda, \mu \in
\Lambda^{+}$ with $\alpha( \lambda ) = \alpha( \mu )$,
\begin{equation*}
[ \overline{Z}^{\lambda} : \overline{L}^{\mu} ] = [ Z^{(\lambda, 0)} :
L_{0}^{\mu} ] = [ W^{\lambda} : L^{\mu} ].
\end{equation*}
\end{thm}

\begin{proof}
Assume that $\alpha(\lambda) = \alpha(\mu)$ for $\lambda$, $\mu \in
  \Lambda^{+}$. In view of Theorem \ref{thm:On a composition
  multiplicity for the cyclotomic
  q-Schur algebra 1}, it is enough to show that $[
    W^{\lambda} :
    L^{\mu} ] \le [ Z^{(\lambda, 0)} :
    L_{0}^{\mu} ]$. Take $\lambda \in \Lambda^{+}$ and let 
\begin{equation*}
0 = W_{0}
\varsubsetneqq W_{1} \varsubsetneqq \cdots \varsubsetneqq
W_{c} = W^{\lambda} 
\end{equation*}
 be a composition series of
$W^{\lambda}$ and assume that
$W_{i+1} / W_{i} \simeq
L^{{\mu}_{i}}$. We regard $W_{i}$, $L^{{\mu}_{j}}$ as
$\mathcal{S}^{0}$-modules, via the
restriction. By Lemma \ref{lem:the key Lemma2}, there exists a
submodule $W'_{i}$ of $W_{i+1}$ containing $W_{i}$ such that $W'_{i} /  
W_{i} \simeq L_{0}^{{\mu}_{i}}$. Moreover, by Proposition
\ref{prop:the 
  property of L^lambda and radZ^(lambda,0)}, $L_{0}^{\mu} \neq 0$ for 
all 
$\mu \in \Lambda^{+}$. Therefore, $W_{i}
\varsubsetneqq W'_{i} \subset W_{i+1}$. 
Recall the injective $\mathcal{S}^{0}$-module map $f_{\lambda} :
 Z^{(\lambda, 0)} 
 \rightarrow W^{\lambda}$ given in Lemma \ref{lem:injective
 homomorphism f_lambda from Z^(lambda,0) to W^lambda_natural} and put 
 $M_{i} = f_{\lambda}( Z^{(\lambda, 0)} ) \cap W_{i}$,
 $M'_{i} = f_{\lambda}( Z^{(\lambda, 0)} ) \cap W'_{i}$. Then we
 obtain a filtration of $\mathcal{S}^{0}$-modules 
\begin{equation*}
0 = M_{0} \subset M'_{0} \subset M_{1} \subset M'_{1} \subset \cdots 
\subset M_{c-1} \subset M'_{c-1} \subset M_{c} = f_{\lambda}(
Z^{(\lambda, 0)} ) 
\end{equation*} 
of $f_{\lambda}( Z^{(\lambda, 0)} )$. Since $f_{\lambda}$ is an
injection from $Z^{(\lambda, 0)}$ to $W^{\lambda}$, we can regard the
above filtration as
a filtration of $Z^{(\lambda, 0)}$. We claim that 
\begin{equation}\label{eq:M_i neq M'_i if alpha(mu_i)= alpha(lambda)}  
M_{i} \neq M'_{i} ~ \text{ if } ~ \alpha( \mu_{i} ) = \alpha(
\lambda).
\end{equation}
We show \eqref{eq:M_i neq M'_i if alpha(mu_i)= alpha(lambda)}. Assume 
that $\alpha( \mu_{i} ) = \alpha(
\lambda )$. Note that, as in the proof of Proposition \ref{prop:the
  property of L^lambda and radZ^(lambda,0)},
$\varphi^{0}_{T^{{\mu}_{i}}} + \text{rad}{Z^{({\mu}_{i}, 0)}}$ is a
non-zero element of $L_{0}^{{\mu}_{i}}$. Then, there exists
$x_{{\mu}_{i}} 
\in W'_{i} / W_{i}$ corresponding to $\varphi^{0}_{T^{{\mu}_{i}}} +
\text{rad}{Z^{({\mu}_{i}, 0)}}$. Take a representative
$\widetilde{x}_{{\mu}_{i}} \in W'_{i} \setminus W_{i}$ of
$x_{{\mu}_{i}}$. Take $\varphi_{{\mu}_{i}} \in \mathcal{S}^{0}$
(cf. \eqref{eq:identity map
  on M^mu_natural}). Since $\widetilde{x}_{{\mu}_{i}} \in
W'_{i}$ and $W'_{i}$ is an $\mathcal{S}^{0}$-submodule of
$W^{\lambda}$, we have $\widetilde{x}_{{\mu}_{i}} \varphi_{{\mu}_{i}}
\in 
W'_{i}$. Moreover, $\widetilde{x}_{{\mu}_{i}} \varphi_{{\mu}_{i}}$
corresponds to 
\begin{equation*}
(\varphi^{0}_{T^{{\mu}_{i}}} + \text{rad}{Z^{({\mu}_{i}, 0)}} ) \cdot 
\varphi_{{\mu}_{i}} = \varphi^{0}_{T^{{\mu}_{i}}} +
\text{rad}{Z^{({\mu}_{i}, 0)}} \neq 0
\end{equation*}
Therefore, $\widetilde{x}_{{\mu}_{i}}
 \varphi_{{\mu}_{i}} \in W'_{i} \setminus W_{i}$. 
On the other hand, since $\widetilde{x}_{{\mu}_{i}} \in
 W_{i+1} \subset W^{\lambda}$, we can write $\widetilde{x}_{{\mu}_{i}}   
 = \sum_{T \in \mathcal{T}_{0}(\lambda)} r_{T} \varphi_{T} ~ ( r_{T}
 \in R )$ and hence $\widetilde{x}_{{\mu}_{i}}
 \varphi_{{\mu}_{i}} = \sum_{T \in \mathcal{T}_{0}(\lambda,
 {\mu}_{i})} r_{T} \varphi_{T}$ since $\varphi_{{\mu}_{i}}$ is an
 identity map on $M^{{\mu}_{i}}$ and is zero otherwise. By noticing
 $\alpha({\mu}_{i}) = \alpha(\lambda)$, we define an element
 $y_{{\mu}_{i}} \in Z^{(\lambda, 0)}$ by $y_{{\mu}_{i}} = \sum_{T \in
 \mathcal{T}_{0}(\lambda, {\mu}_{i})}
 r_{T} \varphi^{0}_{T}$. Then we have $f_{\lambda}( y_{{\mu}_{i}} ) =
 \sum_{T \in 
 \mathcal{T}_{0}(\lambda, {\mu}_{i})}
 r_{T} \varphi_{T} = \widetilde{x}_{{\mu}_{i}}
 \varphi_{{\mu}_{i}}$. It follows that $\widetilde{x}_{{\mu}_{i}}
 \varphi_{{\mu}_{i}} \in \{ f_{\lambda}( Z^{(\lambda, 0)} ) \cap
 W'_{i} \} \setminus \{
f_{\lambda}( Z^{(\lambda, 0)} ) \cap W_{i} \}$. 
Therefore, $M_{i} \neq
 M'_{i}$ and \eqref{eq:M_i neq M'_i if alpha(mu_i)= alpha(lambda)}
 holds.   

By the claim \eqref{eq:M_i neq M'_i if alpha(mu_i)= alpha(lambda)},
the quotient $M'_{i} / M_{i}$ is a non-zero
$\mathcal{S}^{0}$-submodule of $L_{0}^{{\mu}_{i}}$. Hence $M'_{i} /
M_{i} \simeq L_{0}^{{\mu}_{i}}$. This proves that $[ W^{\lambda} :
L^{\mu} ] \le [ Z^{(\lambda, 0)} : L_{0}^{\mu} ]$ for $\lambda$, $\mu
\in \Lambda^{+}$ such that $\alpha(\lambda) = \alpha(\mu)$. The
theorem is proved.
\end{proof}

\addtocounter{subsection}{6}
\subsection{}
We return to the setting in \ref{sec:modified Ariki-Koike
  algebra}. Let $\Lambda = \widetilde{\mathcal{P}}_{n,r}$ under the
  condition \eqref{eq:The condition of the Modified Ariki-Koike} and
  \eqref{con:m_i>=n}. For an
$r$-partition $\lambda \in \mathcal{P}_{n,r}$, we denote by
$\mathcal{S}_{\flat}^{\vee \lambda}$ the $R$-submodule of
$\mathcal{S}^{\flat}(\bs{m}, n)$ spanned by $\varphi^{\flat}_{ST}$
such that $S, T \in \mathcal{T}^{+}_{0}(\alpha)$ with $\alpha \rhd
\lambda$.  Moreover, for an
$r$-partition $\lambda \in \mathcal{P}_{n,r}$, $T^{\lambda} \in
\mathcal{T}^{+}_{0}(\lambda, \lambda)$, and in fact $T^{\lambda}$ is
the unique semistandard $\lambda$-Tableau of type $\lambda$. Moreover,
$\mathfrak{t} = \mathfrak{t}^{\lambda}$ is the unique element in
$\text{Std}(\lambda)$ such that $\lambda( \mathfrak{t} ) =
T^{\lambda}$. Thus, $m^{\flat}_{T^{\lambda} T^{\lambda}} =
m^{\flat}_{\mathfrak{t}^{\lambda} \mathfrak{t}^{\lambda}} =
m^{\flat}_{\lambda}$, and $\varphi^{\flat}_{\lambda} =
\varphi^{\flat}_{T^{\lambda} T^{\lambda}}$ is the identity map on
$M_{\flat}^{\lambda}$. We define the Weyl
module $W_{\flat}^{\lambda}$ as the right $\mathcal{S}^{\flat}(\bs{m},
n)$-submodule of $\mathcal{S}^{\flat}(\bs{m},
n) / \mathcal{S}_{\flat}^{\vee \lambda}$ spanned by the image of
$\varphi^{\flat}_{\lambda}$. For each $T \in
\mathcal{T}^{+}_{0}(\lambda, \mu)$, we denote by $\varphi^{\flat}_{T}$
the image of $\varphi^{\flat}_{T^{\lambda} T}$ in
$\mathcal{S}^{\flat}(\bs{m}, 
n) / \mathcal{S}_{\flat}^{\vee \lambda}$. Then we know that the Weyl
module $W_{\flat}^{\lambda}$ is $R$-free with basis $\{
\varphi^{\flat}_{T} \mid T \in \mathcal{T}^{+}_{0}(\lambda) \}$.
The Weyl module $W_{\flat}^{\lambda}$ enjoys an associative symmetric
bilinear form, defined by the equation
\begin{equation*}
\varphi^{\flat}_{T^{\lambda} S} \varphi^{\flat}_{T T^{\lambda}} \equiv
\langle \varphi^{\flat}_{S}, \varphi^{\flat}_{T} \rangle_{\flat}
 \cdot \varphi^{\flat}_{\lambda} \mod \mathcal{S}_{\flat}^{\vee \lambda}
\end{equation*}
for all $S, T \in \mathcal{T}^{+}_{0}(\lambda)$. Let
$L_{\flat}^{\lambda} = W_{\flat} / \text{rad}{W_{\flat}^{\lambda}}$,
where $\text{rad}{W_{\flat}^{\lambda}} = \{ x \in W_{\flat}^{\lambda}
\mid \langle x, y \rangle_{\flat} = 0 \text{ for all } y \in
W_{\flat}^{\lambda} \}$. By \cite[Proposition 5.11]{SawS}, we
know that, for all $r$-partition $\lambda \in \mathcal{P}_{n,r}$,
$L_{\flat}^{\lambda}$ is an absolutely irreducible and $\{
L_{\flat}^{\lambda} \mid 
\lambda \in \mathcal{P}_{n,r} \}$ is a complete set of non-isomorphic
irreducible $\mathcal{S}^{\flat}(\bs{m}, n)$-modules. Furthermore, for 
$\lambda, 
\mu \in \mathcal{P}_{n,r}$, we denote by $[
  W_{\flat}^{\lambda} : L_{\flat}^{\mu} ]$ the composition
  multiplicity of $L_{\flat}^{\mu}$ in $W_{\flat}^{\lambda}$.

Note that the above definition of the Weyl module
  $W_{\flat}^{\lambda}$ coincides with the
definition of the Weyl module $\overline{Z}^{\lambda}$ when
  $\mathcal{S}^{\flat}(\bs{m}, n)$ is isomorphic to
  $\mathcal{S}^{\overline0}(\Lambda)$ under the isomorphism $\flat$ in
  Proposition \ref{prop:the relationship between S^overline0 and
  S^0}. Consequently, under the isomorphism $\flat$, we have $[
  W_{\flat}^{\lambda} : L_{\flat}^{\mu} ] = [
  \overline{Z}^{\lambda} : \overline{L}^{\mu} ]$ for every $\lambda,
  \mu \in \mathcal{P}_{n,r}$. 

On the other hand, note that in the case where $r = 1$, the notation
for 
$\mathcal{S}^{\flat}(\bs{m}, n)$ coincides with the standard notation
for $q$-Schur algebras discussed as in \cite[Chapter 4]{M1}. So,
we use freely such a notation. For $\lambda, \mu \in
\mathcal{P}_{n,r}$, we denote by $[ W^{\lambda^{( i )}} : L^{\mu^{( i
    )}} ] ~ ( 1 \le i \le r )$ is defined as the composition
multiplicity of $L^{\mu^{( i )}}$ in $W^{\mu^{( i )}}$ for $\lambda =
( \lambda^{( 1 )}, \ldots , \lambda^{( r )} )$ and $\mu = ( \mu^{( i
  )}, \ldots , \mu^{( r )} )$. 

\addtocounter{thm}{1}
\begin{prop}[{\cite[Proposition 5.14]{SawS}}]\label{prop:on
    decomposition number}
Let $\Lambda = \widetilde{\mathcal{P}}_{n,r}$. Suppose that $R$ is a
field, and that \eqref{eq:The condition of the Modified Ariki-Koike}
and \eqref{con:m_i>=n} are satisfied. Let $\lambda$, $\mu \in
\mathcal{P}_{n,r}$. Then under the isomorphism in Theorem \ref{weak
  separation condition}, we have 
\begin{equation*}
[ W^{\lambda} : L^{\mu} ] = 
\begin{cases}
\prod_{i = 1}^{r} [ W^{{\lambda}^{( i )}} : L^{{\mu}^{( i )}} ] &
\text{ if } \alpha(\lambda) = \alpha(\mu), \\
0 & \text{ otherwise}.
\end{cases}
\end{equation*} 
\end{prop}

Then, combining Theorem \ref{thm:On a composition multiplicity
  for the cyclotomic q-Schur algebra 2} and Proposition \ref{prop:the
  relationship 
  between S^overline0 and S^0} and Proposition \ref{prop:on
    decomposition number}, we have the following. 

\begin{cor}
Let $\Lambda = \widetilde{\mathcal{P}}_{n,r}$. Suppose that $R$ is a
field, and that \eqref{eq:The condition of the Modified Ariki-Koike}
and \eqref{con:m_i>=n} are satisfied. Then, for
all $\lambda, \mu \in
\mathcal{P}_{n,r}$ with $\alpha( \lambda ) = \alpha( \mu )$, we have  
\begin{equation*}
[ W^{\lambda} : L^{\mu} ] =
\overset{r}{\underset{i=1}{\prod}} [ W^{\lambda^{( i )}} : L^{\mu^{( i
    )}} ]. 
\end{equation*} 
\end{cor}

\bigskip

\providecommand{\bysame}{\leavevmode\hbox to3em{\hrulefill}\thinspace}

\end{document}